\input amstex
\documentstyle{amsppt}
\magnification=\magstep 1
\document
\vsize6.7in
\input xy
\xyoption{all}
\NoBlackBoxes

\chardef\oldatsign=\catcode`\@
\catcode`\@=11
\newif\ifdraftmode			
\global\draftmodefalse


%
\font@\twelverm=cmr12 
\font@\twelvei=cmmi12 \skewchar\twelvei='177 
\font@\twelvesy=cmsy10 scaled\magstep1 \skewchar\twelvesy='060 
\font@\twelveex=cmex10 scaled\magstep1 
\font@\twelvemsa=msam10 scaled\magstep1 
\font@\twelvemsb=msbm10 scaled\magstep1 
\font@\twelvebf=cmbx12 
\font@\twelvett=cmtt12 
\font@\twelvesl=cmsl12 
\font@\twelveit=cmti12 
\font@\twelvesmc=cmcsc10 scaled\magstep1 
%
%
\font@\ninerm=cmr9 
\font@\ninei=cmmi9 \skewchar\ninei='177 
\font@\ninesy=cmsy9 \skewchar\ninesy='60 
\font@\ninemsa=msam9
\font@\ninemsb=msbm9
\font@\ninebf=cmbx9
%
%
%
\font@\ttlrm=cmbx12 scaled \magstep2 
\font@\ttlsy=cmsy10 scaled \magstep3 
\font@\tensmc=cmcsc10 
%
%
\def\normaltype{
	\def\pointsize@{12}%
	\abovedisplayskip18\p@ plus5\p@ minus9\p@
	\belowdisplayskip18\p@ plus5\p@ minus9\p@
	\abovedisplayshortskip1\p@ plus3\p@
	\belowdisplayshortskip9\p@ plus3\p@ minus4\p@
	\textonlyfont@\rm\twelverm
	\textonlyfont@\it\twelveit
	\textonlyfont@\sl\twelvesl
	\textonlyfont@\bf\twelvebf
	\textonlyfont@\smc\twelvesmc
	\ifsyntax@
		\def\big##1{{\hbox{$\left##1\right.$}}}%
	\else
		\let\big\twelvebig@
 \textfont0=\twelverm \scriptfont0=\ninerm \scriptscriptfont0=\sevenrm
 \textfont1=\twelvei  \scriptfont1=\ninei  \scriptscriptfont1=\seveni
 \textfont2=\twelvesy \scriptfont2=\ninesy \scriptscriptfont2=\sevensy
 \textfont3=\twelveex \scriptfont3=\twelveex  \scriptscriptfont3=\twelveex
 \textfont\itfam=\twelveit \def\it{\fam\itfam\twelveit}%
 \textfont\slfam=\twelvesl \def\sl{\fam\slfam\twelvesl}%
 \textfont\bffam=\twelvebf \def\bf{\fam\bffam\twelvebf}%
 \scriptfont\bffam=\ninebf \scriptscriptfont\bffam=\sevenbf
 \textfont\ttfam=\twelvett \def\tt{\fam\ttfam\twelvett}%
 \textfont\msafam=\twelvemsa \scriptfont\msafam=\ninemsa
 \scriptscriptfont\msafam=\sevenmsa
 \textfont\msbfam=\twelvemsb \scriptfont\msbfam=\ninemsb
 \scriptscriptfont\msbfam=\sevenmsb
	\fi
 \normalbaselineskip=\twelvebaselineskip
 \setbox\strutbox=\hbox{\vrule height12\p@ depth6\p@
      width0\p@}%
 \normalbaselines\rm \ex@=.2326ex%
}
%
%
%
\def\smalltype{
	\def\pointsize@{10}%
	\abovedisplayskip12\p@ plus3\p@ minus9\p@
	\belowdisplayskip12\p@ plus3\p@ minus9\p@
	\abovedisplayshortskip\z@ plus3\p@
	\belowdisplayshortskip7\p@ plus3\p@ minus4\p@
	\textonlyfont@\rm\tenrm
	\textonlyfont@\it\tenit
	\textonlyfont@\sl\tensl
	\textonlyfont@\bf\tenbf
	\textonlyfont@\smc\tensmc
	\ifsyntax@
		\def\big##1{{\hbox{$\left##1\right.$}}}%
	\else
		\let\big\tenbig@
	\textfont0=\tenrm \scriptfont0=\sevenrm \scriptscriptfont0=\fiverm 
	\textfont1=\teni  \scriptfont1=\seveni  \scriptscriptfont1=\fivei
	\textfont2=\tensy \scriptfont2=\sevensy \scriptscriptfont2=\fivesy 
	\textfont3=\tenex \scriptfont3=\tenex \scriptscriptfont3=\tenex
	\textfont\itfam=\tenit \def\it{\fam\itfam\tenit}%
	\textfont\slfam=\tensl \def\sl{\fam\slfam\tensl}%
	\textfont\bffam=\tenbf \def\bf{\fam\bffam\tenbf}%
	\scriptfont\bffam=\sevenbf \scriptscriptfont\bffam=\fivebf
	\textfont\msafam=\tenmsa
	\scriptfont\msafam=\sevenmsa
	\scriptscriptfont\msafam=\fivemsa
	\textfont\msbfam=\tenmsb
	\scriptfont\msbfam=\sevenmsb
	\scriptscriptfont\msbfam=\fivemsb
		\textfont\ttfam=\tentt \def\tt{\fam\ttfam\tentt}%
	\fi
 \normalbaselineskip 14\p@
 \setbox\strutbox=\hbox{\vrule height10\p@ depth4\p@ width0\p@}%
 \normalbaselines\rm \ex@=.2326ex%
}

\def\titletype{
	\def\pointsize@{17}%
	\textonlyfont@\rm\ttlrm
	\ifsyntax@
		\def\big##1{{\hbox{$\left##1\right.$}}}%
	\else
		\let\big\twelvebig@
		\textfont0=\ttlrm \scriptfont0=\twelverm
		\scriptscriptfont0=\tenrm
		\textfont2=\ttlsy \scriptfont2=\twelvesy
		\scriptscriptfont2=\tensy
	\fi
	\normalbaselineskip 25\p@
	\setbox\strutbox=\hbox{\vrule height17\p@ depth8\p@ width0\p@}%
	\normalbaselines
	\rm
	\ex@=.2326ex%
}

\def\tenbig@#1{
	{%
		\hbox{%
			$%
			\left
			#1%
			\vbox to8.5\p@{}%
			\right.%
			\n@space
			$%
		}%
	}%
}

\def\twelvebig@#1{%
	{%
		\hbox{%
			$%
			\left
			#1%
			\vbox to10.2\p@{}
			\right.%
			\n@space
			$%
		}%
	}%
}

%
%
%
%
%
\newif\ifl@beloutopen
\newwrite\l@belout
\newread\l@belin

\global\let\currentfile=\jobname

\def\getfile#1{%
	\immediate\closeout\l@belout
	\global\l@beloutopenfalse
	\gdef\currentfile{#1}%
	\input #1%
	\par
	\newpage
}

\def\getxrefs#1{%
	\bgroup
		\def\gobble##1{}
		\edef\list@{#1,}%
		\def\gr@boff##1,##2\end{
			\openin\l@belin=##1.xref
			\ifeof\l@belin
			\else
				\closein\l@belin
				\input ##1.xref
			\fi
			\def\list@{##2}%
			\ifx\list@\empty
				\let\next=\gobble
			\else
				\let\next=\gr@boff
			\fi
			\expandafter\next\list@\end
		}%
		\expandafter\gr@boff\list@\end
	\egroup
}

\def\testdefined#1#2#3{%
	\expandafter\ifx
	\csname #1\endcsname
	\relax
	#3%
	\else #2\fi
}

\def\document{%
	\minaw@11.11128\ex@ 
	\def\alloclist@{\empty}%
	\def\fontlist@{\empty}%
	\openin\l@belin=\jobname.xref	
	\ifeof\l@belin\else
		\closein\l@belin
		\input \jobname.xref
	\fi
}

\def\getst@te#1#2{%
	\edef\st@te{\csname #1s!#2\endcsname}%
	\expandafter\ifx\st@te\relax
		\def\st@te{0}%
	\fi
}

\def\setst@te#1#2#3{%
	\expandafter
	\gdef\csname #1s!#2\endcsname{#3}%
}

\outer\def\setupautolabel#1#2{%
	\def\newcount@{\global\alloc@0\count\countdef\insc@unt}	
	\def\newtoks@{\global\alloc@5\toks\toksdef\@cclvi}
	\expandafter\newcount@\csname #1Number\endcsname
	\expandafter\global\csname #1Number\endcsname=1%
	\expandafter\newtoks@\csname #1l@bel\endcsname
	\expandafter\global\csname #1l@bel\endcsname={#2}%
}

\def\reflabel#1#2{%
	\testdefined{#1l@bel}
	{
		\getst@te{#1}{#2}%
		\ifcase\st@te
			???
			\message{Unresolved forward reference to
				label #2. Use another pass.}%
		\or	
			\setst@te{#1}{#2}2
			\csname #1l!#2\endcsname 
		\or	
			\csname #1l!#2\endcsname 
		\or	
			\csname #1l!#2\endcsname 
		\fi
	}{
		{\escapechar=-1 
		\errmessage{You haven't done a
			\string\\setupautolabel\space for type #1!}%
		}%
	}%
}

{\catcode`\{=12 \catcode`\}=12
	\catcode`\[=1 \catcode`\]=2
	\xdef\Lbrace[{]
	\xdef\Rbrace[}]%
]%

\def\setlabel#1#2{%
	\testdefined{#1l@bel}
	{
		\edef\templ@bel@{\expandafter\the
			\csname #1l@bel\endcsname}%
		\def\@rgtwo{#2}%
		\ifx\@rgtwo\empty
		\else
			\ifl@beloutopen\else
				\immediate\openout\l@belout=\currentfile.xref
				\global\l@beloutopentrue
			\fi
			\getst@te{#1}{#2}%
			\ifcase\st@te
			\or	
			\or	
				\edef\oldnumber@{\csname #1l!#2\endcsname}%
				\edef\newnumber@{\templ@bel@}%
				\ifx\newnumber@\oldnumber@
				\else
					\message{A forward reference to label 
						#2 has been resolved
						incorrectly.  Use another
						pass.}%
				\fi
			\or	
				\errmessage{Same label #2 used in two
					\string\setlabel s!}%
			\fi
			\expandafter\xdef\csname #1l!#2\endcsname
				{\templ@bel@}
			\setst@te{#1}{#2}3%
			\immediate\write\l@belout 
				{\string\expandafter\string\gdef
				\string\csname\space #1l!#2%
				\string\endcsname
				\Lbrace\templ@bel@\Rbrace
				}%
			\immediate\write\l@belout 
				{\string\expandafter\string\gdef
				\string\csname\space #1s!#2%
				\string\endcsname
				\Lbrace 1\Rbrace
				}%
		\fi
		\templ@bel@	
		\expandafter\ifx\envir@end\endref 
			\gdef\marginalhook@{\marginal{#2}}%
		\else
			\marginal{#2}
		\fi
		\expandafter\global\expandafter\advance	
			\csname #1Number\endcsname
			by 1 %
	}{
		{\escapechar=-1
		\errmessage{You haven't done a \string\\setupautolabel\space
			for type #1!}%
		}%
	}%
}


\newcount\SectionNumber
\setupautolabel{t}{\number\SectionNumber.\number\tNumber}
\setupautolabel{r}{\number\rNumber}
\setupautolabel{T}{\number\TNumber}

\define\rref{\reflabel{r}}
\define\tref{\reflabel{t}}

\define\tnum{\setlabel{t}}
\define\rnum{\setlabel{r}}

%
\def\strutdepth{\dp\strutbox}%
\def\strutheight{\ht\strutbox}%

\newif\iftagmode
\tagmodefalse

\let\old@tagform@=\tagform@
\def\tagform@{\tagmodetrue\old@tagform@}

\def\marginal#1{%
	\ifvmode
	\else
		\strut
	\fi
	\ifdraftmode
		\ifmmode
			\ifinner
				\let\Vorvadjust=\Vadjust
			\else
				\let\Vorvadjust=\vadjust
			\fi
		\else
			\let\Vorvadjust=\Vadjust
		\fi
		\iftagmode	
			\llap{%
				\smalltype
				\vtop to 0pt{%
					\pretolerance=2000
					\tolerance=5000
					\raggedright
					\hsize=.72in
					\parindent=0pt
					\strut
					#1%
					\vss
				}%
				\kern.08in
				\iftagsleft@
				\else
					\kern\hsize
				\fi
			}%
		\else
			\Vorvadjust{%
				\kern-\strutdepth 
				{%
					\smalltype
					\kern-\strutheight 
					\llap{%
						\vtop to 0pt{%
							\kern0pt
							\pretolerance=2000
							\tolerance=5000
							\raggedright
							\hsize=.5in
							\parindent=0pt
							\strut
							#1%
							\vss
						}%
						\kern.08in
					}%
					\kern\strutheight
				}%
				\kern\strutdepth
			}
		\fi
	\fi
}


\newbox\Vadjustbox

\def\Vadjust#1{
	\global\setbox\Vadjustbox=\vbox{#1}%
	\ifmmode
		\ifinner
			\innerVadjust
		\fi		
	\else
		\innerVadjust
	\fi
}

\def\innerVadjust{%
	\def\nexti{\aftergroup\innerVadjust}%
	\def\nextii{%
		\ifvmode
			\hrule height 0pt 
			\box\Vadjustbox
		\else
			\vadjust{\box\Vadjustbox}%
		\fi
	}%
	\ifinner
		\let\next=\nexti
	\else
		\let\next=\nextii
	\fi
	\next
}%

\global\let\marginalhook@\empty

\def\endref{%
\setbox\tw@\box\thr@@
\makerefbox?\thr@@{\endgraf\egroup}%
  \endref@
  \endgraf
  \endgroup
  \keyhook@
  \marginalhook@
  \global\let\keyhook@\empty 
  \global\let\marginalhook@\empty 
}

\catcode`\@=\oldatsign

\nologo

\newcount\zero \zero=0
\newcount\prelim \prelim=1
\newcount\MC \MC=1
\newcount\PC \PC=2
\newcount\DP \DP=3
\newcount\Mac \Mac=4
\newcount\Lpq \Lpq=2
\newcount\GEN \GEN=3
\newcount\TMT \TMT=4
\newcount\exmpls \exmpls=5
\newcount\MR \MR=6
\newcount\nonempty \nonempty=7

\vsize6.7in
\topmatter

\title Artinian  Gorenstein algebras with  linear resolutions\endtitle
  \leftheadtext{El Khoury and  Kustin}
\rightheadtext{Artinian  Gorenstein algebras with  linear resolutions}
\author Sabine El Khoury\footnote{Part of this work was done while the author was on a research leave at the University of South Carolina. \phantom{XXX} \phantom{XXX} \phantom{XXX} \phantom{XXX} \phantom{XXX} \phantom{XXX} \phantom{XXX} \phantom{XXX} \phantom{XXX} \phantom{XXX} \phantom{XXX}  \phantom{XXX} \phantom{XXXXXXXXX}} and  
Andrew R. Kustin\footnote{Supported in part by the National Security Agency Grant 
and by the Simons Foundation. \phantom{XXX}}\endauthor
\address
Mathematics Department,
American University of Beirut,
Riad el Solh 11-0236,
Beirut,
Lebanon
\endaddress
\email se24\@aub.edu.lb \endemail
 \address
Mathematics Department,
University of South Carolina,
Columbia, SC 29208\endaddress
\email kustin\@math.sc.edu \endemail

\keywords  Artinian rings, Buchsbaum-Eisenbud ideals, Build resolution directly from inverse system, Compressed algebras, Flat family of Gorenstein Algebras,  Gorenstein rings, Linear presentation, Linear resolution, Macaulay inverse system, Parameterization of Gorenstein ideals with linear resolutions, Parameterization of linearly presented Gorenstein Algebras,  Pfaffians, Resolutions  
\endkeywords
\subjclass \nofrills{2010 {\it Mathematics Subject Classification.}} 13H10, 13E10, 13D02, 13A02
\endsubjclass

\endtopmatter

\document
{\eightpoint \flushpar{\bf Abstract.} 
For each pair of positive integers $n,d$, we construct a complex $\widetilde{\Bbb G}'(n)$ of modules
over  the bi-graded polynomial ring $\widetilde{R}=\Bbb Z[x_1,\dots,x_d,\{t_M\}]$, where $M$ roams over all monomials of degree $2n-2$ in $\{x_1,\dots,x_d\}$.
The complex $\widetilde{\Bbb G}'(n)$ has the following universal property. Let $P$ be the polynomial ring $\pmb k[x_1,\dots,x_d]$, where $\pmb k$ is a field, and let $\Bbb I_n^{[d]}(\pmb k)$ be the set of  
homogeneous ideals    $I$ in $P$, which are generated by  generated by forms of degree $n$, and for which  $P/I$
is an Artinian Gorenstein algebra with a linear resolution. If $I$ is an ideal from  $\Bbb I_n^{[d]}(\pmb k)$, then there exists
a homomorphism $\widetilde{R}\to P$, so that $P\otimes_{\widetilde{R}} \widetilde{\Bbb G}'(n)$ is a minimal homogeneous
resolution of $P/I$ by free $P$-modules.

The construction of $\widetilde{\Bbb G}'(n)$ is equivariant and explicit. We give the differentials of
$\widetilde{\Bbb G}'(n)$ as well as the modules. On the other hand, the homology of $\widetilde{\Bbb G}'(n)$ is unknown as are the  the properties of the modules that comprise $\widetilde{\Bbb G}'(n)$. Nonetheless, there is
an ideal $\widetilde{I}$ of $\widetilde{R}$ and an element ${\pmb \delta}$ of $\widetilde{R}$ so that
$\widetilde{I} \widetilde{R}_{{\pmb \delta}}$ is a Gorenstein ideal of $\widetilde{R}_{{\pmb \delta}}$ and $\widetilde{\Bbb G}'(n)_{{\pmb \delta}}$
is a resolution of $\widetilde{R}_{{\pmb \delta}}/\widetilde{I} \widetilde{R}_{{\pmb \delta}}$ by projective $\widetilde{R}_{{\pmb \delta}}$-modules. 

The complex $\widetilde{\Bbb G}'(n)$ is obtained from a less complicated complex $\widetilde{\Bbb G}(n)$ which
is built directly, and in a polynomial   manner, from the coefficients of a generic Macaulay inverse system $\Phi$. Furthermore, $\widetilde{I}$ is the ideal of $\widetilde{R}$ determined by $\Phi$. The modules of $\widetilde{\Bbb G}(n)$ are Schur and Weyl modules corresponding to hooks. 
The complex $\widetilde{\Bbb G}(n)$ is  bi-homogeneous and every entry of every  matrix in $\widetilde{\Bbb G}(n)$ is a  monomial. 

If $m_1,\dots,m_N$ is a list of the monomials in $x_1,\dots,x_d$ of degree $n-1$, then    $\pmb \delta$ is the determinant of the $N\times N$ matrix $(t_{m_im_j})$. The previously listed results  exhibit  a flat family of 
$\pmb k$-algebras parameterized by $\Bbb I_{n}^{[d]}(\pmb k)$: 
$$\tsize \pmb k[\{t_M\}]_{\pmb \delta} \to \left(\frac{\pmb k\otimes_{\Bbb Z}\widetilde{R}}{\vphantom{\widetilde{\widetilde{I}}}\widetilde{I}}\right)_{\pmb \delta}.\tag*$$ Every algebra $P/I$, with $I\in \Bbb I_{n}^{[d]}(\pmb k)$, is a fiber of (*). We simultaneously resolve all of these algebras $P/I$.

The natural action of $\operatorname{GL}_d(\pmb k)$ on $P$ induces an action of 
$\operatorname{GL}_d(\pmb k)$ on $\Bbb I_n^{[d]}(\pmb k)$. We prove that  if  
   $d=3$, $n\ge 3$, and the characteristic of $\pmb k$ is zero, then $\Bbb I^{[d]}_n(\pmb k)$ decomposes into at least four disjoint, non-empty orbits under this group action.}

\heading  Table of Contents \endheading

\halign{
#\hfil&\quad#\hfil&\quad#\hfil&\quad#\hfil&\quad#\hfil\cr
0.&Introduction.\cr
1.&Terminology, notation, and preliminary results.\cr
2.&The complexes $\Bbb L(\Psi,n)$ and $\Bbb K(\Psi,n)$.\cr
3.&The generators.\cr
4.&The main theorem.\cr
5.&Examples of the resolution $\widetilde{\Bbb G}(n)$.\cr
6.&The minimal resolution.\cr 
7.&Non-empty disjoint sets of orbits.\cr
}

 \SectionNumber=\zero\tNumber=1
\heading Section \number\SectionNumber. \quad Introduction.
\endheading 

Fix a  pair of positive integers $d$ and $n$.  We create a ring $\widetilde{R}$ and a complex $\widetilde{\Bbb G}'(n)$ of $\widetilde{R}$-modules with the following universal property. Let $P=\pmb k[x_1,\dots,x_d]$ be a polynomial ring in $d$ variables over the field $\pmb k$ and let $I$ be a grade $d$ Gorenstein ideal in $P$ which is generated by homogeneous forms of degree $n$.
If the resolution of $P/I$ by free $P$-modules is (Gorenstein) linear, then there exists a ring homomorphism $\widehat{\phi}\: \widetilde{R}\to P$ such that $P\otimes_{\widetilde{R}}\widetilde{\Bbb G}'(n)$  is a minimal homogeneous resolution of $P/I$ by free $P$-modules. Our construction is coordinate free. 

We briefly describe our construction, many more details will be given later. Let $U$ be a free Abelian group of rank $d$. The ring $\widetilde{R}$ is equal to $$\operatorname{Sym}_{\bullet}^{\Bbb Z}(U\oplus \operatorname{Sym}_{2n-2}^{\Bbb Z}U),$$ and the complex  
$\widetilde{\Bbb G}'(n)$ has the form
$$0\to Y @>   >> X_{d-1,n}@> \operatorname{Kos}^{\Psi}>> X_{d-2,n}@>\operatorname{Kos}^{\Psi}>>\dots @>\operatorname{Kos}^{\Psi}>>
X_{1,n}@> \operatorname{Kos}^{\Psi}>> X_{0,n}@> \widehat{\Psi} >>\widetilde{R},$$
where $X_{p,n}$ is a submodule of $\widetilde{R}\otimes _{\Bbb Z} \bigwedge_{\Bbb Z}^pU\otimes_{\Bbb Z}\operatorname{Sym}_n^{\Bbb Z}U$ and $Y$ is a free $\widetilde{R}$-module of rank one. The $\Bbb Z$-module homomorphism $\Psi\:U\to \widetilde{R}$ is inclusion; $\widehat{\Psi}\:\operatorname{Sym}_\bullet^{\Bbb Z}U\to \widetilde{R}$ is the $\Bbb Z$-algebra homomorphism induced by $\Psi$; and $$\tsize \operatorname{Kos}^\Psi
\:\widetilde{R}\otimes_{\Bbb Z} \bigwedge^p_{\Bbb Z}U\to \widetilde{R}\otimes_{\Bbb Z} \bigwedge^{p-1}_{\Bbb Z}U$$ is the Koszul complex map induced by $\Psi$. 

If one chooses a basis $x_1,\dots,x_d$ for the free Abelian group $U$ and a basis $\{t_M\}$, for the free Abelian group $\operatorname{Sym}_{2n-2}^{\Bbb Z}$, as $M$ roams over all monomials in $x_1,\dots,x_d$ of degree $2n-2$, then  one may think of $\widetilde{R}$ as the polynomial ring $\Bbb Z[x_1,\dots,x_d,\{t_M\}]$. There is a distinguished element $\pmb \delta=\det(t_{m_imj})$ of  $\widetilde{R}$, where $m_1,\dots,m_N$ is a list of the monomials in $x_1,\dots,x_d$ of degree $n-1$. The ${\widetilde{R}}_{\pmb \delta}$-modules $(X_{p,r})_{\pmb \delta}$ are projective and the complex ${\widetilde{\Bbb G}'(n)}_{\pmb \delta}$ is a resolution of ${\widetilde{R}}_{\pmb \delta}/\widetilde{I}{\widetilde{R}}_{\pmb \delta}$, where $\widetilde{I}{\widetilde{R}}_{\pmb \delta}$ is a grade $d$ Gorenstein ideal in ${\widetilde{R}}_{\pmb \delta}$.

Let $P$ be a standard graded polynomial ring $\pmb k[x_1,\dots,x_d]$ over a field $\pmb k$, and $I$ be a homogeneous ideal of $P$ so that the quotient ring $P/I$ is Artinian and Gorenstein. Recall that $P/I$ is said to have a {\it $($Gorenstein$)$ linear resolution over $P$} if    the minimal homogeneous resolution of $P/I$ by free $P$-modules 
has the form 
$$\alignedat1 0\to P(-2n-d+2)\to P(-n-d+2)^{\beta_{d-1}}&\to \dots \to P(-n-j+1)^{\beta_j}\to \dots \\
&\to P(-n-1)^{\beta_2}\to   P(-n)^{\beta_1} \to P,\endalignedat$$
for some positive integer $n$. (The Betti numbers
$\beta_j$, with $1\le j\le d-1$, may be computed using the Herzog-K\"uhl formulas \cite{\rref{HK}}.) The hypothesis that $P/I$ is Gorenstein forces the entries in the first and last matrices in the minimal homogeneous resolution  to have the same degree; all of the entries of the other matrices are linear forms.   If $\pmb k$ is a field and $d$ and $n$ are positive integers, then let $\Bbb  I_n^{[d]}(\pmb k)$ be the following set of ideals:
$$ \Bbb  I_n^{[d]}(\pmb k)=\left\{I\left \vert \matrix\format\l\\  \text{$I$ is a homogeneous, grade $d$, Gorenstein ideal}\\  \text{of $P=\pmb k[x_1,\dots,x_d]$ with a linear resolution,}\\  \text{and $I$ is generated by forms of degree $n$}\endmatrix \right. \right\}.\tag\tnum{In}$$
 If $I$ is in $ \Bbb  I_n^{[d]}(\pmb k)$, then the socle degree of $P/I$ is $2n-2$ and the Hilbert Function of $P/I$ is
$$\dim_{\pmb k} [P/I]_i=\dim_{\pmb k}[P]_i=\dim_{\pmb k} [P/I]_{2n-2-i}\quad \text{for $0\le i\le n-1$}.\tag\tnum{comp}$$ 
Tony Iarrobino \cite{\rref{I84}} initiated the use of the word ``compressed'' to describe the rings of (\tref{comp}).   Among all Artinian Gorenstein $\pmb k$-algebras with the specified embedding codimension and the specified socle degree, these have the largest total length. (The set of Artinian Gorenstein algebras with a linear resolution is a proper subset of the set of compressed Artinian Gorenstein algebras.)

We consider the following questions and projects.
  
\smallskip \noindent{\bf Project  \tnum{EK1}.} Parameterize the elements of $ \Bbb  I_n^{[d]}(\pmb k)$, and their resolutions, in a reasonable manner.

\smallskip \noindent{\bf Project \tnum{EK2}.} The group $\operatorname{GL}_{d}\pmb k$ acts on   the set $ \Bbb  I_n^{[d]}(\pmb k)$ (by way of the action of $\operatorname{GL}_{d}\pmb k$ on the vector space with basis $x_1,\dots,x_d$). Decompose  $ \Bbb  I_n^{[d]}(\pmb k)$   into disjoint orbits under this group action.

\smallskip \noindent{\bf Project \tnum{EK3}.} (This is a different way to say Project \tref{EK2}.) Classify all graded,   Artinian, Gorenstein, $\pmb k$-algebras with linear resolutions, of embedding codimension $d$ and socle degree $2n-2$.

\smallskip \noindent{\bf Question \tnum{EK4}.} If the complete answer to (\tref{EK2}) and (\tref{EK3}) is elusive, maybe one can answer: How many classes are there?

 \smallskip \noindent{\bf  Question \tnum{EK5}.} 
If the complete answer to (\tref{EK4}) is elusive, maybe one can answer:  Is there more than one class? 

\smallskip  Question \tref{EK5} is already interesting when $d=3$. According to Buchsbaum and Eisenbud \cite{\rref{BE}}, every grade three Gorenstein ideal is generated by the maximal order Pfaffians of an odd sized alternating matrix $X$. We consider such ideals when each entry of $X$ is a linear form from $P=\pmb k[x,y,z]$.  Let $\pmb k$ be a fixed field and $n$ be a fixed positive integer. Consider 
 $$ \Bbb  X_n(\pmb k)=\left \{X\left\vert\matrix \format\l\\ \text{$X$ is a $(2n+1)\times (2n+1)$ alternating matrix of linear}\\ \text{forms from $P=\pmb k[x,y,z]$ such that the ideal generated}\\ \text{by the maximal order Pfaffians of $X$ has grade $3$}\endmatrix  \right.\right\}.\tag{\tnum{p1}}$$
 If $X$ is in $ \Bbb  X_n(\pmb k)$, then the ideal generated by the maximal order Pfaffians of $X$ is in $ \Bbb  I_n^{[3]}(\pmb k)$. Buchsbaum and Eisenbud exhibited an element $H_n\in  \Bbb  X_n(\pmb k)$ for all $n$.  The first two matrices $H_n$ are 
$$\tsize H_1=\bmatrix 0&x&z\\-x&0&y\\-z&-y&0\endbmatrix\qquad \text{and}\qquad H_2=\left[\smallmatrix 0&x&0&0&z\\-x&0&y&z&0\\0&-y&0&x&0\\0&-z&-x&0&y\\-z&0&0&-y&0\endsmallmatrix\right]\tag \tnum{Hn}$$
One follows the same pattern to build all of the rest of the $H_n$.

\smallskip  \noindent{\bf  Question \tref{EK5}, when $d=3$} is: Can every  element of $ \Bbb  X_n(\pmb k)$   be put in the form of $H_n$
after row and column operations and linear change of variables?

\smallskip 
The present paper  contains a complete solution  to Project \tref{EK1}, a complete answer to Question \tref{EK5} when $d=3$, and a partial answer to Question \tref{EK4} when $d=3$. We remain very interested in Projects \tref{EK2} and \tref{EK3}.

To parameterize the elements of $\Bbb I_n^{[d]}(\pmb k)$, we use Macaulay inverse systems. View the polynomial ring $P=\pmb k[x_1,\dots,x_d]$ as the symmetric algebra $\operatorname {Sym}^{\pmb k}_{\bullet}U$, where $U$ is the $\pmb k$-vector space with basis $x_1,\dots,x_d$. The divided power algebra $D_{\bullet}^{\pmb k}(U^*)$ is a module over  the polynomial ring $\operatorname {Sym}_{\bullet}^{\pmb k}U$. (See Subsection \number\prelim.\number\DP.) Let $I$ be a grade $d$ Gorenstein ideal in $\operatorname {Sym}_{\bullet}^{\pmb k}U$. One application of Macaulay's Theorem (see Theorem \tref{MT}) is that the annihilator of $I$ in $D_{\bullet}^{\pmb k}(U^*)$ is a cyclic $\operatorname {Sym}_{\bullet}^{\pmb k}U$-module, denoted $\operatorname{ann} I$, and called the Macaulay inverse system of $I$.
If $I$ is in $\Bbb I_n^{[d]}(\pmb k)$, then the Macaulay inverse system of $I$ is generated by an element $\phi$ of $D_{2n-2}^{\pmb k}(U^*)$. 
In Corollary \tref{open} we identify an open subset of the affine space $D_{2n-2}^{\pmb k}(U^*)$ which parameterizes $\Bbb I_n^{[d]}(\pmb k)$.

As an intermediate step toward   our solution  of Project \tref{EK1}, we exhibit a complex $\widetilde{\Bbb G}(n)$ that depends only 
on the positive integer $n$. The complex $\widetilde{\Bbb G}(n)$ is built over the bi-graded polynomial ring 
$$\widetilde{R}=\Bbb Z\left[x_1,\dots,x_d,\left\{t_M\left\vert\matrix\format\l\\ \text{$M$ is a monomial in $\{x_1,\dots,x_d\}$}\\\text{of degree $2n-2$}\endmatrix\right.\right\}\right];\tag\tnum{gener}$$ where the variables $x_i$ each have degree $(1,0)$ and the variables $t_M$ each have degree $(0,1)$. 
 Each element   $\phi\in D^{\pmb k}_{2n-2}(U^*)$  induces a $\Bbb Z[x_1,\dots,x_d]$-algebra homomorphism   $\widehat{\phi}\:\widetilde{R}\to \pmb k[x_1,\dots,x_d]=P$, with $\widehat{\phi}(t_M)= \phi(M)$. In practice, $\phi=\sum \tau_MM^*$ in $D^{\pmb k}_{2n-2}(U^*)$, where $M$ varies over the monomials of $P=\operatorname {Sym}_{\bullet}^{\pmb k}(U)$ of degree $2n-2$, and $\{M^*\}$ is the basis for $D^{\pmb k}_{2n-2}(U^*)$ which is dual to $\{M\}$. The map $\widehat{\phi}\: \widetilde{R}\to P$, which is induced by $\phi$, sends the variable $t_M$ to $\tau_M$, which is an element of the field $\pmb k$. 

\proclaim{Theorem \tnum{MT2}}  Fix  $n$ and $\widetilde{\Bbb G}(n)$ as described in {\rm(\tref{gener})}. Let $m_1,\dots,m_N$ be a list of the monomials in $x_1,\dots,x_d$ of degree $n-1$.
\item{\rm(a)} If  $\pmb \delta$ is the determinant  of the $N\times N$ matrix $T=(t_{m_im_j})$, then the localization $\widetilde{\Bbb G}(n)_{\pmb \delta}$ is a resolution of a Gorenstein ring $\widetilde{R}_{\pmb \delta}/\widetilde{I}\widetilde{R}_{\pmb \delta}$ by free 
$\widetilde{R}_{\pmb \delta}$-modules.
\item{\rm(b)}  Let $\pmb k$ be a field, $P$ be the polynomial ring $\pmb k[x_1,\dots,x_d]$, and $U$ be the $\pmb k$-vector space with basis $x_1,\dots,x_d$. 
If  $I=\operatorname{ann}\phi$ is an element of $\Bbb I_n^{[d]}(\pmb k)$ and $P$ is an $\widetilde{R}$-algebra by way of  $\widehat{\phi}$, then $P\otimes_{\widetilde R}\widetilde{\Bbb G}(n)$ is a resolution of $P/I$ by free $P$-modules.
\item{\rm(c)}  Some of the features of $\widetilde{\Bbb G}(n)$ are {\bf (1)} it is bi-graded; {\bf (2)} there is one $\widetilde{\Bbb G}(n)$ for each $n$; and {\bf (3)}  every entry of every  matrix in $\widetilde{\Bbb G}(n)$ is a  monomial.
\endproclaim

Parts (a) and (b) of Theorem \tref{MT2} are established in Corollaries \tref{EK-K-g} and \tref{EK-K-2}, respectively; the bi-graded Betti numbers of $\widetilde{\Bbb G}(n)$ are exhibited in Remark \tref{bih}; assertion (c.3) is Theorem \tref{mon}. 

 Feature (c.2) is interesting because, even when $d$ is only $3$, there are at least 4 disjoint families of minimal resolutions; see Theorem \tref{atlst}; but $\widetilde{\Bbb G}(n)$ specializes to all of these families.  Feature (c.3) is interesting because, for example, the matrix which  presents $I=\operatorname{ann} \phi$ is monomial and linear. Often the matrix that presents a module is more important than the generators of the module. One uses the presentation of $M$ to compute $F(M)$
for any right exact functor $F$ and any module. One uses the presentation to compute $\operatorname {Sym}_{\bullet}M$, and $\operatorname {Sym}_{\bullet}I$ is the first step toward studying the blow-up algebras -- in particular the Rees algebra -- associated to $I$. We have an uncomplicated presenting matrix!  

On the other hand, $P\otimes _{\widetilde{R}}\widetilde{\Bbb G}(n)$ is not minimal. 
The final step in our solution  of Project \tref{EK1} is the complex $\widetilde{\Bbb G}'(n)$ which is described at the beginning of the Introduction and takes  place in Section \number\MR.

An alternate phrasing of Theorem \tref{MT2} is also given in Corollary \tref{EK-K-2}
where we exhibit  a flat family of 
$\pmb k$-algebras parameterized by $\Bbb I_{n}^{[d]}(\pmb k)$: 
$$\tsize \pmb k[\{t_M\}]_{\pmb \delta} \to \left(\frac{\pmb k\otimes_{\Bbb Z}\widetilde{R}}{\vphantom{\widetilde{\widetilde{I}}}\widetilde{I}}\right)_{\pmb \delta}.\tag\tref{bblop}$$ Every algebra $\pmb k[x_1,\dots,x_d]/I$, with $I\in \Bbb I_{n}^{[d]}(\pmb k)$, is a fiber of (\tref{bblop}). We simultaneously resolve all of these algebras $\pmb k[x_1,\dots,x_d]/I$.

\bigskip In order to provide some insight into the structure of the complex $\widetilde{\Bbb G}(n)$, 
we next describe the complex $P\otimes_{\widetilde{R}}\widetilde{\Bbb G}(n)$, where the polynomial ring $P=\pmb k[x_1,\dots,x_d]=\operatorname{Sym}_\bullet^{\pmb k} U$, the  Macaulay inverse system $\phi\in D_{2n-2}^{\pmb k}(U^*)$, and the $\Bbb Z[x_1,\dots,x_d]$-algebra homomorphism $\widehat{\phi}\: \widetilde{R}\to P$ have all been fixed. (In this discussion, $U$ is the vector space over $\pmb k$ with basis $x_1,\dots,x_d$.) The complex $P\otimes_{\widetilde{R}}\widetilde{\Bbb G}(n)$ is the mapping cone of $$ \eightpoint \matrix   & &0&\to& L_{d-1,n}&\to &\cdots &\to&L_{1,n}&\to&L_{0,n}&\to&P\\  & &\downarrow & &\downarrow&&& & \downarrow&&\downarrow&\\ 0&\to &P\otimes_{\pmb k}{\tsize\bigwedge}^d_{\pmb k}U&\to&K_{d-1,n-2}&\to&\cdots &\to&K_{1,n-2}&\to&K_{0,n-2},
\endmatrix\tag\tnum{pi}$$where the top complex of (\tref{pi}) is a resolution of $P/J^n$ by free $P$-modules for $J=(x_1,\dots,x_d)$, and the bottom complex of (\tref{pi}) is a resolution of the canonical module of $P/J^{n-1}$ by free $P$-modules. The $P$-modules $L_{i,n}$  and $K_{i,n-2}$ are   contained in $P\otimes_{\pmb k}\bigwedge^i_{\pmb k} U\otimes _{\pmb k}\operatorname{Sym}_n^{\pmb k} U$ and  $P\otimes_{\pmb k}\bigwedge^i_{\pmb k} U\otimes _{\pmb k}D_{n-2}^{\pmb k} (U^*)$, respectively; and the vertical map $L_{i,n}\to K_{i,n-2}$ is induced by the map  $\pmb p_n^{\phi}\: \operatorname{Sym}_n^{\pmb k} U\to D_{n-2}^{\pmb k} (U^*)$ which  sends the element $u_n$ of $\operatorname{Sym}_n^{\pmb k} U$ to $u_n(\phi)$ in $D_{n-2}^{\pmb k} (U^*)$. The top complex of  (\tref{pi}) is a well-known complex; it is called $L^1_n(\Psi)$ in \cite{\rref{BE75}, Cor.~3.2} (see also \cite{\rref{S}, Thm.~2.1}), where $\Psi\:P\otimes_{\pmb k} U\to P$ is the  map given by multiplication in $P$. Of course, the bottom complex of (\tref{pi}) is isomorphic to $\operatorname{Hom}_P(L^1_{n-1}(\Psi),P)$. The properties of mapping cones automatically yield that $P\otimes_{\widetilde{R}}\widetilde{\Bbb G}(n)$, which  is the mapping cone of (\tref{pi}), is a resolution. The interesting step involves calculating the zeroth homology of this complex; see the proof of Theorem \tref{EK-K'}. 

The bi-graded complex $\widetilde{\Bbb G}(n)$ is very similar to the graded complex $P\otimes_{\widetilde{R}}\widetilde{\Bbb G}(n)$. Indeed, we can use the double complex (\tref{pi}) to indicate the bi-degrees and the bi-homogeneous Betti number in $\widetilde{\Bbb G}(n)$. Recall the bi-degrees of the variables of $\widetilde{R}$ from (\tref{gener}): the $x_i$ have bi-degree $(1,0)$ and the $t_M$ have bi-degree $(0,1)$. The top complex of (\tref{pi}) only involves the $x_i$ (even in the bi-homogeneous  construction). The matrix on the far right has bi-homogeneous entries of degree $(n,0)$; the other matrices in the top complex have bi-homogeneous entries of degree $(1,0)$. The bottom  complex of (\tref{pi}) also only involves the $x_i$ (even in the bi-homogeneous  construction). The matrix on the far left has bi-homogeneous entries of degree $(n-1,0)$; the other matrices in the bottom complex have bi-homogeneous entries of degree $(1,0)$.  The vertical maps in the bi-homogeneous construction have bi-homogeneous entries of degree $(0,1)$. There are many ways to compute the Betti numbers in (\tref{pi}) (hence in $\widetilde{\Bbb G}(n)$). A formula is given in the Buchsbaum-Eisenbud paper \cite{\rref{BE75}}; the $L$'s and $K$'s  are Schur modules and Weyl modules, respectively, and one  can use combinatorial techniques to calculate the ranks; or one can use the Herzog-K\"uhl formula \cite{\rref{HK}}. Our answer is given in Remark \tref{bih}.

We denote the top complex of (\tref{pi}) by  $\Bbb L(\Psi,n)$ and the bottom complex by  $\Bbb K(\Psi,n-1)$, where, again, $\Psi\:P\otimes _{\pmb k} U\to P$ is multiplication in $P$. Our discussion of the complexes $\Bbb L$ and $\Bbb K$ is contained in Section \number\Lpq. If one replaces  $\Bbb L(\Psi,n)$ and  $\Bbb K(\Psi,n-1)$ in (\tref{pi}) by $\Bbb L(\Psi,n+\rho)$ and  $\Bbb K(\Psi,n-1-\rho)$, respectively, for some positive integer $\rho$, then   the mapping cone of the resulting double complex is a resolution  $P/J^\rho I$. The ideal $J^\rho I$ is a ``truncation'' of $I$ in the sense that $J^\rho I$ is equal to $I_{\ge n+\rho}$; in other words,  $J^\rho I$ is generated by all elements in $I$ of degree at least $n+\rho$. Our construction   works for  $0\le \rho\le n-2$. On the other hand, $\rho\le n-2$ is not a restrictive constraint because $I_{\ge r}=J^{r}$ for $2n-1\le r$, since the socle degree of $P/I$  is  $2n-2$, and $\Bbb L(\Psi,r)$ is a resolution of $P/J^r$.

\bigskip
We now turn our attention to  Questions \tref{EK5} and \tref{EK4}, when $d=3$; namely ``How many orbits does $\Bbb X_n(\pmb k)$ have under the action of  $\operatorname{GL}_{2n+1}\pmb k\times \operatorname{GL}_3\pmb k$?'' If $n=2$, then it is well-known (and not particularly hard to see) that  $ \Bbb  X_n(\pmb k)$ has only one class; see, for example, Observation \tref{n=2}.  Indeed, $\Bbb I_2$ consists of those ideals in $\pmb k[x,y,z]$ which  define homogeneous Gorenstein rings of minimal multiplicity. The $n=2$ case causes one to consider multiplication tables, which is one of the key ideas in our work. The $n=2$ case also reminds us of the resolutions  of Kurt Behnke \cite{\rref{Be}, \rref{Be'}} and Eisenbud, Riemenschneider, and Schreyer \cite{\rref{ERS}}.

Fix  integers $n$ and $\mu$ with  $3\le n$ and $0\le \mu\le 3$, 
let $$ \Bbb  I^{[3]}_{n,\mu}(\pmb k)=\left\{I\in  \Bbb  I_n^{[3]}(\pmb k)\left\vert \matrix\format\l\\ \exists \text{ linearly independent linear forms}\\\text{$\ell_1,\dots,\ell_\mu$ in $P_1$ with $\ell_1^n,\dots,\ell_\mu^n$ in $I$}\\\text{and 
$\not\exists$ $\mu+1$ such forms}\endmatrix\right. \right\}.\tag\tnum{p4}$$
It is clear that $ \Bbb  I_n^{[d]}(\pmb k)$ is the disjoint union of $\bigcup\limits_{\mu=0}^3 \Bbb  I^{[3]}_{n,\mu}(\pmb k)$ and each $ \Bbb  I^{[3]}_{n,\mu}(\pmb k)$ is closed under the action of $\operatorname{GL}_{2n+1}\pmb k\times \operatorname{GL}_{3}\pmb k$.  

\proclaim{Theorem \tnum{atlst}} If $n\ge 3$ and the characteristic of $\pmb k$ is zero, then $ \Bbb  I^{[3]}_{n,\mu}(\pmb k)$ is non-empty for $0\le \mu\le 3$. In particular, if $n\ge 3$, then  $ \Bbb  X_n(\pmb k)$ has at least four non-empty, disjoint orbits in the sense of Project {\rm\tref{EK2}}.
 \endproclaim 

The proof of Theorem \tref{atlst} is carried out in Section \number\nonempty. To prove the result, we exhibit an element of $\Bbb I^{[3]}_{n,\mu}(\pmb k)$ for each $\mu$. It turns out that the ideal $\operatorname{BE}_n$, generated by the maximal order Pfaffians of the Buchsbaum-Eisenbud matrix $H_n$ of  (\tref{Hn}), is in $\Bbb I^{[3]}_{n,2}(\pmb k)$;  $$J_{n,n-1}=(x^{n},y^{n},z^{n}):(x+y+z)^{n-1}$$ is in $\Bbb I^{[3]}_{n,3}(\pmb k)$; and by modifying a homogeneous generator of the  Macaulay inverse system for $\operatorname{BE}_n$ we produce  ideals in $\Bbb I^{[3]}_{n,0}(\pmb k)$ and $\Bbb I^{[3]}_{n,1}(\pmb k)$.    We do {\bf not} claim that every ideal in $ \Bbb  I^{[3]}_{n,\mu}(\pmb k)$ may be converted into any other ideal in $ \Bbb  I^{[3]}_{n,\mu}(\pmb k)$ by using $\operatorname{GL}_{2n+1}\pmb k\times \operatorname{GL}_{3}\pmb k$. 
So, Projects \tref{EK2} and \tref{EK3} are far from resolved, even when $d=3$.

 There are at least two motivations for these projects. First of all, the ideals $J_{n,n-1}$ arise naturally in the study of the Weak Lefschetz Property (WLP)  for monomial complete intersections in characteristic $p$.  Let $\pmb k$ be an infinite field and   $A$ be a standard graded Artinian $\pmb k$-algebra. The ring $A$ has the WLP   if, for every general linear form $\ell$, multiplication by $\ell$ from $[A]_i$ to $[A]_{i+1}$  is a  map of maximal rank  for all $i$. (That is, each of these maps is  either injective or surjective.) Stanley \cite{\rref{St80}} used the hard Lefschetz Theorem from Algebraic Geometry to show that, if $\pmb k$ is the field of complex numbers, then every monomial complete intersection has the WLP. 
Alternate proofs of Stanley's Theorem, requiring only that $\pmb k$ have characteristic zero, and using techniques from other branches of mathematics, may be found in   \cite{\rref{RRR}} and \cite{\rref{HMNW}}. The story is much different in positive characteristic. Let $A(\pmb k,m,n)=\pmb k[x_1,\dots,x_m]/{(x_1^n,\dots,x_m^n)}$, where $\pmb k$ is an infinite field of positive characteristic $p$. If $m=3$, then Brenner and Kaid \cite{\rref{BK}} have identified all $n$, as a function of $p$,  for which $A$ has the WLP. For a given prime $p$, there are intervals of $n$ for which  $A(\pmb k,3,n)$ has the WLP and the position of these intervals is related to the Hilbert-Kunz multiplicity of the Fermat ring $\pmb k[x,y,z]/(x^n+y^n+z^n)$ as studied by Han \cite{\rref{han}} and Monsky \cite{\rref{M}}. If $n=4$, then it is shown in \cite{\rref{WLP}} that $A(\pmb k,4,n)$ rarely has the WLP and only for discrete values of $n$. Furthermore,    the following statement holds.

\proclaim{Observation \tnum{hype}}  Let $\pmb k$ be an  infinite field and, for each pair of positive integers $(a,b)$, let $J_{a,b}$ be the ideal $(x^a,y^a,z^a)\!:\!(x+y+z)^b$ of $\pmb k[x,y,z]$. If $J_{n,n-1}\in \Bbb I_n^{[3]}(\pmb k)$ or $J_{n,n+1}\in \Bbb I_{n-1}^{[3]}(\pmb k)$, then 
$A(\pmb k,4,n)$ has the {\rm WLP}. \endproclaim

\smallskip \flushpar For a complete, up-to-date, history of the WLP see \cite{\rref{MN}}.

Also, 
there is much recent work concerning  the  equations
that define the Rees algebra of   ideals which are primary to the maximal ideal; see, for example, \cite{\rref{HSV},$\allowbreak$\rref{CHW},$\allowbreak$\rref{Bu},$\allowbreak$\rref{CKPU}}. The driving force behind this work is the desire to understand the singularities of parameterized curves or surfaces; see \cite{\rref{SCG},\rref{CWL},\rref{BB},\rref{Bo},\rref{BD'A}} and especially \cite{\rref{CKPU}}. One of the key steps in \cite{\rref{CKPU}} is the decomposition of the space of balanced Hilbert-Burch matrices into disjoint orbits under the action of $\operatorname{GL}_3\pmb k\times \operatorname{GL}_2\pmb k$. The present paper includes a preliminary step toward obtaining a comparable decomposition of the space of syzygy matrices for the set of linearly presented grade three Gorenstein ideals.

\SectionNumber=\prelim\tNumber=1
\heading Section \number\SectionNumber. \quad Terminology, notation, and preliminary results.
\endheading

This section contains preliminary material. It is divided into four subsections: Miscellaneous information, Pfaffian conventions, Divided power structures, and Macaulay inverse systems.
\heading Subsection \number\prelim.\number\MC\quad Miscellaneous   information \endheading 
In this paper, $\pmb k$ is always a field.  Unless otherwise noted, the polynomial ring $\pmb k[x_1,\dots,x_d]$ is assumed to be a standard graded $\pmb k$-algebra; that is, each variable has degree one. 
For each   graded module $M$ we use $[M]_i$ to denote the {\it homogeneous component of $M$ of degree $i$}.



If $\alpha$ is a real number then $\lfloor \alpha\rfloor$ is the {\it round down} of $\alpha$; that is $\lfloor \alpha\rfloor$ is equal to the integer $n$ with $n\le \alpha<n+1$. We use   $\Bbb Z$ and $\Bbb Q$ to represent the ring of integers  and  the field of rational numbers, respectively.

Recall that if $A$, $B$, and $C$ are $R$-modules, then the $R$-module homomorphism $F\: A\otimes_R B\to C$ is a {\it perfect pairing} if the induced $R$-module homomorphisms $A\to \operatorname{Hom}_R(B,C)$ and $B\to \operatorname{Hom}_R(A,C)$, which are given by $a\mapsto F(a\otimes \underline{\phantom{x}})$ and $b\mapsto F(\underline{\phantom{x}}\otimes b)$, are isomorphisms.

If $V$ is a free $R$-module, then $$\tsize \Delta\: \bigwedge_R^\bullet V \to  \bigwedge_R^\bullet V \otimes_R\bigwedge_R^\bullet V \tag\tnum{delta}$$ is the usual co-multiplication map in the exterior algebra. In particular, the component $\Delta\: \bigwedge_R^a V \to  \bigwedge_R^1 V \otimes_R\bigwedge_R^{a-1} V$  of (\tref{delta}) sends
$$v_1\wedge \dots \wedge v_a \quad\text{to}\quad \sum_{i=1}^a (-1)^{i+1}v_i\otimes v_1\wedge\dots \widehat{v_i}\dots \wedge v_a,$$ for $v_1,\dots,v_a\in V$.

Let $\pmb k$ be a field and $A=\bigoplus_{0\le i}[A]_i$ be a graded Artinian $\pmb k$-algebra, with $[A]_0=\pmb k$ and maximal ideal $\frak m_A=\bigoplus_{1\le i} [A]_i$. The {\it socle} of $A$ is the  $\pmb k$-vector space $$0:_A\frak m_A=\{a\in A\mid a\frak m_A=0\}.$$ The Artinian ring $A$ is {\it Gorenstein} if the socle of $A$ has dimension one. In this case, the degree of a generator of the socle  of $A$ is called the {\it socle degree} of $A$. If $A$ is a graded Artinian  Gorenstein $\pmb k$-algebra with socle degree $s$, then the multiplication map 
$$[A]_i\otimes_{\pmb k} [A]_{s-i}\to [A]_s$$ is a perfect pairing for $0\le i\le s$. 
  The homogeneous ideal $I$ of the polynomial ring $P=\pmb k[x_1,\dots,x_d]$ is a {\it grade $d$ Gorenstein ideal} if $P/I$ is an Artinian Gorenstein $\pmb k$-algebra. 

Let  $R$ be an arbitrary commutative Noetherian ring.   The {\it grade} of a proper ideal $I$ in a  $R$ is the length of a maximal $R$-sequence contained in $I$. The ideal $I$ is called {\it perfect} if the grade of $I$ is equal to the projective dimension of $R/I$ as an $R$-module. (The inequality $\operatorname{grade}I\le \operatorname{proj.\ dim.}_R R/I$ always holds.) An ideal $I$ of grade $g$ is a {\it Gorenstein ideal} if $I$ is perfect and $\operatorname{Ext}^g_R(R/I,R)$ is a cyclic $R/I$-module.

If $\Bbb D$ is a double complex:
$$\CD 
@. \vdots@. \vdots@.\vdots\\
@. @VVV @VVV @VVV\\
\cdots@>>>D_{2,1}@>>>D_{1,1}@>>>D_{0,1}\\
@. @VVV @VVV @VVV\\
\cdots@>>>D_{2,0}@>>>D_{1,0}@>>>\phantom{,}D_{0,0},
\endCD$$
 then $\Bbb T(\Bbb D)$ is the {\it total complex} of $\Bbb D$. The modules of  $\Bbb T(\Bbb D)$ are
$$\cdots@>>>\Bbb T(\Bbb D)_i@>>> \Bbb T(\Bbb D)_{i-1}@>>> \cdots, $$with $\Bbb T(\Bbb D)_i=\bigoplus_{a+b=i}D_{a,b}$. The maps of $\Bbb T(\Bbb D)$ are the maps of $\Bbb D$ with the signs adjusted.

If $S$ is a statement, then we define
$$\chi(S)=\cases 1&\text{if $S$ is true}\\0&\text{if $S$ is false}.\endcases\tag\tnum{chi}$$

\heading Subsection \number\prelim.\number\PC\quad  Pfaffian Conventions \endheading
An  {\it alternating} matrix is a square skew-symmetric matrix whose entries on the main diagonal are zero. Fix an alternating 
  $n\times n$ matrix $Z$, with entries $(z_{i,j})$. 
The Pfaffian of $Z$, denoted $\operatorname{Pf}(Z)$, is a square root of a determinant of $Z$. We recall that $\operatorname{Pf}(Z)=0$ if $n$ is odd, and $\operatorname{Pf}(Z)=z_{1,2}$ if $n=2$.  
If $a_1,\dots,a_s$ are  integers between $1$ and $n$, then we let $\sigma(a_1,\dots,a_s)$ be zero if there is some repetition  among the indices $a_1,\dots, a_s$, and 
 if the indices $a_1,\dots, a_s$  are distinct, then $\sigma(a_1,\dots,a_s)$ is the sign of the permutation that puts the indices into ascending order. Let
$Z_{a_1,\dots,a_s}$ denote   $\sigma(a_1,\dots,a_s)$ times the Pfaffian of the matrix obtained by deleting rows and columns $a_1,\dots,a_s$ from $Z$. We recall that   $\operatorname{Pf}(Z)$ may be computed along any row or down any column; that is,
$$\operatorname{Pf}(Z)=\cases\sum\limits_j (-1)^{i+j+1}z_{i,j}Z_{i,j},&\text{with $i$ fixed}\\
\sum\limits_i (-1)^{i+j+1}z_{i,j}Z_{i,j}, &\text{with $j$ fixed.}\endcases$$

\heading Subsection \number\prelim.\number\DP \quad Divided power structures\endheading
Let $V$ be a free module of finite rank over the commutative Noetherian ring  $R$. Form 
the polynomial ring $$\operatorname {Sym}_{\bullet}^R(V)=\bigoplus\limits_{0\le r}\operatorname {Sym}_{r}^R(V)$$ and the divided power $R$-algebra $$D^R_{\bullet}(V^*)=\bigoplus\limits_{0\le r}D_{r}^R(V^*),$$ where the functor $^*$ is equal to $\operatorname{Hom}_{R}(\underline{\phantom{x}},R)$.
Each graded component $\operatorname {Sym}_{r}^R(V)$ of $\operatorname {Sym}_{\bullet}^R(V)$ and 
$D_{r}^R(V^*)$ of $D_{\bullet}^R(V^*)$ is a free module of finite rank over $R$. The rules for a divided power algebra are recorded in \cite{\rref{GL},~section 7} or  \cite{\rref{E95},~Appendix~2}. (In practice these rules say that $w^{(a)}$ behaves like $w^a/(a!)$ would behave if $a!$ were a unit in $R$.) One makes   $D_{\bullet}^R(V^*)$ become a $\operatorname {Sym}_{\bullet}^RV$-module by decreeing that each element $v_1$ of $V$ acts like a divided power derivation on $D_{\bullet}^R(V^*)$. That is, ${v_1(w_1^{(a)})=v_1(w_1)\cdot w_1^{(a-1)}}$, $v_1$ of a product follows the product rule from calculus, and, once one knows how $V$ acts on $D_{\bullet}^R(V^*)$, then one knows how all of $\operatorname {Sym}_{\bullet}^RV$ acts on $D_{\bullet}^R(V^*)$. 

In a similar manner one makes   $\operatorname {Sym}_{\bullet}^RV$ become  a module over the divided power $R$-algebra $D_{\bullet}^R(V^*)$. If $w_1$ is in $D_1^R(V^*)$, then $w_1$ is a derivation on $\operatorname {Sym}_{\bullet}^RV$, and $w_1^{(a)}$ acts on $\operatorname {Sym}_{\bullet}^RV$ exactly like $\frac 1{a!}w_1^a$ would act if $a!$ were a unit in $R$.  The $R$-algebra $D_{\bullet}^R(V^*)$ is generated by elements of the form $w_1^{(a)}$, with $w_1\in V^*$; once one knows how these elements act on $\operatorname {Sym}_{\bullet}^RV$, then one knows  how every element of $D_{\bullet}^R(V^*)$ acts on $\operatorname {Sym}_{\bullet}^RV$.

The above actions are defined in a coordinate free manner. It makes sense to see what happens when one picks bases. Let $x_1,\dots, x_d$ be a basis for $V$ and $x_1^*,\dots,x_d^*$ be the corresponding  dual basis for $V^*$. (We have $x_i^*(x_j)$ is equal to the Kronecker delta $\delta_{i,j}$.) The above rules show that 
$$x_1^{a_1}\cdots x_d^{a_d}\left(x_1^{*(b_1)}\cdots x_d^{*(b_d)}\right)=x_1^{*(b_1-a_1)}\cdots x_d^{*(b_d-a_d)}\in D_{\sum b_i-\sum a_i}^R(V^*)$$
and 
$$\tsize x_1^{*(b_1)}\cdots x_d^{*(b_d)}\left(x_1^{a_1}\cdots x_d^{a_d}\right)=\binom {a_1}{b_1}\cdots \binom{a_d}{b_d}x_1^{a_1-b_1}\cdots x_d^{a_d-b_d}\in \operatorname {Sym}_{\sum a_i-\sum b_i}^R(V),$$
where $x_i^{r}=0$ and $x_i^{*(r)}=0$ whenever $r<0$. Notice in particular, that if ${\sum a_i=\sum b_i}$, then
$$x_1^{a_1}\cdots x_d^{a_d}\left(x_1^{*(b_1)}\cdots x_d^{*(b_d)}\right)\quad\text{and}\quad x_1^{*(b_1)}\cdots x_d^{*(b_d)}\left(x_1^{a_1}\cdots x_d^{a_d}\right)$$ are both equal to
$$\cases 1&\text{if $a_j=b_j$ for all $j$}\\0&\text{otherwise}.\endcases$$
In particular, $\operatorname {Sym}_a^RV$ and $D_a^R(V^*)$ are {\bf naturally} dual to one another. 

We make use of the dual of  the natural evaluation map 
$$\operatorname{ev}_j\:D_j^R(V^*)\otimes_R\operatorname {Sym}_j^RV\to R.$$ 
Indeed, $\operatorname{ev}_j^*(1)$ is a well-defined coordinate-free element of $\operatorname {Sym}_j^R(V)\otimes_R D_j^R(V^*)$.
If $\{A_i\}$ and $\{A_i^*\}$ are a pair of dual bases for the free modules  $\operatorname {Sym}_j^RV$ and $D_j^R(V^*)$, then
$$\operatorname{ev}_j^*(1)=\sum_iA_i\otimes A_{i}^*\in \operatorname {Sym}_j^R(V)\otimes_R D_j^R(V^*).\tag\tnum{ev*}$$

\heading Subsection \number\prelim.\number\Mac \quad Macaulay inverse systems\endheading

\proclaim {Theorem \tnum{MT}} {\bf (Macaulay, \cite{\rref{M16}})} Let $U$ be a vector space of dimension $d$ over the field $\pmb k$. Then
there exists a one-to-one correspondence between the set of non-zero homogeneous grade d Gorenstein ideals of $\operatorname {Sym}_{\bullet}^{\pmb k}U$ and the set of non-zero
homogeneous cyclic  submodules of $D_{\bullet}^{\pmb k}(U^*)$ $:$ $$\matrix  \left\{\matrix\text{homogeneous grade $d$ Gorenstein}\\\text{ideals of $\operatorname {Sym}_{\bullet}^{\pmb k}U$}\endmatrix\right\} &\leftrightarrow & \left\{\matrix\text{non-zero homogeneous cyclic}\\\text{submodules of $D_{\bullet}^{\pmb k}(U^*)$}\endmatrix\right\}.\endmatrix $$If $I$ is an ideal from the set on the left, then the corresponding submodule of $D_{\bullet}^{\pmb k}(U^*)$ is $$\operatorname{ann}I=\{w \in D_{\bullet}^{\pmb k}(U^*)\mid rw=0 \text{ in $D_{\bullet}^{\pmb k}(U^*)$, for all  $r$  in $I$}\},$$ and if $M$ is a homogeneous cyclic submodule from the set on the right, then the corresponding ideal is 
$$\operatorname{ann}(M)=\{r\in   \operatorname {Sym}_{\bullet}^{\pmb k}U\mid rM=0 \text{ in }D_{\bullet}^{\pmb k}(U^*)\}.$$
Furthermore, the socle degree of $\frac{\operatorname {Sym}_{\bullet}^{\pmb k}U}I$ is equal to the degree of a homogeneous generator of $\operatorname{ann}I$.  
\endproclaim

\remark{Definition} In the language of Theorem \tref{MT}, the homogeneous cyclic submodule $\operatorname{ann}(I)$ of $D_{\bullet}^{\pmb k}(U^*)$ is the {\it Macaulay inverse system} of the grade $d$ Gorenstein ideal $I$ of $\operatorname {Sym}^{\pmb k}_{\bullet}U$.\endremark

\demo{The outline of a Proof of Theorem {\rm\tref{MT}}}  Let $P$ denote $\operatorname {Sym}_\bullet^{\pmb k} U$ and $^*$ denote $\pmb k$-dual. 
Let $I$ be an ideal from the set on the left and let $s$ be  the socle degree of $P/I$. Fix an isomorphism $[P/I]_s@>\simeq>> \pmb k$. Define $w\in D_s^{\pmb k}(U^*)=\operatorname{Hom}_{\pmb k}(\operatorname {Sym}_s^{\pmb k}(U),\pmb k)$ to be the homomorphism $w\:\operatorname {Sym}^{\pmb k}_s(U) \to \pmb k$ which is  the composition $$\operatorname {Sym}_s^{\pmb k}(U)= P_s\to [P/I]_s @>\simeq >> \pmb k.$$ 
The definition of $w$ shows that $[P]_s\cap\operatorname{ann} w=[I]_s$. The ring $P/I$ is Gorenstein so multiplication $[P/I]_i\otimes_{\pmb k} [P/I]_{s-i}\to [P/I]_s$ is a perfect pairing for all $i$. It follows that $\operatorname{ann}(w)=I$. 

Now let $w$ be a non-zero element of $D_s^{\pmb k} (U^*)$ and let $I=\operatorname{ann} (w)$. The homomorphism $w\: \operatorname {Sym}_s^{\pmb k} U\to \pmb k$ is non-zero; so the vector space $[P/I]_s$ has dimension $1$ and this vector space is contained in the socle of $P/I$. On the other hand, if $r$ is a homogeneous element of $P$ with $r\notin I$ and $\deg r<s$, then the hypothesis $rw\neq 0$ guarantees that there is an homogeneous element  $r'$ of $P$ with $r'r \in P_s$ and $r'rw\neq 0$. In particular, $r$ is not in the socle of $P/I$. Thus, the socle of $P/I$ has dimension one and $P/I$ is Gorenstein.   \qed \enddemo

Let $I$ be a fixed homogeneous grade $d$ Gorenstein ideal in $$P=\pmb k[x_1,\dots,x_d]=\operatorname {Sym}_{\bullet}^{\pmb k}U,$$ 
where $U$ is the $d$-dimensional vector space $\bigoplus_{i=1}^d\pmb kx_i$, and let $\phi\in D_{\bullet}^{\pmb k}(U^*)$  be a homogeneous generator for the Macaulay inverse system of $I$. Proposition \tref{J18} gives many ways to test   if   $I$   is in $\Bbb I_n^{[d]}(\pmb k)$. We are particularly interested in tests that involve $\phi$. First of all, Proposition \tref{J18} shows that $\phi$ must be in $D_{2n-2}^{\pmb k}(U^*)$. The other condition that $\phi$ must satisfy involves a map $\pmb p_{n-1}^\phi$ or a matrix $T_\phi$. These notions are defined whenever $\phi$ is a homogeneous element of $D_{\bullet}^R(V^*)$ of even degree, $R$ is a commutative Noetherian ring, and $V$ is a free $R$-module of finite rank. 

\definition{Definition \tnum{T}} Let   $R$ be a commutative Noetherian ring, $V$ be a free $R$-module of finite rank $d$,    and $\phi$ be an element of $D_{2e}^R(V^*)$, for some positive integer $e$.
\roster\item For each integer $i$, with $0\le i\le 2e$, define the homomorphism $$ \pmb p_i^\phi:\operatorname {Sym}_i^RV\to D_{2e-i}^R(V^*)$$ by $\pmb p_i^\phi(v_i)=v_i(\phi)$, for all $v_i\in \operatorname {Sym}_i^RV$.
\item Fix a basis $x_1,\dots, x_d$ for $V$. Let $N=\binom{e+d-1}e$ and $m_1,\dots,m_N$ be a list of the monomials of degree $e$ in $\operatorname {Sym}_e^RV$. Define $T_\phi$ to be the $N\times N$ matrix $T_\phi=[\phi(m_im_j)]$; that is, the entry of $T_\phi$ in row $i$ and column $j$ is the element $\phi(m_im_j)$ of $R$.\endroster \enddefinition
\remark{\bf Remarks \tnum{R1}} (1) The matrix $T_\phi$ is the matrix for  $\pmb p_e^\phi$ with respect to the basis $m_1,\dots,m_N$ for $\operatorname {Sym}_e^RV$ and the dual basis $m_1^*,\dots,m_N^*$ for $D_e^R(V^*)$. 

\medskip
\flushpar (2) The entries of $T_{\phi}$ are the coefficients of $\phi$ as an element of $D_{2e}^R(V^*)$. Let $\nu$ equal $\binom{2e+d-1}{2e}$. If  $M_1,\dots,M_{\nu}$ is a list of the monomials of degree $2e$ in $\operatorname {Sym}_{2e}^RV$ and $M_1^*,\dots,M_{\nu}^*$ is the basis for $D_{2n-2}^R(V^*)$ which is dual to $M_1,\dots,M_{\nu}$, then $\phi=\sum_i\phi(M_i)\cdot M_i^*$. Of course, each entry of $T_{\phi}$ is a coefficient of $\phi$ because each $m_im_j$ is equal to some $M_k$. 
\medskip
\flushpar (3) The element $\det T_{\phi}$ of $R$ is known as the ``determinant of the symmetric bilinear form'' 
$$\operatorname {Sym}_e^RV\times \operatorname {Sym}_e^RV \to R,$$which sends $(v_{e},v'_e)$ to $\phi(v_ev'_e)$. A change of basis for $V$ changes $\det T_{\phi}$ by a unit of $R$. However, in practice, we only use $\det T_{\phi}$, up to unit. In particular, the phrases ``provided $\det  T_{\phi}$ is a unit of $R$'' and ``$R$ localized at the element $\det  T_{\phi}$'' are meaningful even in a coordinate-free context.

\endremark

\example{Example \tnum{Ex1}} Let $\operatorname{BE}_2$ be the ideal of $\pmb k[x,y,z]$ which is generated by the maximal order Pfaffians of the matrix $H_2$ from (\tref{Hn}). We have $\operatorname{BE}_2=
(x^2,y^2,xz,yz,xy+z^2)$ and the Macaulay inverse system for $\operatorname{BE}_2$ is generated by $\phi=x^*y^*-{z^*}^{(2)}$.  If $m_1=x$, $m_2=y$, $m_3=z$, then  $$T_{\phi}=\bmatrix \phi(x^2)&\phi(xy)&\phi(xz)\\\phi(yx)&\phi(y^2)&\phi(yz)\\\phi(zx)&\phi(zy)&\phi(z^2)\endbmatrix=\bmatrix 0&1&0\\1&0&0\\0&0&-1\endbmatrix.$$  
\endexample

\proclaim{Proposition \tnum{J18}} Let $\pmb k$ be a field, $I$ be a homogeneous, grade $d$, Gorenstein ideal in $P=\pmb k[x_1,\dots,x_d]$, $U$ be the $d$-dimensional vector space $[P]_1$, and $\phi\in D_{\bullet}^{\pmb k}(U^*)$  be a homogeneous generator for the Macaulay inverse system of $I$. Then the  following statements are equivalent\,{\rm:}
\roster
\item $I$ is in $\Bbb I_n^{[d]}(\pmb k)$,
\item the minimal homogeneous resolution of $P/I$ by free $P$-modules has the form
$$0\to P(-2n-d+2)\to P(-n-d+2)^{\beta_{d-1}}\to \dots \to P(-n-1)^{\beta_2}\to   P(-n)^{\beta_1} \to P,$$
with $$  \beta_i =\frac{2n+d-2}{n+i-1}\binom{n+d-2}{i-1}\binom{n+d-i-2}{n-1},$$for $1\le i\le d-1$,
\item all of the minimal generators of $I$ have degree  $n$  and the socle of $P/I$ has degree $2n-2$,
\item $[I]_{n-1}=0$ and $[P/I]_{2n-1}=0$,
\item $\phi\in D_{2n-2}^{\pmb k}(U^*)$ and the homomorphism $\pmb p_{n-1}^\phi\:\operatorname {Sym}_{n-1}^{\pmb k}U\to D_{n-1}^{\pmb k}(U^*)$ of Definition {\rm\tref{T}} is an isomorphism, and 
\item $\phi\in D_{2n-2}^{\pmb k}(U^*)$ and the $\binom{d+n-2}{d-1} \times \binom{d+n-2}{d-1}$ matrix $T_\phi$ of Definition {\rm\tref{T}} is invertible.
\endroster
\endproclaim

\remark{Remark} The Betti numbers $\beta_i$, from (2), also are equal to 
 $$\beta_i=\binom {d+n-1}{n+i-1}\binom{i+n-2}{i-1}-\binom{d+n-2}{i-1}\binom{d-i+n-2}{d-i},$$for $1\le i\le d-1$. The present formulation is given in Theorem \tref{ek-k-m}. A quick calculation shows that the two formulations are equal. \endremark

\demo{Proof of Proposition {\rm \tref{J18}}} $(1)\Rightarrow (2)$. The fact that $P/I$ has a linear resolution is part of the definition of $\Bbb I^{[d]}_n(\pmb k)$. The  Betti number $\beta_i$ comes from the Herzog-K\"uhl formula \cite{\rref{HK}}.

\medskip \flushpar $(2)\Rightarrow (1).$ One can read from (2) that $I$ is generated by forms of degree $n$ and that $P/I$ has a linear resolution.  Thus, $I$ is in $\Bbb I_n^{[d]}(\pmb k)$, as defined in (\tref{In}).

\medskip \flushpar $(2)\Rightarrow (3).$ One can read from (2) that $I$ is generated by forms of degree $n$. Furthermore, the socle degree of $P/I$ is $b+a(P)$,  where $a(P)$ is the $a$-invariant of $P$ and $b$ is the ``back twist'' in the $P$-free resolution of $P/I$. (This observation is well-known. It amounts to computing $\operatorname{Tor}_d^P(P/I,\pmb k)$ in each component;
see, for example, the proof of \cite{\rref{KV},~Prop.~1.5}.) In the present situation, $b=2n+d-2$ and $a(P)=-d$; so the socle degree of $P/I$ is $2n-2$. 

\medskip \flushpar $(3)\Rightarrow (4)$. Statement  (3) asserts that $[I]_i=0$ for $i\le n-1$ and $[P/I]_i=0$ for $2n-1\le i$.

\medskip \flushpar $(4)\Rightarrow (2)$. Let $\Bbb A:\quad 0\to A_d\to \dots \to A_2\to A_1\to P\to P/I\to 0$ be a minimal homogeneous resolution of $P/I$ by free $P$-modules. Write $A_\ell=\bigoplus_{0\le i}P(-i)^{\beta_{\ell,i}}$. This resolution is self-dual; so, in particular, $A_d=P(-b)$ for some twist $b$. As above, it follows that  the socle degree of $P/I$  is equal to $b-d$. In the present situation, the socle degree of $P/I$ is $2n-2-q$ for some non-negative integer $q$; so $b=2n-2+d-q$. The hypothesis $I_{n-1}=0$ ensures that 
$$\text{$\beta_{1,i}=0$, whenever $i\le n-1$}.\tag\tnum{fst}$$ Duality forces
$$\beta_{d-1,j}=0 \text{ for $n-1+d-q\le j$}.\tag\tnum{blop}$$
The fact that the resolution $\Bbb A$ is  minimal  ensures that 
$$\text{if $\beta_{\ell,j}\neq 0$ for some  $1\le \ell$, then there exists $i\le j-1$ with $\beta_{\ell-1,i}\neq 0$.}\tag\tnum{abv2}$$ One may iterate the  idea of (\tref{abv2}) to see that 
$$\text{if $\beta_{d-1,j}\neq 0$, then there exists $i\le j+2-d$   with $\beta_{1,i}\neq 0$.}\tag\tnum{abv3}$$
Suppose $\beta_{d-1,j}\neq 0$. On the one hand, $j\le n-2+d-q$ by (\tref{blop})
and on the other hand, according to (\tref{abv3}), there exits $i$ with $i\le j+2-d$ and $\beta_{1,i}\neq 0$. Apply (\tref{fst}) to see that $n\le i$. 
Thus,
$$n\le i\le j+2-d\le n-q\le n;$$
hence,  $q=0$, $i=n$, and $j=n+d-2$. 
At this point, we have shown that 
$$A_d=P(-(2n-2+d)),\ \ A_{d-1}=P(-(n+d-2))^{\beta_{d-1,n+d-2}},\ \ \text{and}\ \ A_1=P(-n)^{\beta_{1,n}},$$with
$\beta_{d-1,n+d-2}=\beta_{1,n}$.

We study $\beta_{\ell,j}$ for $2\le \ell\le d-2$. Apply (\tref{abv2}) to $\Bbb A$ repeatedly to see that $\beta_{\ell,j}=0$ for $j\le n+\ell-2$. The ideal $I$ is Gorenstein and $\Bbb A$ is a minimal resolution of $P/I$; therefore, the complex $\Bbb A^*(-2n-d+2)$ is isomorphic to $\Bbb A$,
$$\beta_{a,b}=\beta_{a',b'}\quad \text{whenever $a+a'=d$ and $b+b'=2n+d-2$},$$
 and  $\beta_{\ell,j}$ is non-zero only for $j=n+\ell-1$.   One computes the Betti numbers of $\Bbb A$ by using  the Herzog-K\"uhl formula \cite{\rref{HK}} for the Betti numbers in a pure resolution.

\medskip\flushpar $(3)\Rightarrow (5)$. The ring $P/I$ is graded and Gorenstein with socle degree $2n-2$; hence, multiplication
gives a perfect pairing $[P/I]_{n-1}\otimes_{\pmb k} [P/I]_{n-1} \to [P/I]_{2n-2}$. Statement (3) gives $[I]_{n-1}=0$; hence, $[P/I]_{n-1}$ is equal to $[P]_{n-1}=\operatorname {Sym}_{n-1}^{\pmb k}U$. Statement (3) also gives that the socle degree of $P/I$ is $2n-2$; thus,   Theorem \tref{MT} yields that the function $[P/I]_{2n-2} \to \pmb k$, which sends  the class of $\overline{u_{2n-2}}$ in $[P/I]_{2n-2}$ to $\phi(u_{2n-2})$ in $\pmb k$, is an isomorphism. (We viewed $u_{2n-2}$ as an element of $\operatorname {Sym}_{2n-2}^{\pmb k}U$ and $\overline{u_{2n-2}}$ as the image of $u_{2n-2}$ in $[P/I]_{2n-2}$.)
Thus, statement (3) ensures that the homomorphism
$$\operatorname {Sym}_{n-1}^{\pmb k}U\otimes_{\pmb k} \operatorname {Sym}_{n-1} ^{\pmb k}U\to \pmb k,$$ which is given by
$u_{n-1}\otimes u_{n-1}'\mapsto \phi(u_{n-1}u_{n-1}')$, is a perfect pairing. It follows that  the homomorphism $\operatorname {Sym}_{n-1}^{\pmb k}U\to D_{n-1}^{\pmb k}(U^*)$, which is given by $u_{n-1}\mapsto u_{n-1}(\phi)$, is an isomorphism.

\medskip\flushpar $(5)\Rightarrow (4)$. Apply Theorem \tref{MT} to $\phi\in D_{2n-2}^{\pmb k}(U^*)$ to see that the socle degree of $P/I$ is $2n-2$; and therefore, $[P/I]_{2n-1}=0$. If $v_{n-1}\in \operatorname {Sym}_{n-1}^{\pmb k}U$ is in $I$, then $v_{n-1}(\phi)$ is the zero element of $D_{n-1}^{\pmb k}(U^*)$; thus statement (5) guarantees that  $v_{n-1}$ is zero in $\operatorname {Sym}^{\pmb k}_{n-1}U$ and $[I]_{n-1}=0$. 

 \medskip\flushpar $(5)\Leftrightarrow (6)$. It is clear that these statements are equivalent because $T_\phi$ is the matrix for  $\pmb p_{n-1}^\phi$; see Remark \tref{R1}.1.  \qed\enddemo
 
 The following statement is an immediate consequence of Proposition \tref{J18}; no further proof is necessary.
\proclaim{Corollary \tnum{open}} Let $U$ be a $d$-dimensional vector space over the field $\pmb k$ and $n$ be a positive integer.  If $\phi$ is a homogeneous element of $D_{\bullet}^{\pmb k}(U^*)$, then $\operatorname{ann} \phi \in  \Bbb  I_n^{[d]}(\pmb k)$ if and only if  $\deg \phi =2n-2$   and $\det T_{\phi}\neq 0$ .  
In particular, the open subset 
$$ O=  D_{2n-2}^{\pmb k}(U^*)\setminus \{\text{the variety defined by $\det T_{\phi}=0$}\}$$ of $D_{2n-2}^{\pmb k}(U^*)$  parameterizes  
 $ \Bbb  I_n^{[d]}(\pmb k)$.
 \endproclaim 

\remark{\bf Remark \tnum{Reta}}   We emphasize that, in the language of Corollary \tref{open}, $D_{2n-2}^{\pmb k}(U^*)$ is a vector space of dimension $\nu=\binom{2n+d-3}{d-1}$ and, once a basis $\{b_i\}$ is chosen for this vector space, then $D_{2n+d-3}^{\pmb k}(U^*)$ can be identified with affine $\nu$-space: the point $(a_1,\dots,a_{\nu})$ in affine space corresponds to the element $\sum a_ib_i$ of $D_{2n-2}^{\pmb k}(U^*)$. Under this identification, $\det T_{\phi}$ corresponds to a homogeneous polynomial of degree $\binom{n+d-2}{d-1}$ in the coordinates of affine $\nu$-space; see Remark \tref{R1}.2. Thus, Corollary \tref{open} parameterizes $\Bbb I_n^{[d]}(\pmb k)$ using an open subset of $\nu$-space. The open subset is the complement of a hypersurface. \endremark

One further consequence of Proposition \tref{J18} is the following generalization of the set of ideals $\Bbb I_n^{[d]}(\pmb k)$.

\definition{Definition \tnum{inr}} Let $R_0$ be a commutative Noetherian ring, $U$ be a non-zero free $R_0$-module of finite rank and $n$ be a positive integer. Define $\Bbb I_n(R_0,U)$ to be the following set of ideals in $P=\operatorname{Sym}_\bullet^{R_0}U$:
$$\Bbb I_n(R_0,U)=\{\operatorname{ann}\phi\mid \phi\in D_{2n-2}^{R_0}(U^*) \text{ and } \det T_\phi \text{ is a unit in $R_0$}\}.
$$\enddefinition

We see from Proposition \tref{J18} that if $U$ is a vector space of dimension $d$ over a field $\pmb  k$, then the sets of ideals $\Bbb I_n(R_0,U)$ and $\Bbb I^{[d]}_n(\pmb k)$ are equal. Corollary \tref{EK-K-2}, which is our solution to Project \tref{EK1}, is phrased in terms of $\Bbb I_n(R_0,U)$.  

\SectionNumber=\Lpq\tNumber=1
\heading Section \number\SectionNumber. \quad The complexes $\Bbb L(\Psi,n)$ and $\Bbb K(\Psi,n)$.
\endheading

The following data is in effect throughout this section.
\definition{Data \tnum{data23}}
Let $V$ be a non-zero free module of rank $d$ over the commutative Noetherian ring $R$.\enddefinition

Let $\Psi\: V\to R$ be an $R$-module homomorphism and $n$ be a fixed positive integer. In Theorem \tref{BE} we describe a complexes $\Bbb L(\Psi,n)$ and $\Bbb K(\Psi,n)$, dual to one another, so that $\Bbb L(\Psi,n)$ resolves $R/(\operatorname{im} \Psi)^n$ whenever the image of $\Psi$ has grade   $d$.  Our description of these complexes is very explicit and coordinate free. We use our explicit descriptions in the proof of Theorem \tref{EK-K'}. Of course, the complexes of Buchsbaum and Eisenbud in \cite{\rref{BE75}} resolve $R/(\operatorname{im} \Psi)^n$; so, in some sense, our complex $\Bbb L(\Psi,n)$  ``is in'' \cite{\rref{BE75}}. Our proof of the exactness of $\Bbb L(\Psi,n)$ and $\Bbb K(\Psi,n)$  is self contained, very explicit, and, as far as we can tell, different from the proof found in \cite{\rref{BE75}}. 

Retain Data \tref{data23} and let $a$ and $b$ be integers. Define the $R$-module homomorphisms \vphantom{\tnum{kappa}}
$${\alignedat1&\kappa_{a,b}^V\:\tsize{\bigwedge}_R^aV\otimes_R \operatorname {Sym}^R_bV\to \tsize{\bigwedge}_R^{a-1}V\otimes_R \operatorname {Sym}^R_{b+1}V
\quad \text{and}\\\vspace{5pt} &\eta_{a,b}^V\:\tsize{\bigwedge}_R^aV\otimes_R D^R_b(V^*)\to \tsize{\bigwedge}_R^{a-1}V\otimes_R D^R_{b-1}(V^*)\endalignedat}\tag{\tref{kappa}}$$ to be the compositions
$$\tsize{\bigwedge}_R^aV\otimes_R \operatorname {Sym}^R_bV@> \Delta\otimes 1 >>\tsize{\bigwedge}_R^{a-1}V\otimes_R V\otimes_R\operatorname {Sym}^R_{b}V@> 1\otimes \text{mult} >>\tsize{\bigwedge}_R^{a-1}V\otimes_R \operatorname {Sym}^R_{b+1}V \quad\text{and}$$
$$\eightpoint \tsize{\bigwedge}_R^aV\otimes_R D^R_b(V^*)@> \Delta\otimes 1 >>\tsize{\bigwedge}_R^{a-1}V\otimes_RV\otimes_R D^R_{b}(V^*)@> 1\otimes \text{module-action} >> \tsize{\bigwedge}_R^{a-1}V\otimes_R D^R_{b-1}(V^*),$$ respectively; and define the $R$-modules
$$L^R_{a,b}(V)=\ker \kappa_{a,b}^V\quad \text{and}\quad  K^R_{a,b}(V)=\ker \eta_{a,b}^V.$$ (In the future, we will often write 
$\kappa$ and $\eta$ in place of $\kappa_{a,b}^V$ and $\eta_{a,b}^V$.)
The $R$-modules $L^R_{a,b}(V)$ and $K^R_{a,b}(V)$  have been used by many authors in many contexts. In particular, they are studied extensively in \cite{\rref{BE75}}; although our indexing conventions are different than the conventions of \cite{\rref{BE75}}; that is,
$$\text{the module we call $L^R_{a,b}(V)$ is called ${L_b}^{a+1}(V)$ in \cite{\rref{BE75}}.}$$
The complex 
$$\split 0\to \tsize{\bigwedge}_R^dV\otimes_R \operatorname {Sym}^R_{b-d}V@>\kappa^V_{d,b-d} >> \tsize{\bigwedge}_R^{d-1}V\otimes_R \operatorname {Sym}^R_{b-d+1}V@>\kappa^V_{d-1,b-d+1} >>&\\\vspace{5pt}  
\cdots @>\kappa^V_{2,b-2} >> \tsize{\bigwedge}_R^1V\otimes_R \operatorname {Sym}^R_{b-1}V@>\kappa^V_{1,b-1}>>\tsize{\bigwedge}_R^0V\otimes_R \operatorname {Sym}^R_{b}V\to 0,&\endsplit$$ which is a homogeneous strand of an acyclic  Koszul complex, is split exact for all positive integers $b$; hence, $L^R_{a,b}(V)$ is a projective   $R$-module. In fact, $L^R_{a,b}(V)$ is a free $R$-module of rank 
$$\operatorname{rank}_R L^R_{a,b}(V)=\binom{d+b-1}{a+b}\binom{a+b-1}a;\tag\tnum{rank}$$see \cite{\rref{BE75}, Prop. 2.5}.
The perfect pairing
$$\left(\vphantom{E^{E^E}}{\tsize \bigwedge}_R^aV\otimes_R\operatorname {Sym}_b^RV\right)\otimes_R \left({\tsize \bigwedge}_R^{d-a}V\otimes_RD_b^R(V^*)\right)\to {\tsize \bigwedge}_R^dV,$$ which is given by
$$(\theta_a\otimes v_b)\otimes (\theta_{d-a}\otimes w_b) \mapsto v_b(w_b)\cdot \theta_a\wedge \theta_{d-a},$$
induces a perfect pairing 
$$L^R_{a,b}(V)\otimes K^R_{d-a-1,b-1}(V) \to {\tsize \bigwedge}^d_R(V).\tag\tnum{p207}$$
The indices in (\tref{p207}) are correct because the dual of the presentation
$${\tsize \bigwedge}_R^{a+2}V\otimes_R\operatorname {Sym}_{b-2}^RV@> \kappa^V_{a+2,b-2} >>{\tsize \bigwedge}_R^{a+1}V\otimes_R\operatorname {Sym}_{b-1}^RV@> \kappa^V_{a+1,b-1} >> L^R_{a,b}(V)\to 0$$
is
$$0 \to L^R_{a,b}(V)^*@>>> ({\tsize \bigwedge}_R^{a+1}V\otimes_R\operatorname {Sym}_{b-1}^RV)^*@>>> ({\tsize \bigwedge}_R^{a+2}V\otimes_R\operatorname {Sym}_{b-2}^RV)^*,$$ which is isomorphic to 

$$  \split& 0 \to K^R_{d-a-1,b-1}(V)\otimes_R {\tsize \bigwedge}_R^dV^* \to \\\vspace{5pt} & ({\tsize \bigwedge}_R^{d-a-1}V\otimes_R D_{b-1}^R(V^*))\otimes_R{\tsize \bigwedge}_R^dV^*  \to  ({\tsize \bigwedge}_R^{d-a-2}V\otimes_R D_{b-2}^R(V^*))\otimes_R{ \tsize \bigwedge}_R^dV^*.\endsplit $$

\remark{\bf Remark \tnum{schur}}The modules $L^R_{a,b}(V)$ and $K^R_{a,b}(V)$ may also be thought of as the Schur modules $L_\lambda(V)$ and Weyl modules $K_\lambda(V^*)$ which correspond to certain hooks $\lambda$. We use the notation of Examples 2.1.3.h and 2.1.17.h in Weyman \cite{\rref{W}} to see that 
the module we call $L^R_{a,b}(V)$ is also the Schur module $L_{(a+1,1^{b-1})}(V)$ and the module we call  $K^R_{a,b}(V)$ is also $K^R_{(b+1,1^{d-a-1})}(V^*)\otimes {\tsize \bigwedge}^d V$, where $K^R_{(b+1,1^{d-a-1})}(V^*)$ is a Weyl module. We pursue this line of reasoning in Section \number\exmpls.
\endremark

\definition{Definition \tnum{long}} Let $V$ be a non-zero free module of rank $d$ over the commutative Noetherian ring $R$, $n$ be a positive integer,  and ${\Psi\: V\to R}$ be an $R$-module homomorphism. We define the complexes 
$$\eightpoint \split&\Bbb L(\Psi,n)\:  0\to L^R_{d-1,n}(V)@>\operatorname{Kos}^{\Psi}\otimes 1 >> L^R_{d-2,n}(V)@>\operatorname{Kos}^{\Psi}\otimes 1 >>\dots @>\operatorname{Kos}^{\Psi}\otimes 1 >>
L^R_{0,n}(V) @> \widehat{\Psi}>> R,\ \text{and}\\\vspace{5pt}
& \Bbb K(\Psi,n)\:  0\to {\tsize \bigwedge}_R^d V@> >>  K^R_{d-1,n-1}(V)@>\operatorname{Kos}^{\Psi}\otimes 1 >> K^R_{d-2,n-1}(V)@>\operatorname{Kos}^{\Psi}\otimes 1 >>\cdots \\&\phantom{ \Bbb K(\Psi,n)\:  } \cdots@>\operatorname{Kos}^{\Psi}\otimes 1 >>
K^R_{0,n-1}(V)\endsplit $$
which appear in Theorems \tref{BE} and \tref{EK-K'}. 

The ordinary Koszul complex associated to $\Psi$ is 
$$0\to {\tsize \bigwedge}_R^dV@> \operatorname{Kos}^{\Psi}>> {\tsize \bigwedge}_R^{d-1}V@> \operatorname{Kos}^{\Psi}>> \dots
@> \operatorname{Kos}^{\Psi}>>{\tsize \bigwedge}_R^2V@> \operatorname{Kos}^{\Psi}>>{\tsize \bigwedge}_R^1V@> \operatorname{Kos}^{\Psi}>>R,\tag\tnum{oKc}$$ where for each index $i$, $\operatorname{Kos}^\Psi$ is the composition 
$${\tsize \bigwedge}_R^iV@> \Delta >> V\otimes_R {\tsize \bigwedge}_R^{i-1}V@> \Psi\otimes 1 >> R\otimes_R {\tsize \bigwedge}_R^{i-1}V={\tsize \bigwedge}_R^{i-1}V.$$ (The co-multiplication map $\Delta$ is discussed in (\tref{delta}).) The maps of (\tref{kappa}) combine with the maps $\operatorname{Kos}^\Psi$ to form  double complexes
$$ \eightpoint \CD 
@. \vdots @. \vdots\\
@.   @V  \kappa  VV@V  \kappa  VV\\
\cdots @> \operatorname{Kos}^\Psi\otimes 1 >> \tsize{\bigwedge}_R^aV\otimes_R \operatorname {Sym}^R_bV @> \operatorname{Kos}^\Psi\otimes 1 >> \tsize{\bigwedge}_R^{a-1}V\otimes_R \operatorname {Sym}^R_bV@> \operatorname{Kos}^\Psi\otimes 1 >>\cdots   \\
@. @V  \kappa  VV @V  \kappa  VV\\
\cdots@> \operatorname{Kos}^\Psi\otimes 1 >> \tsize{\bigwedge}_R^{a-1}V\otimes_R \operatorname {Sym}^R_{b+1}V@> \operatorname{Kos}^\Psi\otimes 1 >>  \tsize{\bigwedge}_R^{a-2}V\otimes_R 
\operatorname {Sym}^R_{b+1}V@> \operatorname{Kos}^\Psi\otimes 1 >>\cdots \\
@. @V  \kappa  VV@V  \kappa  VV\\
@.  \vdots @. \vdots\\
\endCD\tag \tnum{DC1}$$
and 
$$ \eightpoint \CD 
@. \vdots @. \vdots\\
@.   @V  \eta  VV@V  \eta  VV\\
\cdots @> \operatorname{Kos}^\Psi\otimes 1 >> \tsize{\bigwedge}_R^aV\otimes_R D^R_b(V^*) @> \operatorname{Kos}^\Psi\otimes 1 >> \tsize{\bigwedge}_R^{a-1}V\otimes_R D^R_b(V^*)@> \operatorname{Kos}^\Psi\otimes 1 >>\cdots   \\
@. @V  \eta  VV @V  \eta  VV\\
\cdots@> \operatorname{Kos}^\Psi\otimes 1 >> \tsize{\bigwedge}_R^{a-1}V\otimes_R D^R_{b-1}(V^*)@> \operatorname{Kos}^\Psi\otimes 1 >>  \tsize{\bigwedge}_R^{a-2}V\otimes_R 
D^R_{b-1}(V^*)@> \operatorname{Kos}^\Psi\otimes 1 >>\cdots \\
@. @V  \eta  VV@V  \eta  VV\\
@.  \vdots @. \phantom{.}\vdots.\\
\endCD\tag \tnum{DC2}$$
Most of the complex $\Bbb L(\Psi,n)$ is induced by the double complex (\tref{DC1}). (Keep in mind that $L^R_{d,n}(V)$, which is equal to $\ker (\kappa\: \bigwedge_R^dV\otimes_R \operatorname {Sym}_n^R V\to \bigwedge^{d-1}_RV\otimes_R \operatorname {Sym}^R_{n+1} V)$, is zero.) It remains to describe the right-most map from $L_{0,n}^R(V)=\operatorname {Sym}_n^R(V)$ to $R$.
According to the definition of  symmetric algebra, the $R$-module homomorphism $\Psi\:V\to R$ induces an $R$-algebra homomorphism $\operatorname {Sym}_{\bullet}^R(V)\to R$. We denote this algebra homomorphism by $\widehat{\Psi}$. Most of the complex $\Bbb K(\Psi,n)$ is induced by the double complex (\tref{DC2}). It remains to describe the left most map:\phantom{\tnum{lmm}}
$$\eightpoint {\alignedat 1 &{\tsize \bigwedge}_R^d V={\tsize \bigwedge}_R^d V\otimes_RR@> 1\otimes_R \operatorname{ev}_n^*>> {\tsize \bigwedge}_R^d V\otimes_R\operatorname {Sym}^R_nV\otimes_R D_n^R(V^*)@> 1\otimes_R \widehat{\Psi} \otimes_R 1 >>\\\vspace{5pt}&\hskip10pt {\tsize \bigwedge}_R^d V\otimes_RR\otimes_R D_n^R(V^*)
={\tsize \bigwedge}_R^d V\otimes_R D_n^R(V^*)@> \eta>\simeq> K^R_{d-1,n-1}(V),\endalignedat} \tag{\tref{lmm}}$$where $\operatorname{ev}_n\: D_n^R(V^*)\otimes_R \operatorname {Sym}^R_nV\to R$ is the natural evaluation map, see (\tref{ev*}). \enddefinition

\example{Example \tnum{E2}} Retain the notation and hypotheses of Definition \tref{long}. The complex $\Bbb K(\Psi,1)$ is equal to the ordinary Koszul complex (\tref{oKc}) and the complex $\Bbb L(\Psi,1)$ is isomorphic to the ordinary Koszul complex (\tref{oKc}) by way of the isomorphism
$$\CD
0@>>>{\tsize \bigwedge}_R^dV@>>> \dots @>>> {\tsize \bigwedge}_R^1V@>\operatorname{Kos}^{\Psi}>> {\tsize \bigwedge}_R^0V\\
@. @V \simeq V \kappa V@. @V \simeq V \kappa V @V\simeq V = V\\
0@>>>L_{d-1,1}^R(V)@>>> \dots @>>> L_{0,1}^R(V)@>\widehat{\Psi}>> R.\endCD$$Complexes ``$\Bbb K(\Psi,0)$'' and ``$\Bbb L(\Psi,0)$'' have not been defined. The recipe of Definition \tref{long} would produce  $$0\to {\tsize\bigwedge}^d_RV\to 0$$ for $\Bbb K (\Psi,0)$ and 
$$0\to  {\tsize \bigwedge}_R^{d-1}V@> \operatorname{Kos}^{\Psi}>> \dots
@> \operatorname{Kos}^{\Psi}>>{\tsize \bigwedge}_R^2V@> \operatorname{Kos}^{\Psi}>>{\tsize \bigwedge}_R^1V@> \operatorname{Kos}^{\Psi}>>{\tsize \bigwedge}_R^0V@> 1>>R$$for $\Bbb L(\Psi,0)$; neither of these objects is very satisfactory.  
\endexample

\proclaim{Theorem \tnum{BE}} {\bf (Buchsbaum-Eisenbud \cite{\rref{BE75},Thm. 3.1})} Let $V$ be a non-zero free module of rank $d$ over the commutative Noetherian ring $R$, $n$ be a positive integer,  ${\Psi\: V\to R}$ be an $R$-module homomorphism, and $J$ be the image of $\Psi$. Let $\Bbb L(\Psi,n)$ and $\Bbb K(\Psi,n)$ be the complexes of Definition {\rm \tref{long}}. 
Then, the following statements hold.  
\roster
\item The complexes $\Bbb K(\Psi,n)$ and $[\operatorname{Hom}_R(\Bbb L(\Psi,n),R)\otimes_R\bigwedge^dV][-d]$
are  isomorphic,
where ``$[-d]$'' describes a shift in homological degree.

\item If the ideal  $J$ has grade at least $d$, then  $\Bbb L(\Psi,n)$ is a resolution of    $R/J^n$ by free $R$-modules and  $\Bbb  K(\Psi,n)$ is a resolution of $\operatorname{Ext}_R^d(R/J^n,R)$ by free $R$-modules,

\item If $R$ is a graded ring and $\Psi$ is homogeneous homomorphism with $V$ equal to $R(-1)^d$, then $\Bbb L(\Psi,n)$ is the homogeneous linear complex 
$$0\to R(-n-d+1)^{\beta_d}\to R(-n-d+2)^{\beta_{d-1}}\to \dots \to R(-n-1)^{\beta_2}\to R(-n)^{\beta_1}\to R,$$with 
$\beta_i=\binom{n+d-1}{n+i-1}\binom{n+i-2}{i-1}$, for $1\le i\le d$. 
\endroster \endproclaim

\topinsert $$\smallmatrix
&&&\cdots&&&&& &   &0&&&&0&&0&&0\\
&&&&&&&& &   &\downarrow&&&&\downarrow&&\downarrow&&\downarrow\\
0&\to&&\cdots&&&&&0&\to&K_{d,n-1}&\to&.&\to&K_{2,n-1}&\to&K_{1,n-1}&\to&K_{0,n-1}\\
&&&&&&&& &   &\downarrow&&&&\downarrow&&\downarrow&&\downarrow\\
0&\to&&\cdots&&&&&0&\to&D_{d,n-1}&\to&.&\to&D_{2,n-1}&\to&D_{1,n-1}&\to&D_{0,n-1}\\
&&&&&&&& &   &\downarrow&&&&\downarrow&&\downarrow&&\downarrow\\
0&\to&&\cdots&&&0&\to&D_{d,n-2}&\to&D_{d-1,n-2}&\to&.&\to&D_{1,n-2}&\to&D_{0,n-2}&\to&0\\
&&&&&&&\pmb .&&&&&&\pmb .&&\pmb .\\
&&&&&&\pmb .&&&&&&\pmb .&&\pmb .\\
&&&&&\pmb .&&&&&&\pmb .&&\pmb .\\
0 &\to &0&\to&D_{d,1}&\to&.&.&.&\to&D_{1,1}&\to&D_{0,1}&\to&0&&\dots&&0\\
&&\downarrow&&\downarrow&&&&&&\downarrow&&\downarrow&&&&&&\downarrow\\
0&\to&D_{d,0}&\to&D_{d-1,0}&\to&.&.&.&\to&D_{0,0}&\to&0&\to&&&\dots&&0\\
&&\downarrow&&\downarrow&&&&&&\downarrow&&\downarrow&&&&&&\downarrow\\
&&0&&0&&&&&&0&&0&&&&\dots&&0\\
\endsmallmatrix$$
{\bf Figure \tnum{Fig1}.} {\smc The double complex $\Bbb D$ which is used in the proof of Theorem \tref{BE} to show that $\Bbb K(\Psi,n)$ is a resolution for $1\le n$.  The modules   $D_{a,b}$ and $K_{a,b}$ represent $\bigwedge^a_RV\otimes_RD^R_b(V^*)$ and    $K^R_{a,b}(V)$, respectively.  The maps are given in (\tref{DC2}).}
\endinsert

\topinsert $$\smallmatrix
&&&\cdots&&&0&&0 &   &0&&0&&0&&0&&0&&0\\
&&&&&&\downarrow&&\downarrow&& \downarrow&   &\downarrow&&\downarrow&&\downarrow&&\downarrow&&\downarrow\\
0&\to&&\cdots&&&0&&0&\to&E_{d,0}&\to&E_{d-1,0}&\to&.&\to&E_{1,0}&\to&E_{0,0}&\to&0\\
&&&&&&\downarrow&&\downarrow &   &\downarrow&&\downarrow&&\downarrow&&\downarrow&&\downarrow&&\downarrow\\
0&\to&&\cdots&&&0&\to&E_{d,1}&\to&E_{d-1,1}&\to&.&\to&E_{1,1}&\to&E_{0,1}&\to&0&\to& 0\\
&&&&&&&\pmb .&&&&&&\pmb .&&\pmb .\\
&&&&&&\pmb .&&&&&&\pmb .&&\pmb .\\
&&&&&\pmb .&&&&&&\pmb .&&\pmb .\\
0 &\to&0&\to &E_{d,n-2}&\to&
.&\to&E_{2,n-2}&\to&E_{1,n-2}&\to&E_{0,n-2}&\to&\dots&&&&&&0\\
&&\downarrow&&\downarrow&&&&\downarrow&&\downarrow&&\downarrow&&&&&&&&\downarrow\\
0 &\to &E_{d,n-1}&\to&E_{d-1,n-1}&\to
&.&\to&E_{1,n-1}&\to&E_{0,n-1}&\to&0&&\dots&&&&&&0\\
&&\downarrow&&\downarrow&&&&\downarrow&&\downarrow&&\downarrow&&&&&&&&\downarrow\\
0&\to&L_{d-1,n}&\to&L_{d-2,n}&\to
&.&\to&L_{0,n}&\to&0&\to&0&&\dots&&&&&&0\\
&&\downarrow&&\downarrow&&&&\downarrow&&\downarrow&&\downarrow&&&&&&&&\downarrow\\
&&0&&0&&&&0&&0&&0&&\dots&&&&&&0\\
\endsmallmatrix$$ 
{\bf Figure \tnum{Fig2}.} {\smc The double complex $\Bbb E$ which is used in the proof of Theorem \tref{BE} to show that $\Bbb L(\Psi,n)$ is a resolution for $1\le n$.  The modules   $E_{a,b}$ and $L_{a,b}$ represent $\bigwedge^a_RV\otimes_R\operatorname {Sym}^R_bV$ and    $L^R_{a,b}(V)$, respectively.  The maps are given in (\tref{DC1}).}
\endinsert

\demo{Proof} Most of the proof of assertion (1) is contained in (\tref{p207}). One can use (\tref{rank}) to prove (3); or one can prove (2) and then appeal to the Herzog-K\"uhl formula \cite{\rref{HK}}. We focus on the proof of (2). Assume that  $J$ has grade $d$. We first show that $\Bbb K(\Psi, n)$ is a resolution. 
Let $\Bbb K'$ be the sub-complex of $\Bbb K(\Psi,n)$ which is obtained by deleting the left-most non-trivial module $\bigwedge^d_RV$. We prove that
$$\operatorname{H}_j(\Bbb K')=\cases 0&\text{if $1\le j\le d-2$}\\\bigwedge^d_RV&\text{if $j=d-1$}\endcases\tag\tnum{light}$$
  and that the image of the map $\bigwedge^d_RV\to K^R_{d-1,n-1}$, which is given in (\tref{lmm}), is exactly equal to
$\operatorname{H}_{d-1}(\Bbb K')$. 

Consider the double complex $\Bbb D$ which is obtained from (\tref{DC2}) by keeping the modules $\bigwedge^a_RV\otimes_RD^R_b(V^*)$ for $b\le n-1$, and adjoining the row of kernels
$$0\to K^R_{d,n-1}(V)\to \dots\to K^R_{0,n-1}(V)\to 0.$$ Keep in mind that each $K^R_{i,n-1}(V)$ is the kernel at the top of a column of our truncation of  (\tref{DC2}). In other words, the maps
 $$\eightpoint \CD 0\\@VVV \\ K_{a,n-1}\\@V \text{incl} VV \\ \bigwedge^a_RV\otimes_RD^R_{n-1}(V^*)\\@V \eta VV\\ \bigwedge^{a-1}_RV\otimes_RD^R_{n-2}(V^*)\endCD$$form an exact sequence for each $a$. We have recorded a picture of $\Bbb D$ in Figure \tref{Fig1}, where  $D_{a,b}$ represents $\bigwedge^a_RV\otimes_RD^R_b(V^*)$ and $K_{a,b}$ represents $K^R_{a,b}(V)$. 
Index $\Bbb D$ so that the total complex $\Bbb T(\Bbb D)$ has the module $\sum_{j=1}^{n-1}D_{0,j}$ in position $0$.

 Each column of $\Bbb D$, except for the left-most non-trivial column, is exact. The hypothesis that $J$ has grade $d$ guarantees that the rows of $\Bbb D$, with the   exception of the top-most non-trivial row, have homology concentrated in the  position of the module $D_{0,*}$. We consider two sub-complexes $\Bbb D'$ and $\Bbb D''$ of $\Bbb D$. Let $\Bbb D'$ be the sub-complex of $\Bbb D$ which is obtained by deleting the left-most non-trivial column and 
   $\Bbb D''$ be the sub-complex  which is obtained by deleting the top-most non-trivial row. 
 Observe that $\Bbb D/\Bbb D'$ is the left-most non-trivial column of $\Bbb D$; hence,
$$\Bbb T(\Bbb D/\Bbb D')_j =\cases 0 &\text{if $j\neq d$}\\D_{d,0}=\bigwedge^d_RV&\text{if $j=d$}.\endcases$$
Observe also that $\Bbb D/\Bbb D''$ is the top-most non-trivial row of $\Bbb D$; so, in particular, 
$$\Bbb T(\Bbb D/\Bbb D'')_j= \left.\cases K_{-1,n-1}=0 &\text{for $j=0$}\\
K_{j-1,n-1} &\text{for $1\le j\le d$}\\
K_{d,n-1}=\operatorname{ker}\eta_{d,n-1}=0&\text{for $j=d+1$}\endcases\right\}=\Bbb K'_{j-1}.$$ 
and $\Bbb T(\Bbb D/\Bbb D'')=\Bbb K'[-1]$. In light of (\tref{light}), we show   that $\Bbb K(\Psi, n)$ is a resolution 
by showing that    
$$\operatorname{H}_j(\Bbb T(\Bbb D/\Bbb D''))=\cases 0&\text{if $2\le j\le d-1 $}\\\bigwedge^d_RV& \text{if $j= d $}  \endcases\tag\tnum{light2}$$and showing that the isomorphism $\bigwedge^d_RV\to \operatorname{H}_d(\Bbb T(\Bbb D/\Bbb D''))$ is given by (\tref{lmm}).

Every column of $\Bbb D'$ is split exact; so, the total complex of $\Bbb D'$,   denoted $\Bbb T(\Bbb D')$, is also split exact. The short exact sequence of total complexes
$$0\to \Bbb T(\Bbb D')\to \Bbb T(\Bbb D) \to \Bbb T(\Bbb D/\Bbb D')\to 0$$ yields that $\operatorname{H}_{\bullet}(\Bbb T(\Bbb D))\simeq \operatorname{H}_{\bullet}(\Bbb T(\Bbb D/\Bbb D'))$.
 Thus, $$\operatorname{H}_{j}(\Bbb T(\Bbb D))\simeq \cases 0&\text{if $j\neq d$}\\\bigwedge_R^dV&\text{if $j= d$};\endcases\tag\tnum{hiso}$$
furthermore, the isomorphism of (\tref{hiso}), when $j=d$, is obtained by lifting each element $\theta$ in $\bigwedge^dV=D_{d,0}$ back to a cycle in $\Bbb T(\Bbb D)$. Recall the canonical map 
$$\operatorname{ev}_j^*\:R\to \operatorname {Sym}_j^R(V)\otimes_R D_j^R(V^*),$$ as described in (\tref{ev*}).
Observe that 
$(1\otimes \widehat{\Psi}\otimes 1)(\theta\otimes \operatorname{ev}_j^*(1))$ is an element of $${\tsize\bigwedge}^d_RV\otimes_R D_j^R(V^*)=D_{d,j}.$$ 
It is easy to see that there exists signs $\sigma_j\in \{1,-1\}$ so that $$\sum_{0\le j}\sigma_j\cdot (1\otimes \widehat{\Psi}\otimes 1)(\theta\otimes \operatorname{ev}_j^*(1))$$ is the unique cycle in (\tref{DC2}) which lifts $\theta\in D_{d,0}$. The complex $\Bbb D$ is obtained by truncating the complex (\tref{DC2}) and then adjoining a row of kernels. 
We conclude that the  isomorphism of (\tref{hiso}), when $j=d$, is given by 
 $$\theta\mapsto \sum_{j=0}^{n-1}\sigma_j\cdot (1\otimes \widehat{\Psi}\otimes 1)(\theta\otimes \operatorname{ev}_j^*(1))+\sigma_{n}\cdot
\eta\left(\vphantom{E^{E^E}}(1\otimes \widehat{\Psi}\otimes 1)(\theta\otimes \operatorname{ev}_n^*(1))\right),
\tag\tnum{psto}$$for $\theta\in \bigwedge^d_R(V)$.

The rows of $\Bbb D''$ all have homology concentrated in position $D_{0,*}$; and therefore, the homology of $\Bbb T(\Bbb D'')$ is concentrated in the position   $0$. The long exact sequence of homology which corresponds to the short exact sequence of complexes 
$$0\to \Bbb T(\Bbb D'')\to \Bbb T(\Bbb D) \to \Bbb T(\Bbb D/\Bbb D'')\to 0$$ yields that 
$$\operatorname{H}_j(\Bbb T(\Bbb D))\simeq \operatorname{H}_j(\Bbb T(\Bbb D/\Bbb D''))\quad \text{for $2\le j$}.\tag\tnum{***}$$
We apply (\tref{hiso}) to conclude that 
$$\operatorname{H}_j(\Bbb T(\Bbb D/\Bbb D''))\simeq \cases 
0&\text{if $2\le j\leq d-1$}\\\bigwedge_R^dV&\text{if $j= d$};\endcases $$
furthermore, the 
composition 
$${\tsize \bigwedge}^d_RV@> \text{(\tref{psto})}>\simeq>  \operatorname{H}_d(\Bbb T(\Bbb D)) @> \text{(\tref{***})}>\simeq> \operatorname{H}_d(\Bbb T(\Bbb D/\Bbb D''))$$ 
sends $\theta$ in ${\tsize \bigwedge}^d_RV$ to 
$$\sigma_{n}\cdot
\eta\left(\vphantom{E^{E^E}}(1\otimes \widehat{\Psi}\otimes 1)(\theta\otimes \operatorname{ev}_n^*(1))\right)\in K_{d-1,n-1}=\Bbb T(\Bbb D/\Bbb D'')_{d}.$$
The constant $\sigma_n\in\{1,-1\}$ is irrelevant. The map
$$\theta \mapsto \eta\left(\vphantom{E^{E^E}}(1\otimes \widehat{\Psi}\otimes 1)(\theta\otimes \operatorname{ev}_n^*(1))\right)$$ is exactly the map of (\tref{lmm}). 
We have accomplished both objectives from (\tref{light2}); hence, we have shown that $\Bbb K(\Psi,n)$ is a resolution. 

We use the same style of argument to show that $\Bbb L(\Psi,n)$ is a resolution. Let $\Bbb E$ be the double complex (\tref{DC1}) truncated to include $\bigwedge_R^aV\otimes_R\operatorname {Sym}_b^RV$ for $0\le b\le n-1$ with a row of cokernels adjoined. We have recorded a picture of $\Bbb E$ in Figure \tref{Fig2} where $E_{a,b}$ represents $\bigwedge_R^aV\otimes_R\operatorname {Sym}_b^RV$ and $L_{a,b}$ represents $L^R_{a,b}(V)$. We index $\Bbb E$ so that $\Bbb T(\Bbb E)_0$ is the module $\sum_{b=0}^{n-1}E_{0,b}\oplus L_{0,n}$. Define the sub-complexes $\Bbb E'$ and $\Bbb E''$ of $\Bbb E$ with $\Bbb E'$ equal to the right-most non-trivial column of $\Bbb E$ and $\Bbb E''$ equal to  the bottom-most non-trivial row of $\Bbb E$. Every column of $\Bbb E/\Bbb E'$ is exact; so, $\operatorname{H}_{\bullet}(\Bbb T(\Bbb E/\Bbb E'))=0$ and the long exact sequence of homology associated to the short exact sequences of complexes
$$0\to \Bbb T(\Bbb E')\to \Bbb T(\Bbb E)\to \Bbb T(\Bbb E/\Bbb E')\to 0$$ yields that 
$$\operatorname{H}_{i}(\Bbb T(\Bbb E))\simeq \operatorname{H}_{i}(\Bbb T(\Bbb E'))=\cases R&\text{if $i=0$}\\0&\text{if $i\neq 0$}.\endcases$$ Furthermore, the long exact sequence of homology yields that the isomorphism $\operatorname{H}_{0}(\Bbb T(\Bbb E))\to R$ is induced by the map that sends the cycle 
$$\sum_{b=0}^n z_b\in  \Bbb T(\Bbb E)_0 \text{ to } z_0\in R,\tag\tnum{abv}$$ 
where $z_b\in E_{0,b}$ for $0\le b\le n-1$ and $z_n\in L^R_{0,n}(V)=\bigwedge^0_RV\otimes \operatorname {Sym}_n^RV$.

The rows of $\Bbb E/\Bbb E''$ have homology concentrated in position $E_{0,*}$. It follows that $\operatorname{H}_i(\Bbb T(\Bbb E/\Bbb E''))=0$ for $1\le i$. The long exact sequence of homology associated to the short exact sequences of complexes
$$0\to \Bbb T(\Bbb E'')\to \Bbb T(\Bbb E)\to \Bbb T(\Bbb E/\Bbb E'')\to 0$$ yields that $$\operatorname{H}_i(\Bbb T(\Bbb E''))=0 \quad\text{for $1\le i$}\tag\tnum{zero}$$ and that the natural map $$\operatorname{H}_0(\Bbb T(\Bbb E''))\to \operatorname{H}_0(\Bbb T(\Bbb E))\quad\text{is an injection.}\tag\tnum{nat}$$  The map (\tref{nat}) sends  the class of a cycle $z_n$ from $L^R_{0,n}(V)=\Bbb T(\Bbb E'')_0$ to the class of the corresponding  cycle in $\Bbb T(\Bbb E)$. If $$z_n=\prod_{i=1}^n\ell i\in L^R_{0,n}(V)={\tsize\bigwedge}^0_RV\otimes_R\operatorname {Sym}_n^R V, \text{  with $\ell_i$ in $V$,}$$ then one lift of  $z_n$ to $\Bbb T(\Bbb E)$ is $\sum_{b=0}^n \pm z_b$, with
$$z_b= \widehat{\Psi}(\ell_{b+1}\cdots \ell_{n})\cdot \ell_1\cdots \ell_b\in E_{0,b}\text{ for $0\le b\le n-1$}.$$ Combine (\tref{zero}), (\tref{nat}), and (\tref{abv}) to conclude that the augmented complex $$\Bbb T(\Bbb E'')@> \widehat{\Psi} >> R$$ is a resolution. This resolution is precisely equal to  $\Bbb L(\Psi,n)$. \qed \enddemo

\SectionNumber=\GEN\tNumber=1
\heading Section \number\SectionNumber. \quad The generators.
\endheading 

Observation \tref{gen'} is an important step in the present paper. In Definition \tref{BBBG} we introduce a family of complexes $\Bbb G(*)$. It is immediately clear that these complexes are resolutions. The main theorem of the paper is Theorem \tref{EK-K'} which identifies the zeroth homology  of the resolutions $\Bbb G(*)$. The proof of Theorem \tref{EK-K'} relies on Observation \tref{gen'}.
\definition{Data \tnum{data20'}} Let $R$ be a commutative Noetherian ring, $V$ be a free $R$-module of rank $d$, $n$ be a positive integer, $\Phi$ be an element of $D^R_{2n-2}(V^*)$,  $\operatorname{ann}\Phi$ be the ideal
$$\operatorname{ann}\Phi=\{\theta\in \operatorname {Sym}_{\bullet}^RV\mid \theta(\Phi)=0\in D_\bullet^R(V^*)\}$$
of the $R$-algebra $\operatorname {Sym}_{\bullet}^RV$, and  $$\pmb p_{n-1}^{\Phi}\: \operatorname {Sym}_{n-1}^RV\to D_{n-1}^R(V^*)\tag\tnum{pn-1}$$ be the $R$-module homomorphism defined by $$\pmb p_{n-1}^{\Phi}(v_{n-1})=v_{n-1}(\Phi),$$ for $v_{n-1}\in  \operatorname {Sym}_{n-1}^RV$.
Assume that $\pmb p_{n-1}^{\Phi}$ is an isomorphism. For each integer $i$, let $[\operatorname{ann} \Phi]_i$ represent $\operatorname{ann}\Phi\cap \operatorname {Sym}_{i}^RV$. Define $\sigma_{n-1}\: D_{n-1}^{R}(V^*) \to \operatorname {Sym}_{n-1}^{R}V$ to be the inverse of $p_{n-1}^\Phi$. In particular,  
$$[\sigma_{n-1}(w_{n-1})](\Phi)=w_{n-1}\quad\text{and} \quad
\Phi\left(v_{n-1}[\sigma_{n-1}(w_{n-1})]\right)=v_{n-1}(w_{n-1})\tag\tnum{q'} $$
for all for all $w_{n-1}\in D_{n-1}^{R}(V^*)$ and $v_{n-1}\in \operatorname {Sym}^{R}_{n-1}V$. 
\enddefinition

In Observation \tref{gen'} we identify a set of  generators for $\operatorname{ann} \Phi$.

\proclaim{Observation \tnum{20.20'}} Adopt Data {\rm\tref{data20'}} and let   $\rho$ be an integer with $0\le \rho\le n-1$. Then the following statements hold.
\item{\rm (1)} The $R$-module homomorphism $\Phi\:\operatorname {Sym}_{2n-2}^RV\to R$ induces an isomorphism $$\frac{\operatorname {Sym}_{2n-2}^RV}{[\operatorname{ann} \Phi]_{2n-2}}\simeq R.$$
\item{\rm (2)} The pairing
$$\operatorname {Sym}_{n-1-\rho}^RV\otimes_R \operatorname {Sym}_{n-1+\rho}^RV\to R,\tag\tnum{pair'}$$ which is given by   $$ v_{n-1-\rho}\otimes v_{n-1+\rho}\mapsto (v_{n-1-\rho}v_{n-1+\rho})(\Phi),$$ induces a perfect pairing
$$\operatorname {Sym}_{n-1-\rho}^RV\otimes_R \frac{\operatorname {Sym}_{n-1+\rho}^RV}{[\operatorname{ann} \Phi]_{n-1+\rho}}\to R.$$ 
\item{\rm (3)} If $x_1,\dots,x_d$ is a basis for $V$, then the elements   $$\tsize \{x_1^\rho\sigma_{n-1}({x_1^*}^{(a_1+\rho)}{x_2^*}^{(a_2)}\cdots {x_d^*}^{(a_d)})\mid \sum_i a_i=n-1-\rho,\ 0\le a_i\} \tag\tnum{out'}$$ of $\operatorname {Sym}_{n-1+\rho}^{R}V$
  are   dual to the monomial basis  $$\tsize \{x_1^{a_1}x_2^{a_2}\cdots x_d^{a_d}\mid \sum_i a_i=n-1-\rho\}$$ of $\operatorname {Sym}_{n-1-\rho}^{R}V$
under the pairing {\rm(\tref{pair'})}. 
\endproclaim 

\demo{Proof} Assertion (1) is a special case of (2); (2) is an immediate consequence of (3); and (3) is obvious. \qed \enddemo 

\remark{\bf Remarks \tnum{20.21'}} (1) If $C_{n-1+\rho}$ is the $R$-submodule of $\operatorname {Sym}_{n-1+\rho}^RV$ which is generated by the elements of (\tref{out'}), then the Observation \tref{20.20'} shows that  $C_{n-1+\rho}$ is a free $R$-module with basis (\tref{out'}) and  $\operatorname {Sym}_{n-1+\rho}^RV$ may be decomposed as the direct sum of two $R$-submodules:
$$\operatorname {Sym}_{n-1+\rho}^RV=[\operatorname{ann} \Phi]_{n-1+\rho}\oplus C_{n-1+\rho}.$$ 

\medskip \flushpar(2) If $0\le \rho\le n-1$, then Observation \tref{20.20'} shows that $$\pmb p_{n-1+\rho}^\Phi\: \operatorname {Sym}_{n-1+\rho}^RV\to D_{n-1-\rho}^R(V^*)$$ is surjective. (Recall the definition of $\pmb p_i^\Phi$ from Definition \tref{T}.) The $R$-module $D_{n-1-\rho}^R(V^*)$ is free; so there exists an $R$-module homomorphism $$\sigma_{n-1-\rho}\: D_{n-1-t}^R(V^*)\to \operatorname {Sym}_{n-1+\rho}^RV$$ which is  a splitting map for $\pmb p_{n-1+\rho}^\Phi$; thus, in particular, the composition 
$$D_{n-1-\rho}^R(V^*)@>\sigma_{n-1-\rho}>> \operatorname {Sym}_{n-1+\rho}^RV@>\pmb p_{n-1+\rho}^\Phi>> D_{n-1-\rho}^R(V^*)\tag\tnum{2'}$$ is the identity map on $D_{n-1-\rho}^R(V^*)$. The map $\sigma_{n-1}$ has been previously defined, in Data \tref{data20'}, in a coordinate-free manner. The maps 
$\sigma_{n-1-\rho}$, for $1\le \rho\le n-1$, depend on the choice of a basis.

\medskip \flushpar(3) For each $\rho$ with $0\le \rho\le n-1$, let $\alpha_{n-1-\rho}\:K^R_{1,n-1-\rho}(V)\to \operatorname {Sym}_{n+\rho}^RV$ be the composition 
$$K^R_{1,n-1-\rho}(V) @> 1\otimes \sigma_{n-1-\rho}>> V\otimes_R \operatorname {Sym}_{n-1+\rho}^RV@> \text{multiplication}>>  \operatorname {Sym}_{n+\rho}^RV,$$and let $A[n+\rho]$ be the image of $\alpha_{n-1-\rho}$ in $\operatorname {Sym}_{n+\rho}^RV$. Notice that $A[n+\rho]$ is defined to be an $R$-module. This $R$-module is independent of the choice of coordinates when $\rho=0$, and depends on the choice of coordinate when $1\le \rho\le n-1$. 
\endremark

\proclaim{Observation \tnum{gen'}} Adopt Data {\rm\tref{data20'}} and the notation of Remarks {\rm\tref{20.21'}}. 
The following statements hold.\roster\item If $0\le \rho\le n-1$, then the $R$-module  $A[n+\rho]$ is contained in the $R$-module $[\operatorname{ann} \Phi]_{n+\rho}$.
\item The ideal $\operatorname{ann}\Phi$ of $\operatorname {Sym}_\bullet^R(V)$ is generated by $[\operatorname{ann} \Phi]_n$. 
\item The $R$-modules $A[n]$ and $[\operatorname{ann} \Phi]_n$ are equal.
\item Let $\Psi\:V\to R$ be an $R$-module homomorphism and $\widehat{\Psi}\:\operatorname {Sym}_\bullet^R(V)\to R$ be the $R$-algebra homomorphism induced by $\Psi$. If $I$ and $J$ are the ideals $I=\widehat{\Psi}(\operatorname{ann} \Phi)$ and $J=\widehat{\Psi}(V)$ of $R$, then the following statements hold for all non-negative integers $\rho $,
\itemitem{\rm(a)}$\widehat{\Psi}(\operatorname{ann} \Phi \,\cap\, \operatorname{Sym}_{n+\rho }^RV)$ is equal to the ideal $J^\rho I$ of $R$, and
\itemitem{\rm(b)} a generating set for the ideal $J^\rho I$ may be obtained in a polynomial manner from the images of the maps $$\Psi\:V\to R\quad \text{and}\quad \Phi\:\operatorname{Sym}_{2n-2}^RV\to R.$$
\endroster
\endproclaim
\demo{Proof} We first prove (1). A typical element of $A[n+\rho]$ has the form $\alpha_{n-1-\rho}(\Theta)$, where 
  $\Theta=\sum_i \ell_i\otimes w_i$ is in $K^R_{1,n-1-\rho}$, with $\ell_i\in V$, $w_i\in D_{n-1-\rho}^R(V^*)$,  and $\sum_i \ell_i(w_i)=0$. The map  $\alpha_{n-1-\rho}$ sends $\Theta$ to $\sum_i \ell_i\cdot \sigma_{n-1-\rho}(w_i)$; and therefore, 
 (\tref{2'}) yields that 
$$[\alpha_{n-1-\rho}(\Theta)](\Phi)=\sum_i \ell_i( [\sigma_{n-1-\rho}(w_i)](\Phi))= \sum_i \ell_i(  w_i )=0,$$ 
and (1) is established.

We prove (2) and (3) simultaneously. 
Let $\Cal I$ be the ideal of $\operatorname {Sym}_{\bullet}^RV$ which is generated by  $A[n]$,  and let $[\Cal I]_i$ be the $R$-module $\Cal I\cap \operatorname {Sym}_i^RV$ for each $i$. It is clear that $\Cal I$ is generated by $[\Cal I]_n$; and assertion (1) shows  that $[\Cal I]_n\subseteq [\operatorname{ann} \Phi]_n$; hence, $\Cal I\subseteq \operatorname{ann}\Phi$. We prove that $\Cal I=\operatorname{ann}\Phi$. 

 Both ideals $\Cal I$ and $\operatorname{ann}\Phi$ are homogeneous; so it suffices to prove the equality
$[\Cal I]_{i}=[\operatorname{ann} \Phi]_{i}$  for one degree $i$ at a time. It is clear that $[\Cal I]_{i}$ is zero for $i\le n-1$. On the other hand,  if $v_i$ is an element of $[\operatorname{ann} \Phi]_{i}$ for some $i\le n-1$, and $x$ is a basis element of $V$, then $x^{n-1-i}v_i \in \ker \pmb p_{n-1}^{\Phi}=0$. But $x^{n-1-i}$ is a regular element in $\operatorname {Sym}_{\bullet}^RV$; so $v_i =0$.

We identify a few critical elements of $[\Cal I]_n$. Let $x_1,\dots,x_d$ be a basis for $V$ and let $a_1,\dots, a_d$ be non-negative integers.
Observe  that
$$\cases \phantom{-}x_i\sigma_{n-1}({x_1^*}^{(a_1)}\cdots{x_i^*}^{(a_i+1)}\cdots {x_j^*}^{(a_j)}\cdots {x_n^*}^{(a_n)})\\
-x_j\sigma_{n-1}({x_1^*}^{(a_1)}\cdots{x_i^*}^{(a_i)}\cdots {x_j^*}^{(a_j+1)}\cdots {x_n^*}^{(a_n)})\endcases\tag\tnum{b'}$$   
is an element of $[\Cal I]_n$ for any pair $i\neq j$ when   $\sum_\ell a_\ell=n-2$,  and 
$$x_i\sigma_{n-1}({x_1^*}^{(a_1)}\cdots {x_n^*}^{(a_n)}) \tag\tnum{c'}$$
   is an element of  $[\Cal I]_n$ whenever $a_i=0$ and  $\tsize{\sum_\ell a_\ell=n-1}$. 

 Fix an integer $\rho$ with $0\le \rho\le n-1$. We prove that $[\operatorname{ann} \Phi]_{n-1+\rho}\subseteq [\Cal I]_{n-1+\rho}$.
Recall the $R$-submodule $C_{n-1+\rho}$ of $\operatorname {Sym}_{n-1+\rho}^RV$  and the direct sum decomposition
$$\operatorname {Sym}_{n-1+\rho}^RV=[\operatorname{ann} \Phi]_{n-1+\rho}\oplus C_{n-1+\rho}$$ from Remark \tref{20.21'}.
  Use the fact that $$\sigma_{n-1}\:D_{n-1}^{R}(V^*)\to \operatorname {Sym}_{n-1}^{R}V$$ is an isomorphism of free $R$-modules to see that  $\sigma_{n-1}(D_{n-1}^R(V^*))=\operatorname {Sym}_{n-1}^RV$; and hence, $\operatorname {Sym}_\rho^RV\cdot \sigma_{n-1}(D_{n-1}^R(V^*))=\operatorname {Sym}_{n-1+\rho}^RV$.
We next show  that 
$$\operatorname {Sym}_{n-1+\rho}^RV\subseteq [\Cal I]_{n-1+\rho}+C_{n-1+\rho}.\tag\tnum{crit'}$$
Fix  non-negative integers $ A_1,\dots,A_d$ and $a_1,\dots,a_d$, with    $\sum_iA_i=\rho$ and 
 $\sum_ia_i$ equal to $n-1$. 
Let ``$\equiv$'' mean congruent mod $\Cal I$.
Apply (\tref{b'}) and (\tref{c'}) repeatedly to see that  
$x_1^{A_1}\cdots x_d^{A_d}\sigma_{n-1}({x_1^*}^{(a_1)}\cdots{x_d^*}^{(a_d)})$ is 
$$\eightpoint  \equiv \cases 
x_1^{A_1-a_1}x_2^{A_2+a_1}x_3^{A_3}\cdots x_d^{A_d} \sigma_{n-1}({x_2^*}^{(a_2+a_1)}\cdots{x_d^*}^{(a_d)})\equiv 0 &\text{if $a_1<A_1$}\\\\
x_1^{A_1+a_i}\cdots x_i^{A_i-a_i}\cdots  x_d^{A_d} \sigma_{n-1}({x_1^*}^{(a_1+a_i)}\cdots {x_i^*}^{(0)}\cdots{x_d^*}^{(a_d)})\equiv 0 &\left\{\matrix \format\l\\\text{if $a_i<A_i$, for}\\\text{some $i\neq 1$}\endmatrix\right.\\\\
x_1^\rho\sigma_{n-1}({x_1^*}^{(a_1-A_1+\rho)} {x_2^*}^{(a_2-A_2)}\cdots{x_d^*}^{(a_d-A_d)})\in C_{n-1+\rho}&\left\{\matrix \format\l\\\text{if $A_i\le a_i$}\\\text{for all $i$.}\endmatrix\right.\\
\endcases\tag\tnum{equiv'}$$
Thus,  (\tref{crit'}) holds and    $$\operatorname {Sym}_{n-1+\rho}^RV\subseteq [\Cal I]_{n-1+\rho}+C_{n-1+\rho}\subseteq [\operatorname{ann} \Phi]_{n-1+\rho}\oplus C_{n-1+\rho}=\operatorname {Sym}_{n-1+\rho}^RV.$$  
It follows that  $$\split [\Cal I]_{n-1+\rho}+C_{n-1+\rho}&{}=[\operatorname{ann} \Phi]_{n-1+\rho}+C_{n-1+\rho},\quad  [\operatorname{ann} \Phi]_{n-1+\rho}  \cap  C_{n-1+\rho}=0,\\ \text{and}\quad [\Cal I]_{n-1+\rho}&{}\subseteq [\operatorname{ann} \Phi]_{n-1+\rho};\endsplit $$and therefore, 
 $[\Cal I]_{n-1+\rho}=[\operatorname{ann}\Phi]_{n-1+\rho}$. 

Finally, we consider $i>2n-2$. It is clear that $[\operatorname{ann} \Phi]_i=\operatorname {Sym}_i^RV$ and, if one makes the calculation analogous to (\tref{equiv'}), then it is not possible for     the bottom case to occur, so $[\Cal I]_i$ is also equal to $\operatorname {Sym}_i^RV$.

\medskip\flushpar(4.a) One consequence of (2) is that the $R$ sub-modules
$$[\operatorname{ann} \Phi]_{n+\rho }, \quad [\operatorname{ann} \Phi]_{n}\cdot \operatorname {Sym}_\rho ^RV\quad \text{and} \quad [\operatorname{ann} \Phi]_{n}\cdot (\operatorname {Sym}_1^RV)^\rho \tag\tnum{jabv}$$of $\operatorname {Sym}_{n+\rho }^RV$ are equal. The $R$-algebra homomorphism $\widehat{\Psi}$ carries this $R$-module to the ideal
  $\widehat{\Psi}([\operatorname{ann} \Phi]_{n+\rho })=\widehat{\Psi}([\operatorname{ann} \Phi]_{n})\cdot \widehat{\Psi}((\operatorname {Sym}_1^RV)^\rho )=IJ^\rho$ of $R$.

\medskip\flushpar(4.b) We will identify a set polynomials $\{p_\alpha\}$ in $\Bbb Z[\{X_1,\dots,X_d\},\{t_M\}]$, where 
$$\text{$M$ roam over the monomials of degree $2n-2$ in $d$ variables}.\tag\tnum{MON}$$
We have already picked a basis $x_1,\dots, x_d$ for $V$. If $M$ is a monomial from (\tref{MON}), then let $M|x$ represent the corresponding element of $\operatorname{Sym}_{2n-2}^RV$; we think of ``$M|x$'' as ``$M$ evaluated at $x_1,\dots,x_d$''. We will choose the $p_{\alpha}$ so that the set of $p_\alpha$, with $X_i$ evaluated at $\Psi(x_i)$ and $t_M$ evaluated at $\Phi(M|x)$,
 generates $J^\rho I$. 

 Our proof of (2) shows that $[\operatorname{ann}\Phi]_n$ is generated by the elements of (\tref{b'}) and (\tref{c'}). Thus, (\tref{jabv}) shows that  $J^\rho I$ is generated by $\widehat{\Psi}$ applied to the elements of
$(\operatorname {Sym}_1^RV)^\rho $ times the elements of 
 (\tref{b'}) and (\tref{c'}). 

Recall from Remark \tref{R1} that the matrix of the map $\pmb p_{n-1}^\Phi$ of (\tref{pn-1}) is $T_{\Phi}=[\Phi(m_im_j)]$, where $m_1,\dots,m_N$ is a list of the monomials of degree $n-1$ in $x_1,\dots,x_d$. Each entry in the matrix $T_\Phi$ is $\Phi(M|x)$ for some $M$ from (\tref{MON}). The matrix $T_{\Phi}$ is invertible over $R$ by the hypothesis that $\pmb p_{n-1}^\Phi$ is an isomorphism; thus, $\det T_{\Phi}$ is a unit in $R$ and the matrix for $\sigma_{n-1}$ is $T_{\Phi}^{-1}=\frac 1{\det T_{\Phi}}\operatorname{Adj}(T_{\Phi})$, where $\operatorname{Adj}(T_{\Phi})$ is the classical adjoint of $T_{\Phi}$. The classical adjoint of is built in a polynomial manner from the entires of $T_{\Phi}$;  the classical adjoint of  
$T_{\Phi}$ is the matrix for the map $(\det T_{\Phi})\sigma_{n-1}$; and the map $(\det T_{\Phi})\sigma_{n-1}$ differs from the map $\sigma_{n-1}$ by a unit of $R$ -- this unit does not affect the ideal generated by the image. Thus,
the ideal $J^\rho I$ is generated by  the elements 
$$\widehat{\Psi}\left((x_1,\dots,x_d)^\rho \cases \phantom{-}x_i(\det T_{\Phi})\sigma_{n-1}({x_1^*}^{(a_1)}\cdots{x_i^*}^{(a_i+1)}\cdots {x_j^*}^{(a_j)}\cdots {x_n^*}^{(a_n)})\\
-x_j(\det T_{\Phi})\sigma_{n-1}({x_1^*}^{(a_1)}\cdots{x_i^*}^{(a_i)}\cdots {x_j^*}^{(a_j+1)}\cdots {x_n^*}^{(a_n)})\endcases\right), $$   
  for   $i\neq j$ and   $\sum_\ell a_\ell=n-2$,  together with the elements
$$\widehat{\Psi}\left((x_1,\dots,x_d)^\rho x_i (\det T_{\Phi})\sigma_{n-1}({x_1^*}^{(a_1)}\cdots {x_n^*}^{(a_n)})\right),$$   with $a_i=0$ and  $\tsize{\sum_\ell a_\ell=n-1}$. Each of these elements of $R$ is  built in a polynomial manner from $\{\Psi(x_i)\}\cup \{\Phi(M|x)\}$.   
\qed \enddemo

\SectionNumber=\TMT\tNumber=1
\heading Section \number\SectionNumber. \quad The main theorem. \endheading

\definition{Data \tnum{a-tmt}} Consider the data $(R,V,n,\Psi,\Phi)$, where $R$ is a commutative Noetherian ring, 
 $V$ is a non-zero free $R$-module of finite rank, $n$ is a positive integer, 
$\Psi\: V\to R $ is an $R$-module homomorphism, and $\Phi$ is an element of $D_{2n-2}^R(V^*)$. Let $d$ be the rank of $V$,
  $\widehat{\Psi}\:\operatorname {Sym}_\bullet^R(V)\to R$ be the $R$-algebra homomorphism induced by $\Psi$, and $I$ and $J$ be the ideals $I=\widehat{\Psi}(\operatorname{ann} \Phi)$ and $J=\widehat{\Psi}(V)$ of $R$.
For each integer $i$ with $0\le i\le 2n-2$, let $\pmb p_{i}^{\Phi}\: \operatorname {Sym}^R_{i}V\to D_{2n-2-i}^R(V^*)$ be the $R$-module homomorphism defined by $\pmb p_{i}^{\Phi}(v_i)=v_i(\Phi)$, for $v_i\in \operatorname {Sym}^R_{i}V$. 
\enddefinition

In this section we produce a family of complexes $\Bbb G(R,V,n,\Psi,\Phi;r)$, one for each integer $r$ with $0\le r\le 2n-2$. The main result is Theorem \tref{EK-K'}, where we prove that if the grade of $J$ is at least $d$ and $\pmb p_{n-1}^\Phi$ is an isomorphism, then each complex $\Bbb G(R,V,n,\Psi,\Phi;r)$ is a resolution and the complexes $\Bbb G(R,V,n,\Psi,\Phi;n-1)$ and  $\Bbb G(R,V,n,\Psi,\Phi;n)$ both resolve $R/I$. 

Corollaries \tref{EK-K-g} and \tref{EK-K-2} are applications of Theorem \tref{EK-K'}. Corollary \tref{EK-K-g} treats the generic resolution $\widetilde{\Bbb G}(r)$ and Corollary \tref{EK-K-2} specializes the generic resolution in order to accomplish Project \tref{EK1}. Section \number\exmpls\ and the beginning of Example \tref{X} give examples of $\widetilde{\Bbb G}$ .

Recall the complexes $\Bbb L$ and $\Bbb K$ from Definition \tref{long}.

\proclaim{Observation \tnum{w-d}} Adopt Data {\rm\tref{a-tmt}} and let $r$ be an integer with $1\le r\le 2n-2$. Then the $R$-module homomorphism 
$\pmb p_{r}^{\Phi}\:\operatorname {Sym}^R_{r}(V)\to D_{2n-2-r}^R(V^*)$ induces a map of complexes $$E(R,V,n,\Psi,\Phi;r):\Bbb L(\Psi,r)\to \Bbb K(\Psi,2n-1-r)[-1]$$ as described:
$$ \eightpoint \matrix   & &0&\to& L_{d-1,r}&\to &\cdots &\to&L_{1,r}&\to&L_{0,r}&\to&R\\  & &\downarrow & &@V 1\otimes \pmb p_r^\Phi VV& & @V 1\otimes \pmb p_r^\Phi VV@V 1\otimes \pmb p_r^\Phi VV&\\ 0&\to &{\tsize\bigwedge}^d_RV&\to&K_{d-1,2n-2-r}&\to&\cdots &\to&K_{1,2n-2-r}&\to&K_{0,2n-2-r},
\endmatrix$$where $L_{a,b}$ and $K_{a,b}$ represent $L_{a,b}^R(V)$ and $K_{a,b}^R(V)$, respectively.\endproclaim

\demo{Proof} The proof is straightforward and short. \qed \enddemo

\remark{\bf Remark \tnum{bih}} If the ring $R$ is bi-graded and the $R$-homomorphisms $$\Psi\:R(-1,0)^d\to R\quad\text{and}\quad \Phi\:R(0,-1)^{\binom{2n+d-3}{d-1}}\to R$$ are bi-homogeneous, then the double complex $E(R,V,n,\Psi,\Phi;r)$ of Observation \tref{w-d} is bi-homogeneous. The entries of the matrix $L_{0,r}\to R$ have degree $(r,0)$; the entries of the matrices $L_{\ell,r}\to L_{\ell,r-1}$ and $K_{\ell,2n-2-r}\to K_{\ell-1,2n-2-r}$ have degree $(1,0)$; the entries of the matrices for $1\otimes \pmb p_r^\Phi$ have degree $(0,1)$; and the entries of the matrix $\bigwedge^d_RV\to K_{d-1,2n-2-r}$ have degree $(2n-1-r,0)$. The ranks of $L_{a,b}$ and $K_{a',b'}$ are given in (\tref{rank}) and (\tref{p207}). (See, also, (\tref{rank'}).) The double complex $E(R,V,n,\Psi,\Phi;r)$  is bi-homogeneous   with $$\matrix\format\l&\c& \l\\ L_{a,r}&=&R(-r-a,0)^{\binom{d+r-1}{a+r}\binom{a+r-1}a}\\\vspace{5pt}
K_{a,2n-2-r}&=&R(-r-a,-1)^{\binom{d+2n-2-r}{a}\binom{2n+d-3-r-a}{d-1-a}}\\\vspace{5pt}
{\tsize\bigwedge}^d_RV&=&R(-2n-d+2,-1)^1.\endmatrix   $$   \endremark

\definition{Definition \tnum{BBBG}} Adopt Data {\rm\tref{a-tmt}} and let $r$ be an integer with ${1\le r\le 2n-2}$.
Let  $\Bbb G(R,V,n,\Psi,\Phi;r)$   be the total complex of the double complex $E(R,V,n,\Psi,\Phi;r)$, from   Observation {\rm \tref{w-d}}, with $R\oplus K^R_{0,2n-2-r}$    in position $0$. \enddefinition

\proclaim{Observation \tnum{dual}} 
 Adopt Data {\rm\tref{a-tmt}} and let $r$ be an integer with $1\le r\le 2n-2$ and $\Bbb G(R,V,n,\Psi,\Phi;r)$ be the complex of Definition {\rm \tref{BBBG}}.
Then the complexes $$\Bbb G(R,V,n,\Psi,\Phi;r)\quad\text{and}\quad  [\operatorname{Hom}_R(\Bbb G(R,V,n,\Psi,\Phi;2n-1-r),R)\otimes_R{\tsize \bigwedge}^dV][-d]$$ of
Definition {\rm \tref{BBBG}} 
are isomorphic, where ``$[-d]$'' describes a shift in homological degree.
\endproclaim
\demo{Proof}Apply assertion (1) of Theorem \tref{BE}.
\qed \enddemo

\proclaim{Observation \tnum{resl}} Adopt Data {\rm\tref{a-tmt}} and let $r$ be an integer with $1\le r\le 2n-2$.
If $J$ has grade at least $d$, then the  complex $\Bbb G(R,V,n,\Psi,\Phi;r)$  of Definition {\rm \tref{BBBG}} is a resolution.
\endproclaim

\demo{Proof} Apply assertion (2) of  Theorem \tref{BE} to see that the complexes 
$$\Bbb L(\Psi,r)\quad  \text{and}\quad \Bbb K(\Psi, 2n-1-r)$$ 
have homology concentrated in position zero, then use the long exact sequence of homology associated to a mapping cone to conclude that $\Bbb G(R,V,n,\Psi,\Phi;r)$ also has homology concentrated in position zero.
\qed \enddemo

\proclaim{Theorem \tnum{EK-K'}} Adopt Data {\rm\tref{a-tmt}}. 
Assume that  $J$ is a proper ideal with grade at least $d$ and   $$\pmb p_{n-1}^{\Phi}\:\operatorname {Sym}^R_{n-1}(V)\to D_{n-1}^R(V^*)$$  is an isomorphism. For each index $r$, with $1\le r\le 2n-2$, let $\Bbb G(r)$ represent the complex $\Bbb G(R,V,n,\Psi,\Phi;r)$, of Definition {\rm \tref{BBBG}}.
The following statements hold.
\roster
\item Each of the ideals $J^\rho I$, with $0\le \rho\le n-2$, is a perfect, grade $d$, ideal of $R$.

\item  The ideal $I$ is a perfect, grade $d$, Gorenstein  ideal of $R$.

\item  The complex $\Bbb G(n+\rho)$ is a resolution of $R/J^\rho I$ by free $R$-modules, for $0\le \rho\le n-2$.

\item  The complexes $\Bbb G(n)$ and $\Bbb G(n-1)$ both resolve $R/I=\operatorname{Ext}^d_R(R/I,R)$.

\item  The complex $\Bbb G(n-1-\rho)$ is a resolution of $\operatorname{Ext}^d_R(R/J^\rho I,R)$ by free $R$-modules, for $0\le \rho\le n-2$.
\endroster
\endproclaim
\remark{Remark} We are primarily interested in the resolution $\Bbb G(n+\rho)$ of $R/J^\rho I$ when $\rho=0$. The argument which produces this resolution of $R/J^\rho I$ for $\rho=0$ also resolves $R/J^\rho I$ for $0\le \rho\le n-2$.
 The ideal $J^\rho I$ is a ``truncation'' of $I$. 
\endremark

\demo{Proof of Theorem {\rm\tref{EK-K'}}} First notice that the ideals $I$ and $J$ of $R$ have the same radical. Recall that $I=\widehat{\Psi}(\operatorname{ann} \Phi)$ and $J=\widehat{\Psi}(V)$. It is clear that, as ideals of $\operatorname {Sym}_{\bullet}^R(V)$,
$$(\operatorname{ann}\Phi)\subseteq  (\operatorname {Sym}_1^RV)\quad \text{and}\quad (\operatorname {Sym}_{2n-1}^RV)\subseteq (\operatorname{ann} \Phi).$$
It follows that
$$J^{2n-1}=(\widehat{\Psi}(V))^{2n-1}=\widehat{\Psi}(\operatorname {Sym}_{2n-1}^R(V)) \subseteq \widehat{\Psi}(\operatorname{ann}(\Phi)) \subseteq \widehat{\Psi}(V) =J;$$and therefore, $J^{2n-1}\subseteq I\subseteq J$. In particular, $$\text{the ideals $J$ and $J^\rho I$ have the same grade for all non-negative integers $\rho$.}\tag\tnum{grd}$$ The hypothesis that $d\le \operatorname{grade} J$ also guarantees that $d\le \operatorname{grade} J^\rho I$, for $0\le \rho$.

We saw in Observation \tref{resl} that $\Bbb G(r)$ is a resolution for every integer $r$ between  $1$ and $2n-2$.  The heart of the argument is the proof that
$$\operatorname{H}_0(\Bbb G(n-1+\lambda ))=\cases 
R/I&\text{if $\lambda =0$}\\
R/J^{\lambda -1}I &\text{if $1\le \lambda \le n-1$}.
\endcases\tag \tnum{hrt}$$ We assume  (\tref{hrt}) for the time-being and complete the proof. Assertion  (3) and the part of  assertion  (4) that $\Bbb G(n)$ and $\Bbb G(n-1)$ both resolve $R/I$ follow immediately from (\tref{hrt}). Fix an index $\rho$, with $0\le \rho\le n-2$. 
We have exhibited that $\Bbb G(n+\rho)$ is a free resolution of $R/J^\rho I$ of length $d$; thus, $$d\le \operatorname{grade} J^\rho I\le \text{the projective dimension of the $R$-module  $R/J^\rho I$}\le d,$$  and  $J^\rho I$ is a perfect, grade $d$, ideal of $R$. Assertion (1) is established and $$\operatorname{Hom}_R(\Bbb G(n+\rho),R)\to \operatorname{Ext}_R^d(R/J^\rho I,R)\to 0$$ is a resolution. Observation \tref{dual} yields assertion (5). 
We know from (4) that $\Bbb G(n-1)$ resolves $R/I$ and from (5) that $\Bbb G(n-1)$ resolves $\operatorname{Ext}^d_R(R/I,R)$. It follows that $R/I=\operatorname{Ext}^d_R(R/I,R)$, and the proof of (2) and (4) are complete.

We now prove (\tref{hrt}).    Write $\Bbb G$ for $\Bbb G(n-1+\lambda )$. We   compute $\operatorname{H}_0(\Bbb G)$. (The same ideas appear again in the proof of Lemma \tref{q-i}.)
The homology $\operatorname{H}_0(\Bbb G)$ is defined to be the cokernel of   
$$\matrix L^R_{0,n-1+\lambda }(V)\\\oplus\\ K^R_{1,n-1-\lambda }(V)\endmatrix @>\bmatrix \widehat{\Psi}&0\\\pmb p_{n-1+\lambda }^{\Phi}&\Psi\otimes 1\endbmatrix >> \matrix R\\\oplus \\ K^R_{0,n-1-\lambda }(V).\endmatrix $$ 
Recall that $L^R_{0,n-1+\lambda }(V)=\operatorname {Sym}_{n-1+\lambda }^RV$ and $K^R_{0,n-1-\lambda }=D_{n-1-\lambda }^R(V^*)$. Recall, also, the splitting map
$$\sigma_{n-1-\lambda }\: D_{n-1-\lambda }^R(V^*)\to \operatorname {Sym}_{n-1+\lambda }^RV$$ 
of Remark \tref{20.21'}.2.
We decompose $L^R_{0,n-1+\lambda }(V)=\operatorname {Sym}^R_{n-1+\lambda }V$ as $$\operatorname{ker} \pmb p_{n-1+\lambda }^{\Phi}\oplus \operatorname{im} {\sigma_{n-1-\lambda }}$$ and we employ the isomorphism
$$\matrix \operatorname{ker} \pmb p_{n-1+\lambda }^{\Phi}\\ \oplus \\K^R_{0,n-1-\lambda }(V)\endmatrix @>\bmatrix 1&0\\0&\sigma_{n-1-\lambda }\endbmatrix>\simeq> \left.\matrix \operatorname{ker} \pmb p_{n-1+\lambda }^{\Phi}\\ \oplus \\ \operatorname{im} {\sigma_{n-1-\lambda }}\endmatrix\right\} =L^R_{0,n-1+\lambda }(V)$$ to see that 
$\operatorname{H}_0(\Bbb G)$ is   the cokernel of 
$$\matrix  \operatorname{ker} \pmb p_{n-1+\lambda }^{\Phi}\\ \oplus \\K^R_{0,n-1-\lambda }(V)\\\oplus\\ K^R_{1,n-1-\lambda }(V)\endmatrix @>\bmatrix \widehat{\Psi}&\widehat{\Psi}\circ \sigma_{n-1-\lambda }&0\\0&1&\Psi\otimes 1\endbmatrix >> \matrix R\\\oplus \\ K^R_{0,n-1-\lambda }(V).\endmatrix  $$  Row and column operations   yield that $\operatorname{H}_0(\Bbb G)$ is the cokernel of 
$$\bmatrix \widehat{\Psi}&0&-\widehat{\Psi}\circ \sigma_{n-1-\lambda }\circ(\Psi\otimes 1)\\0&1&0\endbmatrix.\tag\tnum{map}$$
It is not difficult to see that the diagram
$$\CD K_{1,n-1-\lambda }^R(V) @>\Psi>> D_{n-1-\lambda }^R(V^*)@>\sigma_{n-1-\lambda } >>\operatorname {Sym}_{n-1+\lambda }^RV\\
@V 1\otimes \sigma_{n-1-\lambda } VV @. @V \widehat{\Psi}VV \\
V\otimes _R\operatorname {Sym}_{n-1+\lambda }^RV@>\text{mult}>>\operatorname {Sym}_{n+\lambda }^RV@>\widehat{\Psi}>> R\endCD$$ commutes. 
The clockwise composition from $K_{1,n-1-\lambda }^R(V)$ to $R$ is the map in the upper right-hand corner of (\tref{map}), up to sign.
The counter-clockwise composition  from $K_{1,n-1-\lambda }^R(V)$ to $\operatorname {Sym}_{n+\lambda }^RV$ is called $\alpha_{n-1-\lambda }$ (see Remark \tref{20.21'}.3) and the image of $\alpha_{n-1-\lambda }$ is called $A[n+\lambda ]$. 
It follows that $$\operatorname{H}_0(\Bbb G)=\frac{R}{\widehat{\Psi} 
(\vphantom{E^E}\ker (\pmb p_{n-1+\lambda }^{\Phi}) 
)+ \widehat{\Psi}\left( A[n+\lambda ]\right)}.$$
Use Observation \tref{gen'} to see that 
$$\widehat{\Psi} 
(\vphantom{E^E}\ker (\pmb p_{n-1+\lambda }^{\Phi}))=\widehat{\Psi} 
([\operatorname{ann} \Phi]_{n-1+\lambda })=\cases 0&\text{if $\lambda =0$},\\ J^{\lambda -1}I&\text{if $1\le \lambda \le n-1$},\endcases$$  
$$\matrix \format\l&\quad\l\\
\widehat{\Psi}(A[n+\lambda ])=\widehat{\Psi}(\operatorname{ann}\Phi)=I&\text{if $\lambda =0$, and}\\\vspace{5pt}
 \widehat{\Psi}(A[n+\lambda ]) \subseteq \widehat{\Psi}([\operatorname{ann} \Phi]_{n+\lambda })=J^\lambda I\subseteq J^{\lambda -1}I&\text{if $1\le \lambda \le n-1$}.
\endmatrix$$ 
Thus, (\tref{hrt}) holds 
and the proof is complete.\qed \enddemo

The complexes of Definition \tref{BBBG} can be built generically. (When we consider generic data, there would be no loss of generality if we had written    ``$\Bbb Z$'' every where we have written ``$R_0$''  because, the generic objects constructed over $\Bbb Z$ would become generic objects constructed over $R_0$ once we applied the base change $R_0\otimes_{\Bbb Z}\underline{\phantom{X}}$.)    
\definition{Definition \tnum{19.xms}} Consider    $(R_0,U,n)$, where $R_0$ is a commutative Noetherian ring, 
$U$ is a free $R_0$-module of rank $d$,  and  $n$ is a  positive integer. We define the {\it generic ring} $\widetilde{R}=\widetilde{R}(R_0,U,n)$, the {\it generic ideals} $\widetilde{I}=\widetilde{I}(R_0,U,n)$ and $\widetilde{J}=\widetilde{J}(R_0,U,n)$ of $\widetilde{R}$, and, for each integer $r$ with $1\le r\le 2n-2$, the {\it generic complex} $\widetilde{\Bbb G}(R_0,U,n;r)$. 
 View $U\oplus \operatorname {Sym}_{2n-2}^{R_0}{U}$ as a bi-graded free $R_0$-module, where the elements of $U$ have degree $(1,0)$ and the elements of $\operatorname {Sym}_{2n-2}^{R_0}U$ have degree $(0,1)$. 
Let
 $\widetilde{R}$ be the bi-graded $R_0$-algebra $\operatorname {Sym}_{\bullet}^{R_0}(U\oplus \operatorname {Sym}_{2n-2}^{R_0}U)$. Define the $R_0$-module inclusion homomorphisms $i_1\:U\to  \widetilde{R}$ and $i_2\:\operatorname {Sym}^{R_0}_{2n-2}U\to \widetilde{R}$ by
$$U\to U\oplus 0=[\widetilde{R}]_{(1,0)}\subseteq  \widetilde{R}\quad\text{and}\quad 
\operatorname {Sym}^{R_0}_{2n-2}U\to 0\oplus \operatorname {Sym}^{R_0}_{2n-2}U =[\widetilde{R}]_{(0,1)}\subseteq \widetilde{R},$$ respectively. 
Let $V$ be the  free  $\widetilde{R}$-module $ V= \widetilde{R}\otimes_{R_0} U$ and
define $\Psi\:   V\to   \widetilde{R}$ and $ \Phi\:\operatorname {Sym}_{2n-2}^{\widetilde{R}}V\to  \widetilde{R}$ to be the $ \widetilde{R}$-module homomorphisms

$$\matrix \format \r&\c& \l&\c&\c&\c&\c&\c\\  V&{}={}& \widetilde{R}\otimes_{R_0} U&{}@> 1\otimes i_1 >>{}& \widetilde{R}\otimes _{R_0}  \widetilde{R}&{}@> \text{multiplication}>>{}& \widetilde{R}&\quad \text{and}\\
\operatorname {Sym}_{2n-2}^{\widetilde{R}}V&{}={}&  \widetilde{R}\otimes_{R_0} \operatorname {Sym}_{2n-2}^{R_0} U&{}@> 1\otimes i_2 >>{}&\widetilde{R}\otimes _{R_0} \widetilde{R}&{}@>\text{multiplication}>>{}&\widetilde{R},\endmatrix\tag\tnum{Phi}
$$ respectively. Let 
 $\widehat{\Psi}\:\operatorname {Sym}_\bullet^{\widetilde{R}}(V)\to \widetilde{R}$ be the $\widetilde{R}$-algebra homomorphism induced by $\Psi$,
$\widetilde{I}$ and $\widetilde{J}$ be the ideals $\widetilde{I}=\widehat{\Psi}(\operatorname{ann} \Phi)$ and $\widetilde{J}=\widehat{\Psi}(V)$ of $\widetilde{R}$.
For each integer $r$ with $1\le r\le 2n-2$, define $\widetilde{\Bbb G}(R_0,U,n;r)$ to be the complex $\Bbb G(\widetilde{R},V,n,\Psi,\Phi;r)$ of Definition \tref{BBBG}. \enddefinition

\remark{\bf Remark \tnum{19.3}} Continue with the notation and hypotheses of Definition \tref{19.xms}. Suppose that $U$ has basis $x_1,\dots,x_d$. In this case, we may view $\widetilde{R}$ as the bi-graded polynomial ring $\widetilde{R}=R_0[x_1,\dots,x_d, \{t_{M}\}]$, 
where $$\text{$M$ roams over all monomials of degree $2n-2$ in $\{x_1,\dots,x_d\}$.}\tag\tnum{abc'}$$
The variables $x_i$ have bi-degree $(1,0)$ and each variable $t_{M}$ has bi-degree $(0,1)$. The map $ \Phi\:\operatorname {Sym}_{2n-2}^{\widetilde{R}}V\to  \widetilde{R}$ from (\tref{Phi}) is the same as the element 
$$\sum_{\text{(\tref{abc'})}} t_{x_1^{a_1}\cdots x_d^{a_d}}\otimes {x_1^*}^{(a_1)}\cdots {x_d^*}^{(a_d)}\in \widetilde{R}\otimes_{R_0}D_{2n-2}^{R_0}U^*=D_{2n-2}^{\widetilde{R}}(V^*),\tag\tnum{appr}$$
where $^*$ in $U^*$ means $R_0$-dual and $^*$ in $V^*$ means $\widetilde{R}$-dual.
\endremark

Our first application of Theorem \tref{EK-K'} is to the generic situation.
\proclaim{Corollary \tnum{EK-K-g}} Consider the data   $(R_0,U,n)$, where $R_0$ is a commutative Noetherian ring, 
$U$ is a free $R_0$-module of positive rank $d$,    $n$ is a  positive integers. Let $\widetilde{R}=\widetilde{R}(R_0,U,n)$, $\widetilde{I}=\widetilde{I}(R_0,U,n)$, $\widetilde{J}=\widetilde{J}(R_0,U,n)$, and for each integer $r$, with $1\le r\le 2n-2$, $\widetilde{\Bbb G} (r)=\widetilde{\Bbb G}(R_0,U,n;r)$, be the generic ring, ideals, and complexes of Definition {\rm\tref{19.xms}}, and let $\pmb \delta$ be $\det T_{\Phi}$ for the $\widetilde{R}$-module $\Phi$ of {\rm(\tref{Phi})}.
The following statements hold.
\roster
\item Each of the ideals $\widetilde{J}^\rho \widetilde{I}\widetilde{R}_{\pmb \delta}$, with $0\le \rho\le n-2$, is a perfect, grade $d$, ideal of $\widetilde{R}_{\pmb \delta}$.

\item  The ideal $\widetilde{I}\widetilde{R}_{\pmb \delta}$ is a perfect, grade $d$, Gorenstein  ideal of $\widetilde{R}_{\pmb \delta}$.

\item  The complex $\widetilde{\Bbb G}(n+\rho)_{\pmb \delta}$ is a resolution of $\widetilde{R}_{\pmb \delta}/\widetilde{J}^\rho \widetilde{I}\widetilde{R}_{\pmb \delta}$ by free $\widetilde{R}_{\pmb \delta}$-modules, for $0\le \rho\le n-2$.

\item  The complexes $\widetilde{\Bbb G}(n)_{\pmb \delta}$ and $\widetilde{\Bbb G}(n-1)_{\pmb \delta}$ both are resolutions of  $$\widetilde{R}_{\pmb \delta}/\widetilde{I}\widetilde{R}_{\pmb \delta}=\operatorname{Ext}^d_{\widetilde{R}_{\pmb \delta}}(\widetilde{R}_{\pmb \delta}/\widetilde{I}\widetilde{R}_{\pmb \delta},\widetilde{R}_{\pmb \delta})$$by free $\widetilde{R}_{\pmb \delta}$-modules.

\item The complex $\widetilde{\Bbb G}(n-1-\rho)_{\pmb \delta}$ is a resolution of $\operatorname{Ext}^d(\widetilde{R}_{\pmb \delta}/\widetilde{J}^\rho \widetilde{I}\widetilde{R}_{\pmb \delta},\widetilde{R}_{\pmb \delta})$ by free $\widetilde{R}_{\pmb \delta}$-modules, for   $0\le \rho\le n-2$.
\endroster
\endproclaim

\remark{\bf Remark \tnum{knew}}The sense in which ``$\pmb \delta=\det T_{\Phi}$'' is coordinate-free is explained in Remark \tref{R1}.3. In the language of     Remark {\rm\tref{19.3}}, if $m_1,\dots,m_N$ is a list of the monomials of degree $n-1$ in $x_1,\dots,x_d$, then $\pmb \delta$ is equal to $\det(t_{m_im_j})$. Each $t_{m_im_j}$ is an indeterminate of the polynomial ring $\widetilde{R}$ from (\tref{abc'}). \endremark

\demo{Proof of Corollary {\rm\tref{EK-K-g}}} We apply Theorem \tref{EK-K'}. In the language of Remark \tref{19.3}, the ring $\widetilde{R}$ is the polynomial ring $R_0[x_1,\dots,x_d,\{t_M\}]$ and the ideal $\widetilde{J}$ is generated by the variables $x_1,\dots,x_d$;  thus, $\widetilde{J}$ is a proper ideal of grade $d$. Furthermore, the map $$\pmb p_{n-1}^{\Phi}\: \operatorname{Sym}_{n-1}^{\widetilde {R}_{\pmb \delta}}V_{\pmb \delta}\to  D_{n-1}^{\widetilde {R}_{\pmb \delta}}(V_{\pmb \delta}^*)$$
is invertible because its determinant is the unit $\pmb \delta$; see Definition \tref{T} and Remark \tref{R1}.1. The hypotheses of Theorem \tref{EK-K'} are satisfied. The conclusions follow. 
\qed \enddemo

Our second application of Theorem \tref{EK-K'} is to the study of the ideals in the set  $\Bbb I_n(R_0,U)$ of Definition \tref{inr}.   Fix the data $(R_0,U,n)$, where $R_0$ is a commutative Noetherian ring, $U$ is a free $R_0$-module of positive rank $d$, and $n$ is a positive integer. Let $\widetilde{R}$, $\widetilde{I}$, and  $\widetilde{\Bbb G}(n)_{\pmb \delta}$,   be the generic ring, ideal, and resolution created for the data $(R_0,U,n)$ as described in Definition \tref{19.xms} and Corollary \tref{EK-K-g}. Let $P$ be the polynomial ring $P=\operatorname{Sym}_\bullet^{R_0}U$ and $I=\operatorname{ann}\phi$ be an ideal in $P$ from the set  $\Bbb I_n(R_0,U)$; so, in particular, $\phi\in D_{2n-2}^{R_0}(U^*)$ and $\det T_{\phi}$ (as described in Definition \tref{T}) is a unit of $R_0$. We prove in Corollary \tref{EK-K-2} that $\phi$ naturally induces a $P$-algebra homomorphism $\widehat{\phi}\:\widetilde{R}\to P$   so that 
\roster 
\item $ \widehat{\phi}(\widetilde{I})=I$, and

\vskip5pt 
\item $P\otimes_{\widetilde{R}}\widetilde{\Bbb G}(n)_{\pmb \delta}$ is a resolution of $P/I$ by free $P$-modules.\endroster
We first give the coordinate-free, official, description of $\widehat{\phi}$. The map $\widehat{\phi}$ exists for any $\phi\in D_{2n-2}^{R_0}(U^*)$. The condition that $\det T_\phi$ is a unit in $R_0$ is not needed until we establish properties of $\widehat{\phi}$. The rings $P$ and $\widetilde{R}$ are defined to be 
$$P=\operatorname{Sym}_\bullet^{R_0}U\quad\text{and}\quad \widetilde{R}=\operatorname{Sym}_\bullet^{R_0}(U\oplus \operatorname{Sym}_{2n-2}^{R_0}U); $$ thus, 
$\widetilde{R}$ is also equal to 
$$\widetilde{R}=\operatorname{Sym}_\bullet^{P}(P\otimes_{R_0}\operatorname{Sym}^{R_0}_{2n-2}U).$$
The element $\phi$ of $D_{2n-2}^{R_0}(U^*)$ {\bf IS} an $R_0$-module homomorphism $\operatorname{Sym}_{2n-2}^{R_0}U\to R_0$; therefore $\phi$ induces a $P$-module homomorphism 
$$1\otimes \phi: P\otimes_{R_0}\operatorname{Sym}_{2n-2}^{R_0}U\to P\otimes_{R_0}R_0=P;$$ and, according to the defining property of symmetric algebra, $1\otimes \phi$ induces a $P$-algebra homomorphism 
$$ \widehat{1\otimes\phi}\: \widetilde{R}=\operatorname{Sym}_\bullet^{P}(P\otimes_{R_0}\operatorname{Sym}^{R_0}_{2n-2}U)\to P,\tag\tnum{phihat}$$which we abbreviate as $\widehat{\phi}\: \widetilde{R}\to P$. 
In practice, $\widetilde{R}=P[\{t_M\}]$ for indeterminates $t_M$ as described in (\tref{abc'}) and $\widehat{\phi}$ is the $P$-algebra homomorphism $$\text{$\widehat{\phi}\:\widetilde{R}=P[\{t_M\}]\to P$, with $\widehat{\phi}(t_M)=\phi(M)\in R_0$.}\tag\tnum{???}$$
 Notice that the element $\phi$ of $D_{2n-2}^{R_0}(U^*)$ has the form
$$\phi=\sum_{\text{(\tref{abc'})}} \tau_{x_1^{a_1}\cdots x_d^{a_d}}  {x_1^*}^{(a_1)}\cdots {x_d^*}^{(a_d)}\in D_{2n-2}^{R_0}U^*, \tag\tnum{b-appr}$$ for some elements $\tau_{M}$ in $R_0$. Apply $\phi$ to the monomial $M$ from (\tref{abc'}) to see that $\tau_M=\phi(M)$. 
It follows   that if $\Phi\in \widetilde{R}\otimes_{R_0}D_{2n-2}^{R_0}(U^*)$ is the element of (\tref{appr}) and (\tref{Phi}), and $P$ is an $\widetilde{R}$-algebra by way of $\widehat{\phi}\:\widetilde{R}\to P$, then
$$P\otimes_{\widetilde{R}}\Phi \text{ is equal to }1\otimes \phi \in P\otimes_{R_0}D^{R_0}_{2n-2}(U^*).\tag\tnum{carries}$$

\proclaim{Corollary \tnum{EK-K-2}} Fix $(R_0,U,n)$, where $R_0$ is a commutative Noetherian ring, $U$ is a free $R_0$-module of positive rank $d$, and $n$ is a positive integer.  Let 
$\widetilde{R}=\widetilde{R}(R_0,U,n)$, $\widetilde{I}=\widetilde{I}(R_0,U,n)$, $\widetilde{J}=\widetilde{J}(R_0,U,n)$, and, for each integer $r$ with $1\le r\le 2n-2$,  $\widetilde{\Bbb G}(r)=\widetilde{\Bbb G}(R_0,U,n;r)$ be the generic ring, ideals, and complexes of Definition {\rm\tref{19.xms}}. Let $\pmb \delta$ be the element $\det T_{\Phi}$ in $\widetilde{R}$ as described in Definition {\rm\tref{19.xms}} and Remark {\rm\tref{knew}}, $P$ be the polynomial ring $P=\operatorname{Sym}_\bullet^{R_0}U$, and $J$ be the ideal $(\operatorname{Sym}^{R_0}_1(U))$ of $P$. Let $I=\operatorname{ann} \phi$ be an ideal of $P$ from the set $\Bbb I_{n}(R_0,U)$ of Definition {\rm \tref{inr}} with $\phi\in D_{2n-2}^{R_0}(U^*)$ and $\det T_{\phi}$ a unit of $R_0$.   View $P$ as an $\widetilde{R}$-algebra by way of the $P$-algebra homomorphism $\widehat{\phi}\:\widetilde{R}\to P$ of {\rm(\tref{phihat})} and {\rm(\tref{???})}. The following statements hold.\vphantom{\tnum{flat}}
\roster
\item The ideals $\widehat{\phi}(\widetilde{J}^\rho\widetilde{I})$ and $J^\rho I$ of $P$ are equal, for all non-negative integers $\rho$.

\item
The complex $P\otimes_{\widetilde{R}}\widetilde{\Bbb G}(n+\rho)_{\pmb \delta}$ is a resolution of $P/J^\rho I$ by free $P$-modules, for $0\le \rho\le n-2$.

\item 
The complex $P\otimes_{\widetilde{R}}\widetilde{\Bbb G}(n-1-\rho)_{\pmb \delta}$ is a resolution of $\operatorname{Ext}^d_P(P/J^\rho I,P)$ by free $P$-modules, for $0\le \rho\le n-2$.

\item  
The complexes $P\otimes_{\widetilde{R}}\widetilde{\Bbb G}(n)_{\pmb \delta}$ and $P\otimes_{\widetilde{R}}\widetilde{\Bbb G}(n-1)_{\pmb \delta}$ both resolve $$P/I=\operatorname{Ext}^d_P(P/I,P).$$

\item  Each of the ideals $J^\rho I$, with $0\le \rho\le n-2$, is a perfect, grade $d$, ideal of $P$.

\item  The ideal $I$ is a perfect, grade $d$, Gorenstein  ideal of $P$.

\item Each of the natural ring homomorphisms 
$$\tsize \operatorname{Sym}_\bullet^{R_0}(\operatorname{Sym}_{2n-2}^{R_0}U)_{\pmb \delta} \to \left(\frac{\widetilde{R}}{\vphantom{\widetilde{\widetilde{I}}}\widetilde{J}^\rho\widetilde{I}}\right)_{\pmb \delta},\tag\tref{flat}$$ for $0\le \rho\le n-2$ is flat. 
\endroster 
\endproclaim
 
\remark{\bf Remark \tnum{ff}} In the language of Remark \tref{19.3} the ring homomorphism (\tref{flat}) is
$$\tsize R_0[\{t_M\}]_{\pmb \delta} \to  \left(\frac{R_0[x_1,\dots,x_d,\{t_M\}]}{\vphantom{\widetilde{\widetilde{I}}}\widetilde{J}^\rho\widetilde{I}}\right)_{\pmb \delta}.\tag\tnum{ff-loc}$$In particular, when $\rho=0$ and $R_0=\pmb k $ is a field,
then  
 $$\tsize \pmb k[\{t_M\}]_{\pmb \delta} \to \left(\frac{\widetilde{R}}{\vphantom{\widetilde{\widetilde{I}}}\widetilde{I}}\right)_{\pmb \delta}\tag\tnum{bblop}$$ is a flat family of 
$\pmb k$-algebras parameterized by $\Bbb I_{n}^{[d]}(\pmb k)$ in the sense that every algebra $\pmb k[x_1,\dots,x_d]/I$, with $I\in \Bbb I_{n}^{[d]}(\pmb k)$, is a fiber of (\tref{bblop}).\endremark  

\demo{Proof of Corollary {\rm\tref{EK-K-2}}} The element $\Phi$ of $D_{2n-2}^{\widetilde {R}}(\widetilde {R}\otimes_{R_0}U)$ is defined in (\tref{Phi}). Recall from (\tref{carries}) that $P\otimes_{\widetilde {R}}\Phi=1\otimes \phi$ in $P\otimes_{\widetilde {R}}D_{2n-2}^{R_0}(U^*)$. One consequence is that $\widehat{\phi}$ carries the element $\pmb \delta$ of $\widetilde {R}$ to the unit $\det T_{\phi}$ of $R_0$; therefore, $\widehat{\phi}\:\widetilde {R}\to P$ automatically induces a well-defined $P$-algebra homomorphism  $\widehat{\phi}\:\widetilde {R}_{\pmb \delta}\to P$. The complexes $\widetilde{\Bbb G}(r)_{\pmb \delta}$ are resolutions by Corollary \tref{EK-K-g}. The complexes  $P\otimes_{\widetilde {R}}\widetilde{\Bbb G}(r)_{\pmb \delta}$ are resolutions by  Theorem \tref{EK-K'} because the two hypotheses of Theorem \tref{EK-K'} are satisfied. The first hypothesis concerns the grade of the ideal $\widehat{\phi}(\widetilde{J})$. We see that $\widehat{\phi}(\widetilde{J})$ is equal to the ideal $J$ of $P$ and that this ideal has grade at least $d$ (Indeed, when one uses the language of Remark \tref{19.3}, $J$ is the ideal $(x_1,\dots,x_d)$ in the polynomial ring $P=R_0[x_1,\dots,x_d]$.) 
The second hypothesis concerns the element $\det T_{P\otimes_{\widetilde{R}} \Phi}$ of $P$. We have already observed that 
$P\otimes_{\widetilde{R}} \Phi=1\otimes \phi$. It follows that $\det T_{P\otimes_{\widetilde{R}} \Phi}$ is equal to the unit $\det T_\phi$ of $P$. 

We use the language of Remark \tref{19.3} and (\tref{b-appr}) to prove (1). Let $\rho$ be a non-negative integer. We saw in Observation \tref{gen'}.4b that the ideal $\widetilde{J}^\rho\widetilde{I}$ of $\widetilde{R}$ is built in a polynomial manner from the data $\{x_i\}\cup\{t_M\}$ and that the ideal $ J^\rho I $ of $P$ is built using the same polynomials from the data $\{x_i\}\cup\{\tau_M\}$. The $P$-algebra homomorphism $\widehat{\phi}\:\widetilde {R}\to P$ carries $t_M$ to $\tau_M$. Thus, $\widehat{\phi}(\widetilde{J}^\rho\widetilde{I})=J^\rho I $ and (1) is established. 
Assertions (2) -- (6) now follow immediately from Theorem \tref{EK-K'}. 

Assertion (7) is essentially obvious. Fix $\rho$, with $0\le \rho\le n-2$. We use the language of Remark \tref{19.3}. To prove that (\tref{ff-loc}) is flat, it suffices to prove the result locally and therefore, according to the local criterion for flatness, (see, for example, \cite{\rref{M80},~Thm.~49} or \cite{\rref{E95},~Thm.~6.8}) it suffices to prove that
$$\tsize \operatorname{Tor}_1^{A}(\frac{A}{\frak p A},\frac{B}{\widetilde{J}^\rho\widetilde{I}B})=0,\tag\tnum{tor}$$
where $A= R_0[\{t_M\}]_{\frak p}$, $B=R_0[\{x_i\},\{t_M\}]_{\frak P}$, $\frak P$ is a prime ideal of $R_0[\{x_i\},\{t_M\}]$ which contains $\widetilde{J}^\rho\widetilde{I}$ and does not contain $\pmb \delta$, and $\frak p=R_0[\{t_M\}]\cap \frak P$. Apply Theorem \tref{EK-K'} to $$\tsize \frac{A}{\frak p A}\otimes_{\widetilde{R}_\frak P}\widetilde{\Bbb G}(n+\rho)_{\frak P},\tag\tnum{et}$$ which is the complex $\Bbb G(n+\rho)$ built over the field $\frac{A}{\frak p A}$ with Macaulay inverse system $\frac{A}{\frak p A}\otimes_{\widetilde{R}_\frak P}\Phi$. The Macaulay inverse system still induces an  isomorphism from $\operatorname{Sym}_{n-1}$ to  $D_{n-1}$ (because a ring homomorphism can not send a unit to zero), and the image of $\widetilde{J}$ in $B$ is still a proper ideal (recall, from the proof of Theorem \tref{EK-K'},  that $\widetilde{J}$ and $\widetilde{I}$ have the same radical) of grade at least $d$. Thus, (\tref{et}) is a resolution and (\tref{tor}) holds. \qed \enddemo

\SectionNumber=\exmpls\tNumber=1
\heading Section \number\SectionNumber. \quad Examples of the resolution $\widetilde{\Bbb G}(n)$.
\endheading 
 
In Theorem \tref{mon} we prove that the complex $\widetilde{\Bbb G}(r)$ is a monomial complex. To do this we introduce the standard bases for the Schur and Weyl modules that we call $L_{p,q}$ and $K_{p,q}$. In Example \tref{n=3}, we use these bases to exhibit the matrices of $\widetilde{\Bbb G}(r)$, when $d=n=r=3$. The matrices of $\widetilde{\Bbb G}(r)$, when $d=3$ and $n=r=2$ are given at the beginning of Example \tref{X}.

\definition{Data \tnum{data5}} Consider    $(R_0,U,n)$, where $R_0$ is a commutative Noetherian ring, 
$U$ is a free $R_0$-module of rank $d$,  and  $n$ is a  positive integer.   
Let  $\widetilde{R}=\widetilde{R}(R_0,U,n)$, $\widetilde{I}=\widetilde{I}(R_0,U,n)$, $\widetilde{J}=\widetilde{J}(R_0,U,n)$,  and $\widetilde{\Bbb G}(r)=\widetilde{\Bbb G}(R_0,U,n;r)$, for $1\le r\le 2n-2$,  be the generic ring, the generic ideals, and  the generic complexes of Definition {\rm\tref{19.xms}}. Write $L_{p,q}$ and $K_{p,q}$ for $L_{p,q}^{R_0}U$ and $K_{p,q}^{R_0}U$, respectively. Let $\Phi$ be the element of $\widetilde{R}\otimes_{R_0} D_{2n-2}^{R_0}(U^*)$ which is described in (\tref{Phi}) and (\tref{appr}), and let $\pmb \delta$ be the element $\det T_{\Phi}$ of $\widetilde{R}$ as described in Corollary \tref{EK-K-g} and Remark \tref{knew}.
\enddefinition

\proclaim{Theorem \tnum{mon}} Adopt Data {\rm\tref{data5}}. Then 
one can choose bases for the free modules of $\widetilde{\Bbb G}(r)$ so every entry of every matrix is a signed    monomial. \endproclaim

\demo{Proof} Let $x_1,\dots,x_d$ be a basis for $U$ and think of $\widetilde{R}$ as the polynomial ring $R_0[x_1,\dots,x_d,\{t_M\}]$ as described in Remark \tref{19.3}. 
We see from Definition \tref{BBBG} that $\widetilde{\Bbb G}(r)$ is the mapping cone of\vphantom{\tnum{d}}
$$ \smallmatrix
{\tsize{(\tref{d})}}                                          & &\widetilde{R}\otimes_{R_0}L_{d-1,r}  &\to&\dots&\to&\widetilde{R}\otimes_{R_0}L_{p,r}      &\to&\dots&\to&\widetilde{R}\otimes_{R_0}L_{0,r} & \to&\widetilde{R}\\\vspace{7pt}
                               &   &\downarrow&&     &   &\downarrow                                   &   &     &   &
\downarrow                               &
\\\vspace{7pt}
\widetilde{R}\otimes_{R_0}\bigwedge^d_{R_0}U&\to&\widetilde{R}\otimes_{R_0}K_{d-1,2n-2-r} &\to&\dots&\to&\widetilde{R}\otimes_{R_0}K_{p,2n-2-r}&\to&\dots&\to&\widetilde{R}\otimes_{R_0}K_{0,2n-2-r}.\\
\endsmallmatrix$$
The top complex in the above diagram is a minimal resolution of $\widetilde{R}/(x_1,\dots,x_d)^r$ by free $\widetilde{R}$-modules and the bottom complex is the dual of a minimal resolution of $\widetilde{R}/(x_1,\dots,x_d)^{2n-1-r}$. Both of these complexes are naturally monomial complexes. We will use   well-understood bases for each of the modules; and therefore,  it will not be difficult to demonstrate that the matrices  for the horizontal maps have monomial entries. (We consider the horizontal maps at the end of the proof.) The interesting part of the proof involves the vertical maps.

Let $\omega=x_1\wedge \dots\wedge x_d$ be the basis for $\bigwedge^d_{R_0}U$. 
There are standard bases for $L_{p,q}$ and $K_{p,q}$. 
We will define these bases and make a few remarks. More details may be found in Remark \tref{schur}, \cite{\rref{W}}, \cite{\rref{BoBu}, Sect.~III.1},  and elsewhere. These bases are usually exhibited as tableau; our tableau are simply hooks, so we will simply record the information without distinguishing between the row and the column. The basis for $L_{p,q}$ is
$$\left\{\ell_{\pmb a;\pmb b}\left\vert \matrix \format\l\\ \pmb a\text{ is }a_1<\dots<a_{p+1},\\\pmb b\text{ is }b_1\le \dots\le b_{q-1},\text{ and}\\a_1\le b_1 \endmatrix \right.\right\}\tag\tnum{lpb}$$ and the basis for $K_{p,q}$ is
$$\left\{k_{\pmb a;\pmb b}\left\vert \matrix \format\l\\ \pmb a\text{ is }a_1<\dots<a_{d-p-1},\\\pmb b\text{ is }b_1\le \dots\le b_{q+1},\text{ and}\\b_1< a_1 \endmatrix \right.\right\},\tag\tnum{kpb}$$
where $$\ell_{\pmb a;\pmb b}=\kappa(x_{a_1}\wedge \dots\wedge x_{a_{p+1}}\otimes x_{b_1}\cdot \ldots\cdot x_{b_{q-1}})\in L_{p,q}\subseteq {\tsize\bigwedge}^p_{R_0}U\otimes_{R_0} \operatorname{Sym}^{R_0}_qU\text{ and}$$
$$\eightpoint k_{\pmb a;\pmb b}=\eta\left((x_{a_1}^*\wedge\dots\wedge x_{a_{d-p-1}}^*)(\omega)\otimes {x_1^*}^{(\beta_1)}\cdot \ldots\cdot  {x_d^*}^{(\beta_d)}\right)\in K_{p,q}\subseteq {\tsize\bigwedge}^p_{R_0}U\otimes_{R_0}D^{R_0}_{q}(U^*),$$for
$$\pmb b=(\underbrace{1,\ldots,1}_{\beta_1},\underbrace{2,\ldots,2}_{\beta_2},\cdots,\underbrace{d,\ldots,d}_{\beta_d}), {\tsize\text{ with $\sum\beta_i=q+1$}}.\tag\tnum{alt}$$ The maps $\kappa$ and $\eta$ are defined in (\tref{kappa}).
To show that (\tref{lpb}) is a basis for $L_{p,q}$ one can verify that (\tref{lpb}) contains $\operatorname{rank}L_{p,q}$ elements (see (\tref{rank})) and that (\tref{lpb}) spans $\kappa(\bigwedge^{p+1}_{R_0}U\otimes_{R_0} \operatorname{Sym}_{q-1}^{R_0}U)=L_{p,q}$. The assertion about spanning is obvious. Indeed, if $X=x_{a_1}\wedge\dots\wedge x_{a_{p+1}}\otimes x_{b_1}\cdots x_{b_{q-1}}\in \bigwedge^{p+1}_{R_0}U\otimes_{R_0} \operatorname{Sym}_{q-1}^{R_0}U$ with $a_1<\dots<a_{p+1}$, $b_1\le \dots\le b_{q-1}$, {\bf but} $b_1<a_1$, then \vphantom{\tnum{****'}}
$$\eightpoint \alignat1 0&{}=\kappa\kappa(x_{b_1}\wedge x_{a_1}\wedge\dots\wedge x_{a_{p+1}}\otimes x_{b_2}\cdots x_{b_{q-1}})\\&{}=\kappa(X)+\text{a linear combination of elements of (\tref{lpb}) with coefficients from $\{+1,-1\}$}.\tag{\tref{****'}}\endalignat$$
In a similar manner, one shows that (\tref{kpb}) is a basis for $K_{p,q}$ by showing that (\tref{kpb}) contains
$$\operatorname{rank} K_{p,q}(U)=\binom{d+q}p\binom{d+q-p-1}{q}\tag\tnum{rank'}$$elements (we used (\tref{rank}) and (\tref{p207}) to compute this number) and if $$X=(x_{a_1}^*\wedge\dots\wedge x_{a_{d-p-1}}^*)(\omega)\otimes {x_1^*}^{(\beta_1)}\cdot \ldots\cdot  {x_d^*}^{(\beta_d)}\in {\tsize\bigwedge}^{p+1}_{R_0}U\otimes_{R_0}D^{R_0}_{q+1}(U^*)$$
 with $a_1<\dots<a_{d-p-1}$, $b_1\le \dots\le b_{q+1}$ for  $(b_1,\dots,b_{q+1})$ given in (\tref{alt}), {\bf but} $a_1\le b_1$, then \vphantom{\tnum{****}}
$$ \eightpoint \alignat1 0&{}={\cases\eta\eta\left((x_{a_2}^*\wedge\dots\wedge x_{a_{d-p-1}}^*)(\omega)\otimes {x_{a_1}^*}{x_{b_1}^*}^{(\beta_{b_1})}\cdot \ldots\cdot  {x_d^*}^{(\beta_d)}\right)&\text{if $a_1<b_1$}\\
\eta\eta\left((x_{a_2}^*\wedge\dots\wedge x_{a_{d-p-1}}^*)(\omega)\otimes {x_{b_1}^*}^{(\beta_{b_1}+1)}\cdot \ldots\cdot  {x_d^*}^{(\beta_d)}\right)&\text{if $a_1=b_1$}\\
\endcases}
\\&{}=\eta(X)+\text{a linear combination of elements of (\tref{kpb}) with coefficients from $\{+1,-1\}$}.\tag{\tref{****}}\endalignat$$
The critical calculation involves a careful analysis of (\tref{****}).

\medskip \flushpar{\bf Claim \tnum{CC}.} Fix positive integers $P$,  $Q$, and $a_1<\dots<a_P$, with $P\le d-1$, $Q\le 2n-2$,   and $a_P\le d$. Let $\beta_1,\dots,\beta_d$ vary over all choices of non-negative integers   with $\sum \beta_i=Q$. 
  We claim that when all   elements of the form
$$\eta\left((x_{a_1}^*\wedge\dots \wedge x_{a_P}^*)(\omega)\otimes {x_1^*}^{(\beta_1)}\dots {x_d^*}^{(\beta_d)}\right)\tag\tnum{OM}$$are written in terms of the basis elements $\{k_{\pmb a;\pmb b}\}$ of $K_{d-P-1,Q-1}$, as given in (\tref{kpb}), then any given basis element $k_{\pmb a;\pmb b}$ appears {\bf at most once}.

\medskip\demo{Proof of Claim {\rm\tref{CC}}} Consider $\beta_1,\dots,\beta_d$, with $\sum \beta_i=Q$. Let $b_1$ be the least index with $\beta_{b_1}\neq 0$. 

\medskip
\flushpar If $b_1<a_1$, then the expression (\tref{OM}) is the basis element $k_{\pmb a;\pmb b}$ for $$\pmb a\text{ equal to }a_1<\dots< a_P\text{ and }\pmb b\text{ equal to } \underbrace{b_1,\ldots,b_1}_{\beta_{b_1}},\cdots,\underbrace{d,\ldots,d}_{\beta_d}.\tag\tnum{ans1}$$  

\medskip
\flushpar If $a_1\le b_1$, then the expression (\tref{OM}) involves the basis elements
$k_{\pmb a;\pmb b}$ for:
$$\left\{\left. \matrix\format\l\\\text{$\pmb a$ equal to $j,a_2,\dots,a_P$ written in strictly ascending order}\\\vspace{5pt}\text{$\pmb b$   equal to 
$a_1,\underbrace{b_1,\ldots,b_1}_{\beta_{b_1}},\dots \underbrace{j,\ldots,j}_{\beta_j-1},\cdots,\underbrace{d,\ldots,d}_{\beta_d}$}\endmatrix\right\vert\matrix
a_1<j\\\text{and}\\1\le \beta_j\endmatrix 
\right\}.\tag\tnum{ans2}$$  
Now one notices that given $k_{\pmb a;\pmb b}$ as described in (\tref{ans1}) or (\tref{ans2}), one can recreate the unique $d$-tuple $(\beta_1,\dots,\beta_d)$ so that (\tref{OM}) involves the basis element $k_{\pmb a;\pmb b}$. This completes the proof of Claim \tref{CC}. \enddemo

We show that the matrix for the vertical map 
$$1\otimes \pmb p_{r}^\Phi\: \widetilde{R}\otimes_{R_0} L_{p,r}\to \widetilde{R}\otimes_{R_0} K_{p,2n-2-r}$$
has monomial entries. (In fact, each non-zero entry in the matrix is plus or minus a variable from $\widetilde{R}$.) In this calculation, $0\le p\le d-1$ and $1\le r\le 2n-2$. Let 
$$\ell=\kappa(x_{a_1}\wedge\dots\wedge x_{a_{p+1}}\otimes x_1^{\beta_1}\dots x_d^{\beta_d})$$ be a basis element in $L_{p,r}$; so, in particular, $a_1<\dots<a_{p+1}$ and $\sum \beta_i=r-1$.
(Let $b_1$ be the least index with $\beta_{b_1}\neq 0$. The inequality $a_1\le b_1$ also holds; but this inequality will play no role in the present calculation.) It is easy to see that
$$(1 \otimes \pmb p_{r}^\Phi)\circ \kappa =\eta\circ (1\otimes \pmb p_{r-1}^\Phi).$$
Recall that $$\Phi=\sum_{\sum C_i=2n-2} t_{x_1^{C_1}\cdots x_d^{C_d}} {x_1^*}^{(C_1)}\cdots {x_d^*}^{(C_d)},$$where the sum is taken over all non-integers $C_1,\dots C_d$ with $\sum C_i=2n-2$. Of course, $t_M$ is a monomial (indeed, even a variable) in $\widetilde{R}$  for all monomials ${M=x_1^{C_1}\cdots x_d^{C_d}}$ of degree $2n-2$. It follows that 
$$\allowdisplaybreaks \alignat1 &(1 \otimes \pmb p_{r}^\Phi)(\ell)=(1 \otimes \pmb p_{r}^\Phi) \left(\kappa(x_{a_1}\wedge\dots\wedge x_{a_{p+1}}\otimes x_1^{\beta_1}\dots x_d^{\beta_d})\right)\\&
=\eta\left( (1\otimes \pmb p_{r-1}^\Phi)(x_{a_1}\wedge\dots\wedge x_{a_{p+1}}\otimes x_1^{\beta_1}\dots x_d^{\beta_d})\right)\\&
=\eta\left( x_{a_1}\wedge\dots\wedge x_{a_{p+1}}\otimes (x_1^{\beta_1}\dots x_d^{\beta_d})(\Phi)\right)\\&
=\eta\left( x_{a_1}\wedge\dots\wedge x_{a_{p+1}}\otimes 
\sum_{\sum C_i=2n-2\atop{0\le C_i-\beta_i}} t_{x_1^{C_1}\cdots x_d^{C_d}} {x_1^*}^{(C_1-\beta_1)}\cdots {x_d^*}^{(C_d-\beta_d)}\right)
\\&
=\sum_{\sum c_i=2n-1-r} t_{x_1^{\beta_1+c_1}\dots x_d^{\beta_d+c_d}}\eta\left( x_{a_1}\wedge\dots\wedge x_{a_{p+1}}\otimes{x_1^*}^{(c_1)}\cdot \dots \cdot {x_d^*}^{(c_{d})} \right),\endalignat$$where the most recent sum is taken over all non-negative integers  $c_1,\dots,c_d$ with $\sum c_i=2n-1-r$. Let $A_1<\dots<A_{d-p-1}$ be the complement of $a_1<\dots<a_{p+1}$ in $\{1,\dots,d\}$. Observe that 
$$x_{a_1}\wedge\dots\wedge x_{a_{p+1}}=\pm (x_{A_1}^*\wedge\dots\wedge x_{A_{d-p-1}}^*)(\omega).$$
We have shown that the vertical map $1\otimes \pmb p_r^\Phi$ sends the basis vector $\ell$ to a sum of terms of the form 
$$\text{monomial}\cdot \eta\left((x_{A_1}^*\wedge\dots\wedge x_{A_{d-p-1}}^*)(\omega)\otimes{x_1^*}^{(c_1)}\cdot \dots \cdot {x_d^*}^{(c_{d})} \right),\tag\tnum{bstar}$$
with $A_1,\dots A_{d-p-1}$ fixed and $c_1,\dots,c_d$ allowed to vary under the constraint that $\sum c_i$ is fixed. Claim \tref{CC} shows that when the elements of (\tref{bstar}) are written in terms of the basis $\{k_{\pmb a;\pmb b}\}$ of $K_{p,2n-2-r}$, then any given basis element $k_{\pmb a;\pmb b}$ appears at most once. In other words, the vertical map is a monomial map: each non-zero entry in the resulting matrix is plus or minus a variable. 

As promised, we now consider the horizontal maps. The maps
$$\widetilde{R}\otimes_{R_0} L_{0,r}\to R_0\quad\text{and}\quad \widetilde{R}\otimes_{R_0} {\tsize\bigwedge}^dU\to \widetilde{R}\otimes K_{d-1,2n-2-r}$$
merely list all of the monomials in $x_1,\dots,x_d$ of degree $r$ and degree $2n-1-r$, respectively. The map 
$\widetilde{R}\otimes_{R_0} L_{p,r}\to\widetilde{R}\otimes_{R_0} L_{p-1,r}$ sends the basis element 
$$\ell_{\pmb a;\pmb b} =\kappa(x_{a_1}\wedge\cdots \wedge x_{a_{p+1}}\otimes x_{b_1}\cdots x_{b_{r-1}})$$ 
of $L_{p,r}$ (with $a_1<\dots<a_{p+1}$, $b_1\le \dots\le b_{r-1}$, and $a_1\le b_1$) to $A+B$, with $A= x_{a_1}\cdot \kappa(x_{a_2}\wedge\cdots \wedge x_{a_{p+1}}\otimes x_{b_1}\cdots x_{b_{r-1}})$ and 
$$B=\sum_{i=2}^{p+1}(-1)^i
x_{a_i}\cdot \kappa(x_{a_1}\wedge\cdots \widehat{x_{a_i}}\cdots\wedge x_{a_{p+1}}\otimes x_{b_1}\cdots x_{b_{r-1}})$$
The sum $B$ is already written in terms of basis elements of $L_{p-1,r}$. If necessary, one can use the technique of (\tref{****'}) to write $A$ in terms of basis elements. It is not difficult to see that the basis elements needed for $A$ are distinct from the basis elements used in $B$.    

Finally, we consider the map 
$\widetilde{R}\otimes_{R_0} K_{p,r}\to\widetilde{R}\otimes_{R_0} K_{p-1,r}$  applied to   the basis element 
$$k_{\pmb a;\pmb b}=\eta\left((x_{a_1}^*\wedge\dots\wedge x_{a_{d-p-1}}^*)(\omega)\otimes {x_1^*}^{(\beta_1)}\cdot \ldots\cdot  {x_d^*}^{(\beta_d)}\right)$$ 
of $K_{p,r}$ 
 with $a_1<\dots<a_{d-p-1}$, $b_1\le \dots\le b_{r+1}$,  and $b_1< a_1$   
for
$$(b_1,\dots,b_{r+1})=(\underbrace{1,\ldots,1}_{\beta_1},\underbrace{2,\ldots,2}_{\beta_2},\cdots,\underbrace{d,\ldots,d}_{\beta_d}). $$ The basis element $k_{\pmb a;\pmb b}$ is sent to $A+B$ with
$$\split A&{}=\sum_{i\le b_1} x_i\eta\left((x_i^*\wedge x_{a_1}^*\wedge\dots\wedge x_{a_{d-p-1}}^*)(\omega)\otimes {x_1^*}^{(\beta_1)}\cdot \ldots\cdot  {x_d^*}^{(\beta_d)}\right)\\
B&{}=\sum_{b_1<i} x_i\eta\left((x_i^*\wedge x_{a_1}^*\wedge\dots\wedge x_{a_{d-p-1}}^*)(\omega)\otimes {x_1^*}^{(\beta_1)}\cdot \ldots\cdot  {x_d^*}^{(\beta_d)}\right).\endsplit$$
The sum $B$ is already written in terms of basis elements of $K_{p-1,r}$. One can use the technique of (\tref{****}) to write $A$ in terms of basis elements. It is not difficult to see that the basis elements needed for $A$ are distinct from the basis elements used in $B$. 
\qed\enddemo

\example{Example \tnum{n=3}} Adopt Data \tref{data5} with $d=3$. Let $x,y,z$ be a basis for the free $R_0$-module $U$.  We exhibit the resolution $\widetilde{\Bbb G}(3)=\widetilde{\Bbb G}(R_0,U,3;3)$ over
$\widetilde{R}=R_0[z,y,z,\{t_{M}\}]$, as $M$ roams over the fifteen monomials of degree $4$ in $x,y,z$. The resolution $\widetilde{\Bbb G}(3)$ is the mapping cone of 
$$\eightpoint\matrix
 & &0&\to&\widetilde{R}\otimes_{R_0}L_{2,3}&\to&\widetilde{R}\otimes_{R_0}L_{1,3}&\to&\widetilde{R}\otimes_{R_0}L_{0,3}
&\hskip-1pt \to&\widetilde{R}\\
& &\downarrow & &\downarrow&&\downarrow&&\downarrow
 & & & \\
0&\to&\widetilde{R}\otimes_{R_0}\bigwedge^3_{R_0}U&\to&\widetilde{R}\otimes_{R_0}K_{2,1}&\to&\widetilde{R}\otimes_{R_0}K_{1,1}&\to&\widetilde{R}\otimes_{R_0}K_{0,1}.
\endmatrix$$
We record the matrices using the bases $\ell_{\pmb a;\pmb b}$ and $k_{\pmb a;\pmb b}$ of (\tref{lpb}) and (\tref{kpb}):
$$\matrix\format\l&\quad\l&\quad\l&\quad\l &\quad\l&\quad\l \\
  L_{2,3}       & L_{1,3}      & L_{0,3}    &K_{2,1}&K_{1,1}&K_{0,1}\\
\ell_{1,2,3;1,1}&\ell_{1,2;1,1}&\ell_{1;1,1}&k_{\underline{\phantom{x}};1,1}&k_{2;1,1}&k_{2,3;1,1}\\
\ell_{1,2,3;1,2}&\ell_{1,2;1,2}&\ell_{1;1,2}&k_{\underline{\phantom{x}};1,2}&k_{2;1,2}&k_{2,3;1,2}\\
\ell_{1,2,3;1,3}&\ell_{1,2;1,3}&\ell_{1;1,3}&k_{\underline{\phantom{x}};1,3}&k_{2;1,3}&k_{2,3;1,3}\\
\ell_{1,2,3;2,2}&\ell_{1,2;2,2}&\ell_{1;2,2}&k_{\underline{\phantom{x}};2,2}&k_{3;1,1}\\
\ell_{1,2,3;2,3}&\ell_{1,2;2,3}&\ell_{1;2,3}&k_{\underline{\phantom{x}};2,3}&k_{3;1,2}\\
\ell_{1,2,3;3,3}&\ell_{1,2;3,3}&\ell_{1;3,3}&k_{\underline{\phantom{x}};3,3}&k_{3;1,3}\\
                &\ell_{1,3;1,1}&\ell_{2;2,2}&                               &k_{3;2,2}\\
                &\ell_{1,3;1,2}&\ell_{2;2,3}&                               &k_{3;2,3}\\
                &\ell_{1,3;1,3}&\ell_{2;3,3}&                               \\
                &\ell_{1,3;2,2}&\ell_{3;3,3}\\
                &\ell_{1,3;2,3}\\
                &\ell_{1,3;3,3}\\
                &\ell_{2,3;2,2}\\
                &\ell_{2,3;2,3}\\
                &\ell_{2,3;3,3}\\
\endmatrix$$ We identify $x_1$ with $x$, $x_2$ with $y$, and $x_3$ with $z$. 
We take $x\wedge y\wedge z$ to be the basis for $\bigwedge^3_{R_0}U$.
The resolution $\widetilde{\Bbb G}(3)$ then is  the mapping cone of 
$$\eightpoint \matrix
 & &0&\to&\widetilde{R}(-5,0)^6&@>h_3>>&\widetilde{R}(-4,0)^{15}&@>h_2>>&\widetilde{R}(-3,0)^{10}
&@>h_1>>&\widetilde{R}\\
& &\downarrow & &@V v_3 VV@V v_2 VV @V v_1 VV
 & & & \\
0&\to&\widetilde{R}(-7,-1)&@>h_3'>>&\widetilde{R}(-5,-1)^6&@>h_2'>>&\widetilde{R}(-4,-1)^8&@>h_1'>>&\widetilde{R}(-3,-1)^3
\endmatrix$$
with 
$h_1=[\smallmatrix x^3,x^2y,x^2z,xy^2,xyz,xz^2,y^3,y^2z,yz^2,z^3\endsmallmatrix]$,
$$\eightpoint  h_2=\left[\smallmatrix 
-y& 0& 0& 0& 0& 0&-z& 0& 0& 0& 0& 0& 0& 0& 0\\
 x&-y& 0& 0& 0& 0& 0&-z& 0& 0& 0& 0& 0& 0& 0\\
 0& 0&-y& 0& 0& 0& x& 0&-z& 0& 0& 0& 0& 0& 0\\
 0& x& 0&-y& 0& 0& 0& 0& 0&-z& 0& 0& 0& 0& 0\\
 0& 0& x& 0&-y& 0& 0& x& 0& 0&-z& 0& 0& 0& 0\\
 0& 0& 0& 0& 0&-y& 0& 0& x& 0& 0&-z& 0& 0& 0\\
 0& 0& 0& x& 0& 0& 0& 0& 0& 0& 0& 0&-z& 0& 0\\
 0& 0& 0& 0& x& 0& 0& 0& 0& x& 0& 0& y&-z& 0\\
 0& 0& 0& 0& 0& x& 0& 0& 0& 0& x& 0& 0& y&-z\\
 0& 0& 0& 0& 0& 0& 0& 0& 0& 0& 0& x& 0& 0& y\endsmallmatrix\right],\ \ h_3=\left[\smallmatrix
 z& 0& 0& 0& 0& 0\\
 0& z& 0& 0& 0& 0\\
-x& 0& z& 0& 0& 0\\
 0& 0& 0& z& 0& 0\\
 0&-x& 0& 0& z& 0\\
 0& 0&-x& 0& 0& z\\
-y& 0& 0& 0& 0& 0\\
 x&-y& 0& 0& 0& 0\\
 0& 0&-y& 0& 0& 0\\
 0& x& 0&-y& 0& 0\\
 0& 0& x& 0&-y& 0\\
 0& 0& 0& 0& 0&-y\\
 0& 0& 0& x& 0& 0\\
 0& 0& 0& 0& x& 0\\
 0& 0& 0& 0& 0& x\endsmallmatrix\right],$$
$$v_1=\left[\smallmatrix 
-t_{x^4} &-t_{x^3y}  &-t_{x^3z}  &-t_{x^2y^2}&-t_{x^2yz}&-t_{x^2z^2}&-t_{xy^3}&-t_{xy^2z} &-t_{xyz^2} &-t_{xz^3}\\
-t_{x^3y}&-t_{x^2y^2}&-t_{x^2yz} &-t_{xy^3}  &-t_{xy^2z}&-t_{xyz^2} &-t_{y^4} &-t_{y^3z}  &-t_{y^2z^2}&-t_{yz^3}\\
-t_{x^3z}&-t_{x^2yz} &-t_{x^2z^2}&-t_{xy^2z} &-t_{xyz^2}&-t_{xz^3}  &-t_{y^3z}&-t_{y^2z^2}&-t_{yz^3}  &-t_{z^4}
\endsmallmatrix\right],$$
$v_2$ is $$\hskip-22pt \eightpoint\left[\smallmatrix 
0&0&0&0&0&0&-t_{x^4} &-t_{x^3y}  &-t_{x^3z}  &-t_{x^2y^2}&-t_{x^2yz}&-t_{x^2z^2}&-t_{xy^3}&-t_{xy^2z} &-t_{xyz^2}\\
0&0&0&0&0&0&-t_{x^3y}&-t_{x^2y^2}&-t_{x^2yz} &-t_{xy^3}  &-t_{xy^2z}&-t_{xyz^2} &-t_{y^4} &-t_{y^3z}  &-t_{y^2z^2}\\
0&0&0&0&0&0&-t_{x^3z}&-t_{x^2yz} &-t_{x^2z^2}&-t_{xy^2z} &-t_{xyz^2}&-t_{xz^3}  &-t_{y^3z}&-t_{y^2z^2}&-t_{yz^3}\\
t_{x^4} &t_{x^3y}  &t_{x^3z}  &t_{x^2y^2}&t_{x^2yz}&t_{x^2z^2}&0&0&0&0&0&0&-t_{xy^2z} &-t_{xyz^2}&-t_{xz^3}\\
t_{x^3y}&t_{x^2y^2}&t_{x^2yz} &t_{xy^3}  &t_{xy^2z}&t_{xyz^2} &0&0&0&0&0&0&-t_{y^3z}&-t_{y^2z^2}&-t_{yz^3}\\
t_{x^3z}&t_{x^2yz} &t_{x^2z^2}&t_{xy^2z} &t_{xyz^2}&t_{xz^3}&0&0&0&0&0&0&-t_{y^2z^2}&-t_{yz^3}  &-t_{z^4}\\
t_{x^2y^2}&t_{xy^3}&t_{xy^2z}&t_{y^4}&t_{y^3z}&t_{y^2z^2}&t_{x^2yz}&t_{xy^2z}&t_{xyz^2}&t_{y^3z}&t_{y^2z^2}&t_{yz^3}&0&0&0\\
t_{x^2yz}&t_{xy^2z}&t_{xyz^2}&t_{y^3z}&t_{y^2z^2}&t_{yz^3}&t_{x^2z^2}&t_{xyz^2}&t_{xz^3}&t_{y^2z^2}&t_{yz^3}&t_{z^4}&0&0&0\\
\endsmallmatrix\right],$$
$$v_3=\left[\smallmatrix 
t_{x^4}&t_{x^3y}&t_{x^3z}&t_{x^2y^2}&t_{x^2yz}&t_{x^2z^2}\\
t_{x^3y}&t_{x^2y^2}&t_{x^2yz}&t_{xy^3}&t_{xy^2z}&t_{xyz^2}\\
t_{x^3z}&t_{x^2yz}&t_{x^2z^2}&t_{xy^2z}&t_{xyz^2}&t_{xz^3}\\
t_{x^2y^2}&t_{xy^3}&t_{xy^2z}&t_{y^4}&t_{y^3z}&t_{y^2z^2}\\
t_{x^2yz}&t_{xy^2z}&t_{xyz^2}&t_{y^3z}&t_{y^2z^2}&t_{yz^3}\\
t_{x^2z^2}&t_{xyz^2}&t_{xz^3}&t_{y^2z^2}&t_{yz^3}&t_{z^4}\endsmallmatrix\right],$$
$$h_1'=\left[\smallmatrix 
-z& 0& x& y&-x& 0& 0& 0\\
 0&-z& 0& 0& y& 0&-x& 0\\
 0& 0&-z& 0& 0& y& 0&-x\endsmallmatrix\right], \ \ h_2'=\left[\smallmatrix 
 y&-x& 0& 0& 0& 0\\
 0& y& 0&-x& 0& 0\\
 0& 0& y& 0&-x& 0\\
 z& 0&-x& 0& 0& 0\\
 0& z& 0& 0&-x& 0\\
 0& 0& z& 0& 0&-x\\
 0& 0& 0& z&-y& 0\\
 0& 0& 0& 0& z&-y\endsmallmatrix\right],\ \ \text{and}\ \ h_3'=\left[\smallmatrix 
x^2\\xy\\xz\\y^2\\yz\\z^2\endsmallmatrix\right].  
$$ 
 \endexample  

\SectionNumber=\MR\tNumber=1
\heading Section \number\SectionNumber. \quad The minimal resolution.
\endheading

In Theorem \tref{ek-k-m} we turn the resolution $P\otimes_{\widetilde{R}}\widetilde{\Bbb G}(r)_{\pmb \delta}$ of assertion (2) from Corollary \tref{EK-K-2} into a minimal resolution (when $P$ is a polynomial ring over a field). Our calculations are made at the level of $\widetilde{R}$; and, when $r=n$, our calculations are coordinate-free. The explicit resolution $\widetilde{\Bbb G}'(2)_{\pmb \delta}$ with $d=3$ and $n=r=2$ is recorded as Example \tref{X}.

\definition{Data \tnum{data6}} Consider    $(R_0,U,n,r)$, where $R_0$ is a commutative Noetherian ring, 
$U$ is a free $R_0$-module of rank $d$,    $n$ is a  positive integer, and $r$ is an integer with $n\le r\le 2n-2$.   
Let  $\widetilde{R}=\widetilde{R}(R_0,U,n)$, $\widetilde{I}=\widetilde{I}(R_0,U,n)$, $\widetilde{J}=\widetilde{J}(R_0,U,n)$,  and $\widetilde{\Bbb G}(r)=\widetilde{\Bbb G}(R_0,U,n;r)$  be the generic ring, the generic ideals, and  the generic complexes of Definition {\rm\tref{19.xms}}. Write $L_{p,q}$ and $K_{p,q}$ for $L_{p,q}^{R_0}U$ and $K_{p,q}^{R_0}U$, respectively. Let $\Phi$ be the element of $\widetilde{R}\otimes_{R_0} D_{2n-2}^{R_0}(U^*)$ which is described in (\tref{Phi}) and (\tref{appr}),   $\pmb \delta$ be the element $\det T_{\Phi}$ of $\widetilde{R}$ as described in Corollary \tref{EK-K-g} and Remark \tref{knew},
  $\Psi\:\widetilde{R}\otimes_{R_0} U\to \widetilde{R}$ be the  multiplication map of  {\rm(\tref{Phi})} and $\widehat{\Psi}\:\widetilde{R}\otimes_{R_0} \operatorname{Sym}_{\bullet}^{R_0}U\to \widetilde{R}$ be the $\widetilde{R}$-algebra map induced by the $\widetilde{R}$-module homomorphism $\widetilde{R}\otimes_{R_0}U\to \widetilde{R}$.
\enddefinition

Retain Data \tref{data6}. Apply $\widetilde{R}_{\pmb \delta} \otimes_{\widetilde {R}}\underline{\phantom{X}}$ to the double complex (\tref{d}) to obtain the double complex \vphantom{\tnum{d'}}
$$ \smallmatrix
{\tsize(\tref{d'})}                                        & &\widetilde{R}_{\pmb \delta}\otimes_{R_0}L_{d-1,r}  &\to&\dots&\to&\widetilde{R}_{\pmb \delta}\otimes_{R_0}L_{p,r}      &\to&\dots&\to&\widetilde{R}_{\pmb \delta}\otimes_{R_0}L_{0,r} & \hskip-18.8pt \to&\hskip-7.8pt\widetilde{R}_{\pmb \delta}\\\vspace{7pt}
                                    &   &\downarrow&&     &   &\downarrow                                   &   &     &   &
\downarrow                               &
\\\vspace{7pt}
\widetilde{R}_{\pmb \delta}\otimes_{R_0}\bigwedge^d_{R_0}U&\to&\widetilde{R}_{\pmb \delta}\otimes_{R_0}K_{d-1,2n-2-r} &\to&\dots&\to&\widetilde{R}_{\pmb \delta}\otimes_{R_0}K_{p,2n-2-r}&\to&\dots&\to&\widetilde{R}_{\pmb \delta}\otimes_{R_0}K_{0,2n-2-r}.\\
\endsmallmatrix$$ We know from assertion (3) of Corollary \tref{EK-K-g}  that the mapping cone of (\tref{d'}) is the  resolution 
$\widetilde{\Bbb G}(r)_{\pmb \delta}$ 
of $\widetilde{R}_{\pmb \delta} /\widetilde{J}^{r-n}\widetilde{I}$ by free $\widetilde{R}_{\pmb \delta}$-modules. In this section we show that each of the vertical maps $$\CD \widetilde{R}_{\pmb \delta}\otimes_{R_0}L_{p,r}\\ @V 1\otimes \pmb p_r^{\Phi} VV\\ \widetilde{R}_{\pmb \delta}\otimes_{R_0}K_{p,2n-2-r}\endCD\tag\tnum{vmap}$$
in (\tref{d'}) splits and  we split each of these vertical maps from $\widetilde{\Bbb G}(r)_{\pmb \delta}$ thereby producing a smaller resolution $\widetilde{\Bbb G}'(r)_{\pmb \delta}$ 
of $\widetilde{R}_{\pmb \delta} /\widetilde{J}^{r-n}\widetilde{I}\widetilde{R}_{\pmb \delta}$.

 The $\widetilde{R}$-module homomorphisms  which constitute  $\widetilde{\Bbb G}'(r)$ are introduced in Definition \tref{MR}. The modules $X_{p,r}$ which form the bulk of $\widetilde{\Bbb G}'(r)$ may be found in Definition \tref{SONE}. Most of the  maps of $\widetilde{\Bbb G}'(r)$ are induced by maps of the form ``$\operatorname{Kos}^\Psi$''; these maps are shown to be legitimate in Observation \tref{2.2}. The final ingredient,  a map called $\frak L_r$, may be found in Definition \tref{Br} and in the remarks of  \tref{equiv}. The properties of $\widetilde{\Bbb G}'(r)$, and especially $\widetilde{\Bbb G}'(r)_{\pmb \delta}$, are established in Theorem \tref{ek-k-m}.

\proclaim{Claim \tnum{proj}} Adopt Data {\rm\tref{data6}}. For each index $p$, with ${0\le p\le d-1}$, the $\widetilde{R}_{\pmb \delta}$-module homomorphism  {\rm(\tref{vmap})} is surjective. 
 \endproclaim 
\demo{Proof} Observe first that  the diagram 
$$\eightpoint \CD \widetilde{R}_{\pmb\delta}\otimes_{R_0}{\tsize\bigwedge} ^{p+1}_{R_0}U\otimes_{R_0}\operatorname{Sym}_{r-1}^{R_0}U@>\kappa>> \widetilde{R}_{\pmb\delta}\otimes_{R_0}{\tsize\bigwedge} ^{p}_{R_0}U\otimes_{R_0}\operatorname{Sym}_{r}^{R_0}U\\@V 1\otimes \pmb p_{r-1}^{\Phi}VV @V1\otimes \pmb p_r^{\Phi}VV \\\widetilde{R}_{\pmb\delta}\otimes_{R_0}{\tsize\bigwedge} ^{p+1}_{R_0}U\otimes_{R_0}D_{2n-r-1}^{R_0}(U^*)@>\eta>> \widetilde{R}_{\pmb\delta}\otimes_{R_0}{\tsize\bigwedge} ^{p}_{R_0}U\otimes_{R_0}D_{2n-2-r}^{R_0}(U^*)\endCD\tag\tnum{CD}$$commutes. The parameter $r-1$ is guaranteed to be at least $n-1$; so Remark \tref{20.21'}.2 ensures that  the vertical maps in (\tref{CD}) are surjective.     The domain of (\tref{vmap}) is the image of $\kappa$ in (\tref{CD}), and the target of (\tref{vmap}) is the image of $\eta$ in (\tref{CD}). \qed
\enddemo

\definition{Definition \tnum{SONE}}  Adopt Data {\rm\tref{data6}}.  For each index $p$, with $0\le p\le d-1$, define the $\widetilde{R}$-module  
$X_{p,r}$ to be the kernel of 
$$\CD \widetilde{R}\otimes_{R_0}L_{p,r}\\ @V 1\otimes \pmb p_r^{\Phi} VV\\ \widetilde{R}\otimes_{R_0}K_{p,2n-2-r}.\endCD$$
 \enddefinition 

\remark{\bf Remarks \tnum{Xpr}} 

\flushpar 1. Notice that the $\widetilde{R}$-module $X_{p,r}$ is defined in a coordinate-free manner. 

\bigskip\flushpar 2. The $R_0$-module $K_{p,2n-2-r}$ is free; so it follows, from Claim \tref{proj}, that  
 $(X_{p,r})_{\pmb \delta}$ is a projective $\widetilde{R}_{\pmb\delta}$-module
of rank: $$\operatorname{rank} L_{p,r}U-\operatorname{rank}K_{p,2n-2-r}U.$$ Formulas for these ranks are given in (\tref{rank}) and (\tref{rank'}).

\bigskip\flushpar 3. The module $X_{d-1,n}$ is equal to zero because the diagram 
$$\eightpoint \CD \widetilde{R}\otimes_{R_0}{\tsize\bigwedge} ^{d}_{R_0}U\otimes_{R_0}\operatorname{Sym}_{n-1}^{R_0}U@>\kappa>> \widetilde{R} \otimes_{R_0}{\tsize\bigwedge} ^{ d-1}_{R_0}U\otimes_{R_0}\operatorname{Sym}_{n}^{R_0}U\\@V 1\otimes \pmb p_{n-1}^{\Phi}V\simeq V @V1\otimes \pmb p_n^{\Phi}VV \\\widetilde{R} \otimes_{R_0}{\tsize\bigwedge} ^{d}_{R_0}U\otimes_{R_0}D_{n-1}^{R_0}(U^*)@>\eta>> \widetilde{R}\otimes_{R_0}{\tsize\bigwedge} ^{d}_{R_0}U\otimes_{R_0}D_{n-2}^{R_0}(U^*)\endCD $$commutes.

\endremark

\proclaim{Observation \tnum{2.2}} Adopt Data   {\rm\tref{data6}}.   Then the map $$\operatorname{Kos}^\Psi\: \widetilde{R}\otimes_{R_0}{\tsize\bigwedge}^p_{R_0}U\to \widetilde{R}\otimes_{R_0}{\tsize\bigwedge}^{p-1}_{R_0}U,$$ from {\rm(\tref{long})}, induces an $\widetilde{R}$-module homomorphism  $\operatorname{Kos}^\Psi\:X_{p,r}\to X_{p-1,r}$.\endproclaim

\demo{Proof} The cube
\SelectTips{eu}{12}
$$\eightpoint \xymatrix{\bigwedge^p\otimes S_r\ar[dd]^{\pmb p_r^\Phi}\ar[rd]^{\operatorname{Kos}^{\Psi}}\ar[rr]^{\kappa}&&\bigwedge^{p-1}\otimes S_{r+1}\ar'[d][dd]^{\pmb p_{r+1}^\Phi}\ar[rd]^{\operatorname{Kos}^{\Psi}}\\
&\bigwedge^{p-1}\otimes S_r\ar[dd]^/10pt/{\pmb p_r^\Phi}\ar[rr]^/-25pt/{\kappa}&&\bigwedge^{p-2}\otimes S_{r+1}\ar[dd]^{\pmb p_{r+1}^\Phi}\\
\bigwedge^p\otimes D_{2n-2-r}\ar'[r][rr]^{\eta}\ar[rd]^{\operatorname{Kos}^{\Psi}}&&\bigwedge^{p-1}\otimes D_{2n-3-r}\ar[rd]^{\operatorname{Kos}^{\Psi}}\\
&\bigwedge^{p-1}\otimes D_{2n-2-r}\ar[rr]^{\eta}&&\bigwedge^{p-2}\otimes D_{2n-3-r}}$$
commutes, where ``$\bigwedge^a\otimes S_b$'' represents $\widetilde{R}\otimes_{R_0}\bigwedge_{R_0}^aU\otimes_{R_0} \operatorname{Sym}_b^{R_0}U$
and ``$\bigwedge^a\otimes D_b$'' represents $\widetilde{R}\otimes_{R_0}\bigwedge_{R_0}^aU\otimes_{R_0} D_b^{R_0}(U^*)$. The $\widetilde{R}$-module
$X_{p,r}$ is the set of elements in the back, left, top module which are sent to $0$ in the back face of the cube. The elements of $X_{p,r}$ are sent to elements in the front, left, top module which go to zero in the front face of the cube. Such elements are in $X_{p-1.r}$. 
\qed \enddemo

A snake-like map $\frak L_r\: \widetilde{R}\otimes_{R_0}\bigwedge^d_{R_0}U\to \widetilde{R}\otimes_{R_0}L_{d-2,r}$ from the module on the bottom left of (\tref{d}) to the second (from the left) non-zero module on the top complex of (\tref{d})
plays a very important role in the complex $\widetilde{\Bbb G}'(r)$. Recall that for each index $i$, with $n-1\le i\le 2n-2$, the $\widetilde{R}$-module homomorphism $$\pmb p_i^{\Phi}\: \widetilde{R}_{\pmb \delta}\otimes_{R_0}\operatorname{Sym}_i^{R_0}U \to \widetilde{R}_{\pmb \delta}\otimes_{R_0}D_{2n-2-i}^{R_0}(U^*)$$ is a surjection. A splitting map 
$$\sigma_{2n-2-i}:\widetilde{R}_{\pmb \delta}\otimes_{R_0}D_{2n-2-i}^{R_0}(U^*)\to \widetilde{R}_{\pmb \delta}\otimes_{R_0}\operatorname{Sym}_i^{R_0}U$$ has been named in Remark \tref{20.21'}.2.  Careful examination of Remark \tref{20.21'}.2 shows that the only denominator that occurs in the map $\sigma_{2n-2-i}$ is $\pmb \delta$; so, in fact, $$\pmb \delta\sigma_{2n-2-i}:\widetilde{R}\otimes_{R_0}D_{2n-2-i}^{R_0}(U^*)\to \widetilde{R}\otimes_{R_0}\operatorname{Sym}_i^{R_0}U$$is an $\widetilde{R}$-module homomorphism.

\definition{Definition \tnum{Br}} Retain Data   \tref{data6}. Define the map $$\tsize {\frak L_r\:\widetilde{R}\otimes_{R_0}\bigwedge^d_{R_0}U \to \widetilde{R}\otimes_{R_0}L_{d-2,r}}$$ to be  the following composition of $\widetilde{R}$-module homomorphisms:
$$ \eightpoint \CD
        @. \widetilde{R}\otimes_{R_0}L_{d-1,r}@> \operatorname{Kos}^{\Psi} >> \widetilde{R}\otimes_{R_0}L_{d-2,r}. 
\\ 
                                   @.@A\kappa A\simeq A
\\  @. \widetilde{R}\otimes_{R_0}{\tsize \bigwedge}^d_{R_0}U\otimes_{R_0}\operatorname{Sym}_{r-1}^{R_0}U
       \\                              @.@A 1\otimes1\otimes \pmb \delta\sigma_{2n-1-r}AA
\\ @.\widetilde{R}\otimes_{R_0}{\tsize \bigwedge}^d_{R_0}U\otimes_{R_0}D_{2n-1-r}^{R_0}(U^*)\\
                                    @. @A\eta^{-1}A\simeq A\\ 
\widetilde{R}\otimes_{R_0}\bigwedge^d_{R_0}U @>\text{(\tref{lmm})} >> \widetilde{R}\otimes_{R_0}K_{d-1,2n-2-r} 
\endCD\tag\tnum{snake}$$ \enddefinition 

\remark{\bf Remarks \tnum{equiv}} 
\nopagebreak

\flushpar 1. It is clear from the definition of $\frak L_r$ that the image of $\frak L_r$ is actually contained in $X_{d-2,r}$. Indeed, $X_{d-2,r}$ is the kernel of the vertical map in (\tref{d}) 
emanating from  $\widetilde{R}\otimes_{R_0}L_{d-2,r}$ and (\tref{d}) is a map of complexes. 

\bigskip \flushpar 2. In practice, $$\tsize \frak L_r\:\widetilde{R}\otimes_{R_0}\bigwedge^d_{R_0}U \to X_{d-2,r}$$ does not really involve $\eta^{-1}$ because the map (\tref{lmm}) is the composition
$$\eightpoint \matrix\format\l\\ \tsize \widetilde{R}\otimes_{R_0}\bigwedge^d_{R_0}U\otimes_{R_0}R_0@> 1\otimes1\otimes\operatorname{ev}^* >> 
\widetilde{R}\otimes_{R_0}\bigwedge^d_{R_0}U\otimes_{R_0}\operatorname{Sym}_{2n-1-r}^{R_0}U\otimes_{R_0} D_{2n-1-r}^{R_0}(U^*)\\\vspace{5pt}
 \tsize @> \widehat{\Psi}>> \widetilde{R}\otimes_{R_0}\bigwedge^d_{R_0}U\otimes_{R_0} D_{2n-1-r}^{R_0}(U^*)@>\eta >> 
\widetilde{R}\otimes_{R_0}{\tsize \bigwedge}^d_{R_0}U\otimes_{R_0}D_{2n-2-r}^{R_0}(U^*),\endmatrix\tag\tnum{snake'} $$so there is no reason to compute $\eta^{-1}\circ \eta$. The map   $\operatorname{ev}^*$ is described in   (\tref{ev*}).

\bigskip \flushpar 3. It is important to notice that $\frak L_n$   can be made an equivariant map because $\sigma_{n-1}$ is the uniquely determined inverse of $$\pmb p_{n-1}^\Phi\: \widetilde{R}_{\pmb \delta}\otimes_{R_0}\operatorname{Sym}_{n-1}^{R_0}U \to \widetilde{R}_{\pmb \delta}\otimes_{R_0}D_{n-1}^{R_0}(U^*);$$ furthermore,  $\frak L_n$   can be defined over  $\widetilde{R}$  in a polynomial and equivariant manner by using the classical adjoint of $\pmb p_{n-1}^\Phi$ in place of its inverse. We do not know an equivariant description of $\frak L_i$ for $n+1\le i\le 2n-2$. An equivariant description  of the classical adjoint in the present situation is a little tricky; so we will record the notion completely at this point. Usually,
we will be a little cavalier  on this issue. An interested reader should always think of the domain of $\frak L_n$ as  $$\tsize \widetilde{R}\otimes_{R_0}\bigwedge^d_{R_0}U  \otimes_{R_0} \left(\bigwedge^{\text {top}}_{R_0}(D_{n-1}^{R_0}(U^*))\right)^{\otimes 2}.\tag\tnum{inter}$$
 The classical adjoint of $\pmb p_{n-1}^\Phi\: \operatorname{Sym}_{n-1}^{R_0}U\to D_{n-1}^{R_0}(U^*)$ is the map 
$$\tsize \operatorname{Adj}(\pmb p_{n-1}^\Phi)\: \operatorname{Sym}_{n-1}^{R_0}U\otimes _{R_0}\left(\bigwedge^{\text {top}}_{R_0}(D_{n-1}^{R_0}(U^*))\right)^{\otimes 2}\to D_{n-1}^{R_0}(U^*),$$
with $$\tsize [\operatorname{Adj}(\pmb p_{n-1}^\Phi)](u_{n-1}\otimes \Theta_1\otimes \Theta_2) = [(\bigwedge^{{\text {top}}-1}\pmb p_{n-1}^\Phi)(u_{n-1}(\Theta_1))](\Theta_2),$$for   $u_{n-1}\in \operatorname{Sym}_{n-1}^{R_0}U$  and $\Theta_i\in \bigwedge^{\text{top}}_{R_0}(D_{n-1}^{R_0}(U^*))$, where ``$\text{top}$'' is equal to the rank of $D_{n-1}^{R_0}(U^*)$. That is, ``top'' is equal to $\binom{n+d-2}{n-1}$.
One easily checks that 
$$\split \tsize \pmb p_{n-1}^\Phi([\operatorname{Adj}(\pmb p_{n-1}^\Phi)](u_{n-1}\otimes \Theta_1\otimes \Theta_2))&\tsize{}= 
[(\bigwedge^{{\text{top}}}\pmb p_{n-1}^\Phi)(\Theta_1)](\Theta_2)\cdot u_{n-1}\\
 \tsize ([\operatorname{Adj}(\pmb p_{n-1}^\Phi)]((\pmb p_{n-1}^\Phi)^*(w_{n-1})\otimes \Theta_1\otimes \Theta_2))&\tsize{}= 
[(\bigwedge^{{\text{top}}}\pmb p_{n-1}^\Phi)(\Theta_1)](\Theta_2)\cdot w_{n-1}\endsplit$$
for   $w_{n-1}\in D_{n-1}^{R_0}(U^*)$. Thus, $\operatorname{Adj}(\pmb p_{n-1}^\Phi)\circ  \pmb p_{n-1}^\Phi\:\operatorname{Sym}_{n-1}^{R_0}U \to \operatorname{Sym}_{n-1}^{R_0}U $ and  $\pmb p_{n-1}^\Phi \circ \operatorname{Adj}(\pmb p_{n-1}^\Phi)\: D_{n-1}^{R_0}(U^*)\to D_{n-1}^{R_0}(U^*)$
 are both equal to multiplication by the determinant of $\pmb p_{n-1}^\Phi$.

\endremark

\definition{Definition \tnum{MR}} Adopt Data  \tref{data6}. Let $\widetilde{\Bbb G}'(r)$ be the following sequence of $\widetilde{R}$-modules and  $\widetilde{R}$-module homomorphisms: 
$$\eightpoint 0\to \matrix X_{d-1,r}\\\oplus\\\widetilde{R}\otimes_{R_0}{\tsize \bigwedge}^d_{R_0}U\endmatrix @>[\smallmatrix \operatorname{Kos}^{\Psi}&\frak L_r\endsmallmatrix] >>X_{d-2,r}@> \operatorname{Kos}^{\Psi} >>X_{d-3,r}@> \operatorname{Kos}^{\Psi} >>\dots @> \operatorname{Kos}^{\Psi} >>X_{0,r}@> \widehat{\Psi}>>\widetilde{R}.$$ 
 Of course, according to Remark \tref{Xpr}.3, $\widetilde{\Bbb G}'(n)$ is equal to 
$$\eightpoint 0\to \widetilde{R}\otimes_{R_0}{\tsize \bigwedge}^d_{R_0}U@>\frak L_n >>X_{d-2,n}@> \operatorname{Kos}^{\Psi} >>X_{d-3,n}@> \operatorname{Kos}^{\Psi} >>\dots @> \operatorname{Kos}^{\Psi} >>X_{0,n}@> \widehat{\Psi}>>\widetilde{R}.$$ (The modification of (\tref{inter}) should be applied if one wants $\widetilde{\Bbb G}'(n)$ to be one hundred percent equivariant.)\enddefinition

\proclaim{Theorem \tnum{ek-k-m}}Adopt Data   {\rm\tref{data6}}.
\roster
\item The $\widetilde{R}$-module homomorphisms  $\widetilde{\Bbb G}'(r)$ of Definition {\rm\tref{MR}} form a complex.
\item The complex $\widetilde{\Bbb G}'(r)_{\pmb \delta}$ is a resolution of $\widetilde{R}_{\pmb \delta}/\widetilde{J}^{r-n}\widetilde{I}\widetilde{R}_{\pmb \delta}$ by projective  $\widetilde{R}_{\pmb \delta}$\hskip-.5pt-modules. 
\item If $\pmb k$  is a field and an $R_0$-algebra,  $P$ is the polynomial ring $\operatorname{Sym}_{\bullet}^{\pmb k}U$, 
$I=\operatorname{ann}\phi$ is an ideal of the set $\Bbb I_n(\pmb k,U)$ of Definition {\rm \tref{inr}}, 
$J$ is the ideal of $P$ generated by $\operatorname{Sym}_1^{\pmb k}U$, and $P$ is an $\widetilde{R}$-algebra by way of 
$$\widetilde{R}\to \pmb k\otimes_{R_0}\widetilde{R} @> \widehat{\phi}>> P,$$where
$\widehat{\phi}\:\pmb k\otimes_{R_0}\widetilde{R}\to P$ is the $P$-algebra homomorphism of {\rm(\tref{phihat})} and {\rm(\tref{???})}, then $P\otimes_{\widetilde{R}}\widetilde{\Bbb G}'(r)_{\pmb \delta}$ is a minimal homogeneous resolution of
$P/J^{r-n}I$ by free $P$-modules. 
\item In the situation of {\rm(3)}, with $r=n$, $P\otimes_{\widetilde{R}}\widetilde{\Bbb G}'(n)_{\pmb \delta}$ is an equivariant minimal homogeneous resolution of
$P/I$ by free $P$-modules. 
\item In the situation of {\rm(3)}, the graded Betti numbers of $P/J^{r-n}I$ are described by 
$$ 0\to \matrix P(-r-d+1)^{\beta_d}\\\oplus\\ P(-2n-d+2)\endmatrix\to P(-r-d+2)^{\beta_{d-1}}\to\dots\to P(-r-1)^{\beta_{2}}\to P(-r)^{\beta_{1}}\to P,$$with $\beta_i=\binom{d+r-1}{i-1+r}\binom{i+r-2}{i-1}-\binom{d+2n-r-2}{i-1}\binom{d+2n-r-i-2}{d-i}$. In particular, the strand with the exponents labeled by $\beta_i$ is linear. 
\endroster \endproclaim

\demo{Proof} Assertion (1) is clear from the construction of $\widetilde{\Bbb G}'(r)$. The point of this assertion is that  $\widetilde{\Bbb G}'(r)$ is a well-defined complex over $\widetilde{R}$. It is clear from Lemma \tref{q-i} that  $\widetilde{\Bbb G}'(r)_{\pmb \delta}$ is a well-defined complex over  $\widetilde{\Bbb G}'(r)_{\pmb \delta}$; but in fact, we carefully defined  $\widetilde{\Bbb G}'(r)$ to actually be a complex over $\widetilde{R}$, without any denominators.  Assertion (2) is a consequence of Lemma \tref{q-i} applied to the resolution 
 $\widetilde{\Bbb G}(r)_{\pmb \delta}$ of $\widetilde{R}_{\pmb \delta}/\widetilde{J}^{r-n}\widetilde{I}\widetilde{R}_{\pmb \delta}$, as found in Corollary \tref{EK-K-g}. Claim \tref{proj} establishes that the vertical maps of (\tref{d'}) are surjective. The targets of these maps are free; hence the maps split and the hypotheses of Lemma \tref{q-i} are satisfied.

\medskip\flushpar (3) and (5). We know from (2) that the complex $P\otimes_{\widetilde{R}}\widetilde{\Bbb G}'(r)_{\pmb \delta}$  is a homogeneous resolution of $P/J^{r-n}I$ by free $P$-modules. One readily reads that the graded Betti numbers of this complex are given in (5); see also Remark \tref{Xpr}.2. We conclude that $P\otimes_{\widetilde{R}}\widetilde{\Bbb G}'(r)$ is a minimal resolution.

\medskip\flushpar (4) All of the maps and modules of $\widetilde{\Bbb G}'(r)$ are obviously equivariant, except, possibly, $\frak L_r$. We have explained in Remark \tref{equiv}.3 how to make $\frak L_n$ become an equivariant map. \qed \enddemo

\proclaim{Lemma \tnum{q-i}}Let $R$ be a commutative Noetherian ring and let  $$\eightpoint\CD @.0@>>> T_d@> h_d>> T_{d-1}@>h_{d-1} >> \dots @>h_2>> T_1@> h_1>> T_0\\
@. @. @V v_d VV @V v_{d-1} VV @. @V v_1 VV\\
0@>>> B_d@>h_d'>>B_{d-1} @> h_{d-1}'>> B_{d-2}@>h_{d-2}'>>\dots  @>h_1'>> B_0.\endCD\tag\tnum{dd'}$$ be  a map of complexes of $R$-modules. Suppose that for each $i$, with $0\le i\le d-1$, $\sigma_i\:B_i\to T_{i+1}$ is an $R$-module homomorphism which splits $v_{i+1}\:T_{i+1}\to B_{i}$ in the sense that the composition $v_{i+1}\circ \sigma_i$ is the identity map on $B_i$. Then the following statements hold.
\roster
\item The $R$-module homomorphisms 
$$\eightpoint 0\to \matrix \operatorname{ker} v_d\\\oplus\\ B_d\endmatrix @>[\matrix v_d&h_d\circ \sigma_{d-1}\circ h_d'\endmatrix]>> \operatorname{ker} v_{d-1}@> h_{d-1}>> \dots @> h_{3}>>\operatorname{ker} v_2 @> h_2>> \operatorname{ker} v_1 @> h_1 >>T_0 \tag\tnum{s}$$ form a complex.

\item There is a quasi-isomorphism from the total complex of {\rm (\tref{dd'})} to the complex {\rm (\tref{s})}.
\endroster
\endproclaim

\demo{Proof} Assertion (1) is obvious as soon as one observes that $h_i(\operatorname{ker} v_i)$ is contained in  $\operatorname{ker} v_{i-1}$.
We prove assertion (2). (We follow ideas used in the proof of Theorem \tref{EK-K'}.) Let $(\Bbb M_{\bullet},m_{\bullet})$ be the mapping cone of (\tref{dd'}). Fix $i$ with $1\le i\le d$. Observe that $m_i\:\Bbb M_i\to \Bbb M_{i-1}$ is
$$\matrix T_i\\\oplus\\ B_i\endmatrix @>\left[\smallmatrix h_i&0\\v_i&-h_i'\endsmallmatrix\right] >> \matrix T_{i-1}\\\oplus\\ B_{i-1}.\endmatrix\tag\tnum{mi}$$
Take advantage of the direct sum decomposition $T_i=\operatorname{ker} v_i\oplus \operatorname{im} \sigma_{i-1}$, which is induced by the equation $$\theta=(1-\sigma_{i-1}\circ v_i)(\theta)+(\sigma_{i-1}\circ v_i)(\theta),$$ for all $\theta$ in $T_i$ to write (\tref{mi}) in the form
$$\eightpoint \hskip-80pt\CD {\smallmatrix \operatorname{ker} v_i\\\oplus \\ \operatorname{im} \sigma_{i-1}\\\oplus\\ B_i\endsmallmatrix} @>{\left[\smallmatrix (1-\sigma_{i-2}\circ v_{i-1})\circ h_i&
(1-\sigma_{i-2}\circ v_{i-1})\circ h_i&0\\
\sigma_{i-2}\circ v_{i-1}\circ h_i& \sigma_{i-2}\circ v_{i-1}\circ h_i&0\\
0&v_i&-h_i'\endsmallmatrix\right]} >> 
{\smallmatrix \operatorname{ker} v_{i-1}\\\oplus\\ \operatorname{im} \sigma_{i-2}\\\oplus\\ B_{i-1}\endsmallmatrix}\\
@V=VV @V =VV \\
 {\smallmatrix \operatorname{ker} v_i\\\oplus \\ \operatorname{im} \sigma_{i-1}\\\oplus\\ B_i\endsmallmatrix} @>{\left[\smallmatrix h_i&
h_i-\sigma_{i-2}h_{i-1}'\circ v_{i}&0\\
0& \sigma_{i-2}\circ h_{i-1}'\circ v_{i}&0\\
0&v_i&-h_i'\endsmallmatrix\right]} >> 
{\smallmatrix \operatorname{ker} v_{i-1}\\\oplus\\ \operatorname{im} \sigma_{i-2}\\\oplus\\ B_{i-1}\endsmallmatrix}\\
\endCD$$Use the isomorphisms $\sigma_j\:B_j\to\operatorname{im} \sigma_j$ and $v_{j+1}\:\operatorname{im} \sigma_j \to B_j$ to see that (\tref{mi}) is isomorphic to Figure \tref{Fig3}.
This calculation   works for $1\le i\le d-1$, but it must be modified for $i=d$, because in the isomorphism of Figure \tref{Fig3}, the image of $h_{i+1}\circ \sigma_{i}-\sigma_{i-1}\circ h_{i}'$ is contained in $\operatorname{ker}v_i$, provided $i\le d-1$, since $v_{d+1}\circ \sigma_d$ is not defined (or is zero); but it certainly is not the identity map on $B_d$. The modification when $i=d$ is 
$$\eightpoint \hskip-80pt\CD{\smallmatrix \operatorname{ker} v_d\\\oplus \\ B_{d-1}\\\oplus\\ B_d\endsmallmatrix} @>{\left[\smallmatrix  h_d&
h_d\circ \sigma_{d-1}-\sigma_{d-2}\circ h_{d-1}'&0\\
0& h_{d-1}'&0\\
0&1&-h_d'\endsmallmatrix\right]} >> {\smallmatrix \operatorname{ker} v_{d-1}\\\oplus\\ B_{d-2}\\\oplus\\ B_{d-1}\endsmallmatrix}\\@V\simeq V{\left[\smallmatrix 1&0&0\\0&1& -h_{d}'\\0&0&1\endsmallmatrix\right]} V @V \simeq V{\left[\smallmatrix 1&0&-(h_d\circ \sigma_{d-1}-\sigma_{d-2}\circ h_{d-1}')\\0&1& -h_{d-1}'\\0&0&1\endsmallmatrix\right]}V
\\
{\smallmatrix \operatorname{ker} v_d\\\oplus \\ B_{d-1}\\\oplus\\ B_d\endsmallmatrix} @>{\left[\smallmatrix  
h_d&0&h_d \circ \sigma_{d-1}\circ h_d'\\
0&0&0\\
0&1&0\endsmallmatrix\right]} >> {\smallmatrix \operatorname{ker} v_{d-1}\\\oplus\\ B_{d-2}\\\oplus\\ B_{d-1}.\endsmallmatrix}
\endCD$$It is now clear that (\tref{s}) is obtained from $\Bbb M$ by splitting off    the extraneous summands. \qed
\topinsert
$$\eightpoint \hskip-80pt\CD
 \\{\smallmatrix \operatorname{ker} v_i\\\oplus \\ B_{i-1}\\\oplus\\ B_i\endsmallmatrix} @>
{\left[\smallmatrix h_i&
h_i-\sigma_{i-2}h_{i-1}'\circ v_{i}\circ \sigma_{i-1}&0\\
0& v_{i-1}\circ\sigma_{i-2}\circ h_{i-1}'\circ v_{i}\circ \sigma_{i-1}&0\\
0&v_i\circ \sigma_{i-1}&-h_i'\endsmallmatrix\right]}
>> {\smallmatrix \operatorname{ker} v_{i-1}\\\oplus\\ B_{i-2}\\\oplus\\ B_{i-1}\endsmallmatrix}\\
@V= VV @V = VV
\\
{\smallmatrix \operatorname{ker} v_i\\\oplus \\ B_{i-1}\\\oplus\\ B_i\endsmallmatrix} @>{\left[\smallmatrix  h_i&
h_i\circ \sigma_{i-1}-\sigma_{i-2}\circ h_{i-1}'&0\\
0& h_{i-1}'&0\\
0&1&-h_i'\endsmallmatrix\right]} >> {\smallmatrix \operatorname{ker} v_{i-1}\\\oplus\\ B_{i-2}\\\oplus\\ B_{i-1}\endsmallmatrix}\\@V\simeq V{\left[\smallmatrix 1&0&-(h_{i+1}\circ \sigma_{i}-\sigma_{i-1}\circ h_{i}')\\0&1& -h_{i}'\\0&0&1\endsmallmatrix\right]} V @V \simeq V{\left[\smallmatrix 1&0&-(h_i\circ \sigma_{i-1}-\sigma_{i-2}\circ h_{i-1}')\\0&1& -h_{i-1}'\\0&0&1\endsmallmatrix\right]}V
\\
{\smallmatrix \operatorname{ker} v_i\\\oplus \\ B_{i-1}\\\oplus\\ B_i\endsmallmatrix} @>{\left[\smallmatrix  
h_i&0&0\\
0&0&0\\
0&1&0\endsmallmatrix\right]} >> {\smallmatrix \operatorname{ker} v_{i-1}\\\oplus\\ B_{i-2}\\\oplus\\ B_{i-1}.\endsmallmatrix}
\endCD$$
{\bf Figure \tnum{Fig3}.} {\smc A complex isomorphic to (\tref{mi}). This complex is used in the proof of Lemma \tref{q-i}.}
\endinsert\enddemo

\example{Example \tnum{X}} Adopt Data \tref{data6} with $R_0=\Bbb Z$, $d=3$, and $n=r=2$. In this example we record and explain  the generic resolution $\widetilde{\Bbb G}'(2)_{\pmb\delta}$ of Theorem \tref{ek-k-m} when $(d,n,r)=(3,2,2)$. Theorem \tref{ek-k-m} only promises a resolution by projective $\widetilde{R}_{\pmb \delta}$-modules; however, in fact, this is a resolution by {\bf free} $\widetilde{R}_{\pmb \delta}$-modules. 
Recall that    $\widetilde{R}$ is the polynomial ring $\widetilde{R}=\Bbb Z[x,y,z,t_{x^2},t_{xy},t_{xz},t_{y^2},t_{yz},t_{z^2}]$, $T$ is the matrix $$T=\bmatrix t_{x^2}&t_{xy}&t_{xz}\\t_{xy}&t_{y^2}&t_{yz}\\t_{xz}&t_{yz}&t_{z^2}\endbmatrix,$$and $\pmb \delta=\det T$. Let $Q$ be the classical adjoint of $T$. We note for future reference that $Q$ is a symmetric matrix and $$TQ=QT=\pmb\delta I_3,$$ where $I_r$ represents the $r\times r$ identity matrix.
Define the elements $\lambda_1,\lambda_2,\lambda_3$ of $\widetilde{R}$ by
$$[\lambda_1,\lambda_2,\lambda_3]=[x,y,z]Q\quad\text{or}\quad Q\bmatrix x\\y\\z\endbmatrix = \bmatrix \lambda_1\\\lambda_2 \\\lambda_3\endbmatrix.$$
 View $\widetilde{R}_{\pmb \delta}$ as a graded ring where $x,y,z$ all have degree one and the $t$'s have degree zero. 

\proclaim{Claim \tnum{CLAIM}} The generic resolution $\widetilde{\Bbb G}'(2)_{\pmb\delta}$ is equal to 
$$0\to \widetilde{R}_{\pmb\delta}(-5)@>d_3>> \widetilde{R}_{\pmb\delta}(-3)^5@>d_2>> \widetilde{R}_{\pmb\delta}(-2)^5@>d_1>> \widetilde{R}_{\pmb\delta},$$with
$$d_1=d_3^{\text{\rm T}}=\bmatrix -t_{z^2}x\lambda_1+\pmb\delta z^2&-t_{y^2}x\lambda_1+\pmb \delta y^2&-t_{yz}x\lambda_1+\pmb \delta yz&x\lambda_2&x\lambda_3\endbmatrix,\tag\tnum{d1}$$
and $d_2$ is the alternating matrix which is obtaining by completing:
$$\eightpoint \bmatrix 0&-xQ_{2,3}&-xQ_{3,3}&-xt_{y^2}Q_{1,3}&-xt_{yz}Q_{3,1}+\pmb\delta y\\
*&0&xQ_{2,2}&xt_{yz}Q_{2,1}-\pmb \delta z&xt_{z^2}Q_{2,1}\\
*&*&0&xt_{yz}Q_{3,1}-xt_{y^2}Q_{2,1}+\pmb\delta y&xt_{z^2}Q_{3,1}-xt_{yz}Q_{2,1}-\pmb \delta z\\
*&*&*&0&xQ_{1,1}^2\\
*&*&*&*&0\endbmatrix.\tag\tnum{d2}$$\endproclaim 

\remark{Remark}We  write $d_3^{\text{\rm T}}$ to mean the transpose of the matrix $d_3$. In the present formulation, the matrix $d_1$ is equal to $-\pmb \delta$ times the row vector of signed maximal order Pfaffians of $d_2$.\endremark 

\demo{Proof of Claim {\rm\tref{CLAIM}}} The resolution $\widetilde{\Bbb G}'(2)_{\pmb\delta}$ is obtained using our techniques as described in Example \tref{n=3}. The resolution $\widetilde{\Bbb G}(2)$ is the mapping cone of 
$$\eightpoint\matrix
 & &0&\to&\widetilde{R}\otimes_{\Bbb Z}L_{2,2}&\to&\widetilde{R}\otimes_{\Bbb Z}L_{1,2}&\to&\widetilde{R}\otimes_{\Bbb Z}L_{0,2}
&\hskip-1pt \to&\widetilde{R}\\
& &\downarrow & &\downarrow&&\downarrow&&\downarrow
 & & & \\
0&\to&\widetilde{R}\otimes_{\Bbb Z}\bigwedge^3_{\Bbb Z}U&\to&\widetilde{R}\otimes_{\Bbb Z}K_{2,0}&\to&\widetilde{R}\otimes_{\Bbb Z}K_{1,0}&\to&\widetilde{R}\otimes_{\Bbb Z}K_{0,0}.
\endmatrix$$
We record the matrices of 
$\widetilde{\Bbb G}(2)$
using the bases $\ell_{\pmb a;\pmb b}$ and $k_{\pmb a;\pmb b}$ of (\tref{lpb}) and (\tref{kpb}):
$$\matrix\format\l&\quad\l&\quad\l&\quad\l &\quad\l&\quad\l \\
  L_{2,2}       & L_{1,2}      & L_{0,2}    &K_{2,0}&K_{1,0}&K_{0,0}\\
\ell_{1,2,3;1}&\ell_{1,2;1}&\ell_{1;1}&k_{\underline{\phantom{x}};1}&k_{2;1}&k_{2,3;1}\\
\ell_{1,2,3;2}&\ell_{1,2;2}&\ell_{1;2}&k_{\underline{\phantom{x}};2}&k_{3;1}&\\
\ell_{1,2,3;3}&\ell_{1,2;3}&\ell_{1;3}&k_{\underline{\phantom{x}};3}&k_{3;2}&\\
              &\ell_{1,3;1}&\ell_{2;2}& &\\
              &\ell_{1,3;2}&\ell_{2;3}& &\\
              &\ell_{1,3;3}&\ell_{3;3}& &\\
                &\ell_{2,3;2}& &                               &\\
                &\ell_{2,3;3}& &                               &
\endmatrix\tag\tnum{basis}$$ We identify $x_1$ with $x$, $x_2$ with $y$, and $x_3$ with $z$. 
We take $\omega=x\wedge y\wedge z$ to be the basis for $\bigwedge^3_{\Bbb Z}U$.
The resolution $\widetilde{\Bbb G}(2)$ then is  the mapping cone of 
$$\eightpoint \matrix
 & &0&\to&\widetilde{R}(-4)^3&@>h_3>>&\widetilde{R}(-3)^{8}&@>h_2>>&\widetilde{R}(-2)^{6}
&@>h_1>>&\widetilde{R}\\
& &\downarrow & &@V v_3 VV@V v_2 VV @V v_1 VV
 & & & \\
0&\to&\widetilde{R}(-5)&@>h_3'>>&\widetilde{R}(-4)^3&@>h_2'>>&\widetilde{R}(-3)^3&@>h_1'>>&\widetilde{R}(-2)
\endmatrix$$
with $$\allowdisplaybreaks\eightpoint \alignat1 &h_1=[\smallmatrix x^2,xy,xz,y^2,yz,z^2\endsmallmatrix]\ \  h_2=\left[\smallmatrix 
-y& 0& 0&-z& 0& 0& 0& 0\\
 x&-y& 0& 0&-z& 0& 0& 0\\
 0& 0&-y& x& 0&-z& 0& 0\\
 0& x& 0& 0& 0& 0&-z& 0\\
 0& 0& x& 0& x& 0& y&-z\\
 0& 0& 0& 0& 0& x& 0&  y\endsmallmatrix\right],\ \ h_3=\left[\smallmatrix
 z& 0& 0\\
 0& z& 0\\
-x& 0& z\\
-y& 0& 0\\
 x&-y& 0\\
 0& 0&-y\\
 0& x& 0\\
 0& 0& x\endsmallmatrix\right],
\\&
v_1=\left[\smallmatrix t_{x^2},t_{xy},t_{xz},t_{y^2},t_{yz},t_{z^2}
\endsmallmatrix\right],\ \ 
 v_2=\left[\smallmatrix 
0&0&0&-t_{x^2} &-t_{xy}  &-t_{xz}  &-t_{y^2}&-t_{yz}\\
t_{x^2} &t_{xy}  &t_{xz} &0&0&0&-t_{yz} &-t_{z^2}\\
t_{xy}&t_{y^2}&t_{yz} &t_{xz}  &t_{yz}&t_{z^2} &0&0\endsmallmatrix\right],\ \ v_3=T, 
\\&
h_1'=\left[\smallmatrix -z&y&-x\endsmallmatrix\right], \ \ h_2'=\left[\smallmatrix 
 y&-x& 0\\
 z& 0&-x&\\
 0& z&-y\endsmallmatrix\right],\ \ \text{and}\ \ h_3'=\left[\smallmatrix 
x\\y\\z\endsmallmatrix\right].\endalignat $$ According to Definition \tref{MR}, the resolution $\widetilde{\Bbb G}'(2)_{\pmb\delta}$ 
has the form
$$0\to  \widetilde{R}_{\pmb\delta}\otimes_{\Bbb Z}{\tsize \bigwedge}^3_{\Bbb Z}U@>\frak L_2 >>(X_{1,2})_{\pmb\delta}@> \operatorname{Kos}^{\Psi} >>(X_{0,2})_{\pmb\delta}@> \widehat{\Psi}>>\widetilde{R}_{\pmb\delta}.$$ 
The key step in this proof is that we are able to  identify bases for $(X_{0,2})_{\pmb\delta}=(\operatorname{ker} v_1)_{\pmb\delta}$ and $(X_{1,2})_{\pmb\delta}=(\operatorname{ker} v_2)_{\pmb\delta}$. Indeed, the matrices
$$J_1=\bmatrix Q&0\\0&I_3\endbmatrix\quad\text{and}\quad J_2=\bmatrix Q&0&0\\0&Q&0\\0&0&I_2\endbmatrix$$ are invertible over $\widetilde{R}_{\pmb\delta}$; it is easy to read the kernels of the matrices
$$\eightpoint v_1 J_1=\bmatrix\pmb \delta&0&0&t_{y^2}&t_{yz}&t_{z^2}\endbmatrix\quad\text{and}\quad v_2 J_2=\bmatrix 0&0&0&-\pmb \delta&0&0&-t_{y^2}&-t_{yz}\\\pmb \delta&0&0&0&0&0&-t_{yz}&-t_{z^2}\\
0&\pmb \delta&0&0&0&\pmb\delta&0&0\endbmatrix;$$and, in particular, we conclude that $(X_{0,2})_{\pmb\delta}$ and $(X_{1,2})_{\pmb\delta}$ are both free $\widetilde{R}_{\pmb\delta}$-modules and the columns of 
$$\eightpoint B=\bmatrix -t_{z^2}Q_{*,1}&-t_{y^2}Q_{*,1}&-t_{yz}Q_{*,1}&Q_{*,2}&Q_{*,3}\\
0&\pmb\delta&0&0&0\\
0&0&\pmb\delta&0&0\\
\pmb\delta&0&0&0&0\endbmatrix \ \text{and}\tag\tnum{blah}$$$$\eightpoint B'=\bmatrix Q_{*,3}&0&Q_{*,2}&t_{yz}Q_{*,1}&t_{z^2}Q_{*,1}\\0&-Q_{*,2}
&-Q_{*,3}&-t_{y^2}Q_{*,1}&-t_{yz}Q_{*,1}\\0&0&0&\pmb\delta  &0\\0&0&0&0&\pmb\delta\endbmatrix\tag\tnum{blah'}$$
represent bases for $(X_{0,2})_{\pmb\delta}$ and $(X_{1,2})_{\pmb\delta}$, respectively. (We write $Q_{*,j}$ to represent column $j$ of the matrix $Q$.)  In other words, we have proven that 
$$\matrix \format \l&\quad\quad\l\\
g_1=-t_{z^2}\sum\limits_{j=1}^3 Q_{j,1}\ell_{1;j}+\pmb \delta\ell_{3;3}&
g_2=-t_{y^2}\sum\limits_{j=1}^3 Q_{j,1}\ell_{1;j}+\pmb \delta \ell_{2;2}\\
g_3=-t_{yz}\sum\limits_{j=1}^3 Q_{j,1}\ell_{1;j}+\pmb  \delta\ell_{2;3}&
g_4=\sum\limits_{j=1}^3 Q_{j,2}\ell_{1;j}\\
g_5=\sum\limits_{j=1}^3 Q_{j,3}\ell_{1;j}\quad\text{and}\\
\endmatrix$$
$$\matrix \format \l\\
\gamma_1=\sum\limits_{j=1}^3 Q_{j,3}\ell_{1,2;j}\quad 
\gamma_2=-\sum\limits_{j=1}^3 Q_{j,2}\ell_{1,3;j}\quad
\gamma_3=\sum\limits_{j=1}^3 Q_{j,2}\ell_{1,2;j}-\sum\limits_{j=1}^3 Q_{j,3}\ell_{1,3;j}\\
\gamma_4=t_{yz}\sum\limits_{j=1}^3 Q_{j,1}\ell_{1,2;j}-t_{y^2}\sum\limits_{j=1}^3 Q_{j,1}\ell_{1,3;j}+\pmb \delta \ell_{2,3;2}\\
\gamma_5=t_{z^2}\sum\limits_{j=1}^3 Q_{j,1}\ell_{1,2;j}-t_{yz}\sum\limits_{j=1}^3 Q_{j,1}\ell_{1,3;j}+\pmb \delta \ell_{2,3;3}\\
\endmatrix$$
are bases for the free $\widetilde{R}_{\pmb\delta}$-modules $(X_{0,2})_{\pmb\delta}$ and $(X_{1,2})_{\pmb\delta}$, respectively.
There is no difficulty in seeing that $[\widehat{\Psi}(g_1),\widehat{\Psi}(g_2),\widehat{\Psi}(g_3),\widehat{\Psi}(g_4),\widehat{\Psi}(g_5)]$ is equal to the matrix $d_1$ from (\tref{d1}). Indeed, for example,
$$\allowdisplaybreaks\alignat1  \widehat{\Psi}(g_1)&{}=\widehat{\Psi}\left (-t_{z^2}\sum\limits_{j=1}^3 Q_{j,1}\ell_{1;j}+\pmb \delta\ell_{3;3}\right)\\&{}
=-t_{z^2}\sum\limits_{j=1}^3 Q_{j,1}\widehat{\Psi}(\kappa({x\otimes x_j}))+\pmb \delta\widehat{\Psi}(\kappa({z\otimes z}))\\&{}
=-t_{z^2}\sum\limits_{j=1}^3 Q_{j,1} \widehat{\Psi}(x x_j)+\pmb \delta \widehat{\Psi}(z^2) 
=-t_{z^2}x (x_j\sum\limits_{j=1}^3 Q_{j,1})   +\pmb \delta  z^2 \\&{}=-t_{z^2}x\lambda_1+\pmb \delta  z^2,\endalignat$$ as expected,  where, as always, we write $x$ for $x_1$, $y$ for $x_2$, and $z$ for $x_3$. 

We next   calculate the matrix for $\pmb \delta\operatorname{Kos}^\Psi\: (X_{1,2})_{\pmb \delta}\to(X_{0,2})_{\pmb \delta}$ with respect to the bases $\{g_1,\dots, g_5\}$ and $\{\gamma_1,\dots,\gamma_5\}$. Use the fact that $\ell_{2;1}=\ell_{1;2}$ to see that the column vector for $$\pmb \delta\operatorname{Kos}^\Psi(\gamma_1)=\pmb \delta\sum\limits_{j=1}^3 Q_{j,3}\operatorname{Kos}^\Psi(\ell_{1,2;j})=\pmb \delta\sum\limits_{j=1}^3 Q_{j,3}(x\ell_{2;j}-y\ell_{1;j}),$$ with respect to the basis for $L_{0,2}$ in (\tref{basis}), is
$$\pmb\delta\bmatrix -yQ_{1,3}\\-yQ_{2,3}+xQ_{1,3}\\-yQ_{3,3}\\xQ_{2,3}\\xQ_{3,3}\\0\endbmatrix.$$
Recall the matrix  $B$ of (\tref{blah}) which expresses the basis of $(X_{0,2})_{\pmb \delta}$ in terms of the basis for $L_{0,2}$.
We claim that
$$\pmb\delta\bmatrix -yQ_{1,3}\\-yQ_{2,3}+xQ_{1,3}\\-yQ_{3,3}\\xQ_{2,3}\\xQ_{3,3}\\0\endbmatrix = B\bmatrix 0\\xQ_{2,3}\\xQ_{3,3}\\xt_{y^2}Q_{1,3}\\xt_{yz}Q_{1,3}-\pmb \delta y\endbmatrix.\tag\tnum{calc}$$
Once (\tref{calc}) is established, then  one reads that 
$$\pmb \delta\operatorname{Kos}^\Psi(\gamma_1)=0 g_1+xQ_{2,3}g_2+xQ_{3,3} g_3+xt_{y^2}Q_{1,3}g_4+(xt_{yz}Q_{1,3}-\pmb \delta y)g_5,$$ as expected. The main trick in the calculation of  (\tref{calc}) involves the coefficient of  $x$ in the top three rows. On the right side, this coefficient is 
$$\bmatrix -t_{z^2}Q_{*,1}&-t_{y^2}Q_{*,1}-t_{yz}Q_{*,1}&Q_{*,2}&Q_{*,3}\endbmatrix \bmatrix 0\\ Q_{2,3}\\ Q_{3,3}\\ t_{y^2}Q_{1,3}\\ t_{yz}Q_{1,3}\endbmatrix$$$$
=Q\bmatrix -t_{y^2}Q_{2,3}-t_{yz}Q_{3,3}\\t_{y^2}Q_{1,3}\\t_{yz}Q_{1,3}\endbmatrix
=Q \bmatrix 0&-t_{y^2}&-t_{yz}\\t_{y^2}&0&0\\t_{yz}&0&0\endbmatrix  Q_{*,3} 
$$$$\eightpoint =Q\left(\bmatrix t_{xy}&0&0\\t_{y^2}&0&0\\t_{yz}&0&0\endbmatrix-\bmatrix t_{xy}&t_{y^2}&t_{yz}\\0&0&0\\0&0&0\endbmatrix\right) Q_{*,3}=\bmatrix 0&0&0\\
\pmb\delta&0&0\\
0&0&0
\endbmatrix Q_{*,3}-Q\bmatrix 0\\0\\0\endbmatrix =\bmatrix 0\\\pmb \delta Q_{1,3}\\0\endbmatrix,$$which is equal to the coefficient of $x$ in the top three rows on the left side of (\tref{calc}). When $\operatorname{Kos}^{\Psi}(\gamma_2)$, $\operatorname{Kos}^{\Psi}(\gamma_3)$, $\operatorname{Kos}^{\Psi}(\gamma_4)$, and $\operatorname{Kos}^{\Psi}(\gamma_5)$ are written in terms of the basis for $L_{0,2}$, then the result is the columns of
$$\hskip-12pt  \eightpoint \left[\smallmatrix 
zQ_{1,2}&-yQ_{1,2}+zQ_{1,3}&-yt_{yz}Q_{1,1}+zt_{y^2}Q_{1,1}&-yt_{z^2}Q_{1,1}+zt_{yz}Q_{1,1}\\
zQ_{2,2}&xQ_{1,2}-yQ_{2,2}+zQ_{2,3}&xt_{yz}Q_{1,1}-yt_{yz}Q_{2,1}+zt_{y^2}Q_{2,1}&xt_{z^2}Q_{1,1}-yt_{z^2}Q_{2,1}+zt_{yz}Q_{2,1}\\
zQ_{3,2}-xQ_{1,2}&\ -yQ_{3,2}-xQ_{1,3}+zQ_{3,3}&\ -yt_{yz}Q_{3,1}-xt_{y^2}Q_{1,1}+zt_{y^2}Q_{3,1}&\ -yt_{z^2}Q_{3,1}-xt_{yz}Q_{1,1}+zt_{yz}Q_{3,1}\\
0&xQ_{2,2}&xt_{yz}Q_{2,1}-z\pmb \delta&xt_{z^2}Q_{2,1}\\
-xQ_{2,2}&0&xt_{yz}Q_{3,1}-xt_{y^2}Q_{2,1}+y\pmb \delta&xt_{z^2}Q_{3,1}-xt_{yz}Q_{2,1}-\pmb \delta z\\
-xQ_{3,2}&-xQ_{3,3}&-xt_{y^2}Q_{3,1}&-xt_{yz}Q_{3,1}+\pmb \delta y
\endsmallmatrix\right].$$ Calculations similar to the one we just made    show that $\pmb \delta$ times the above matrix   is equal to 
$$B\left[\smallmatrix -xQ_{2,3}&-xQ_{3,3}&-xt_{y^2}Q_{1,3}&-xt_{yz}Q_{3,1}+\pmb \delta y\\
0&xQ_{2,2}&xt_{yz}Q_{2,1}-\pmb \delta z&xt_{z^2}Q_{2,1}\\
-xQ_{2,2}&0&xt_{yz}Q_{3,1}-xt_{y^2}Q_{2,1}+\pmb\delta y&\ xt_{z^2}Q_{3,1}-xt_{yz}Q_{2,1}-\pmb \delta z\\
-xt_{yz}Q_{1,2}+z\pmb \delta&\ -xt_{yz}Q_{3,1}+xt_{y^2}Q_{2,1}-\pmb\delta y&0&xQ_{1,1}^2\\
-xt_{z^2}Q_{2,1}&-xt_{z^2}Q_{3,1}+xt_{yz}Q_{2,1}+\pmb \delta z&-xQ_{1,1}^2&0\endsmallmatrix\right];$$ and therefore the matrix $d_2$ from (\tref{d2}) is the  matrix for $\pmb\delta\operatorname{Kos}^{\Psi}\:(X_{1,2})_{\pmb \delta}\to (X_{0,2})_{\pmb \delta}$, with respect to the bases $\{g_1,\dots,g_5\}$ and $\{\gamma_1,\dots,\gamma_5\}$.

We use (\tref{snake}) and (\tref{snake'}) to compute $\frak L_2(\omega)$ for $\omega=x\wedge y\wedge z\in \bigwedge^3_{\Bbb Z} U$. Thus,
$\frak L_2(\omega)= (\operatorname{Kos}^{\Psi}\circ \kappa) (\Theta)$
for $$\Theta= x\otimes \omega \otimes \pmb \delta (\pmb p_1^{\Phi})^{-1}(x^*)+
y\otimes \omega \otimes \pmb \delta (\pmb p_1^{\Phi})^{-1}(y^*)+z\otimes \omega \otimes \pmb \delta (\pmb p_1^{\Phi})^{-1}(z^*)\in \widetilde{R}\otimes_{\Bbb Z} {\tsize\bigwedge}^3_{\Bbb Z} U\otimes_{\Bbb Z} U.$$
The map $\pmb \delta (\pmb p_1^{\Phi})^{-1}\: U^*\to \widetilde{R}\otimes_{\Bbb Z} U$ is given by the classical adjoint $Q$ of $T$: $$\matrix\format \l\\ \pmb \delta (\pmb p_1^{\Phi})^{-1}(x^*)=Q_{1,1}x+Q_{2,1}y+Q_{3,1}z,\\ \pmb \delta (\pmb p_1^{\Phi})^{-1}(y^*)=Q_{1,2}x+Q_{2,2}y+Q_{3,2}z,\text{ and}\\ \pmb \delta (\pmb p_1^{\Phi})^{-1}(z^*)=Q_{1,3}x+Q_{2,3}y+Q_{3,3}z;\endmatrix$$
so, $\Theta$ is equal to $$\cases \phantom{+}x\otimes \omega \otimes   (Q_{1,1}x+Q_{2,1}y+Q_{3,1}z)\\+
y\otimes \omega \otimes   (Q_{1,2}x+Q_{2,2}y+Q_{3,2}z)\\+z\otimes \omega \otimes   (Q_{1,3}x+Q_{2,3}y+Q_{3,3}z)\endcases 
=\cases \phantom{+}(xQ_{1,1}+yQ_{1,2}+zQ_{1,3})\cdot (\omega \otimes   x)\\+
(xQ_{2,1}+yQ_{2,2}+zQ_{2,3})\cdot (\omega \otimes   y)\\+(xQ_{3,1}+yQ_{3,2}+zQ_{3,3})\cdot (\omega \otimes   z)\endcases$$
$$=\lambda_1\cdot (\omega \otimes   x)+
\lambda_2\cdot (\omega \otimes   y)+\lambda_3\cdot (\omega \otimes   z).$$
The composition $\operatorname{Kos}^{\Psi}\circ \kappa$ is equal to $\kappa\circ \operatorname{Kos}^{\Psi}$:
$$\CD  
\widetilde{R}\otimes_{\Bbb Z} L_{2,2} @> \operatorname{Kos}^{\Psi} >> \widetilde{R}\otimes_{\Bbb Z} L_{1,2}\\
@A\kappa AA @A \kappa AA\\
\widetilde{R}\otimes_{\Bbb Z} {\tsize\bigwedge}_{\Bbb Z}^3 U\otimes_{\Bbb Z} U @> \operatorname{Kos}^{\Psi} >> \widetilde{R}\otimes_{\Bbb Z} {\tsize\bigwedge}_{\Bbb Z}^2 U\otimes_{\Bbb Z} U.\endCD$$ It follows that 
$$\frak L_2(\omega)=\lambda_1\cdot (\kappa\circ \operatorname{Kos}^{\Psi})( \omega\otimes x)+  \lambda_2\cdot (\kappa\circ \operatorname{Kos}^{\Psi})(\omega\otimes y)+\lambda_3\cdot (\kappa\circ \operatorname{Kos}^{\Psi})( \omega\otimes z).$$
Recall that $$(\kappa\circ \operatorname{Kos}^{\Psi})( \omega\otimes x)=x\cdot \kappa(y\wedge z\otimes x)-y\cdot \kappa(x\wedge z\otimes x)+
z\cdot\kappa(x\wedge y\otimes x)$$ and that
$$\kappa(y\wedge z\otimes x)=\kappa(x\wedge z\otimes y)-\kappa(x\wedge y\otimes z)=\ell_{1,3;2}-\ell_{1,2;3}.$$Thus, 
$$\frak L_2(\omega)=\left\{\matrix\format \l&\ \l&\ \l\\ 
\phantom{+}x\lambda_1(\ell_{1,3;2}-\ell_{1,2;3})&-y\lambda_1 \ell_{1,3;1}&+z\lambda_1\ell_{1,2;1}\\
+x\lambda_2\ell_{2,3;2}&-y\lambda_2 \ell_{1,3;2}&+z\lambda_2\ell_{1,2;2}\\
+x\lambda_3\ell_{2,3;3}&-y\lambda_3 \ell_{1,3;3}&+z\lambda_3\ell_{1,2;3}.\endmatrix \right.$$
When $\pmb \delta \frak L_2(\omega)$ is written in terms of the basis for $L_{0,2}$, we obtain the matrix on the left side of (\tref{calc'}). The right most factor in (\tref{calc'}) is the matrix we have called $d_3$. The proof is complete as soon as we verify 
$$\pmb \delta\bmatrix \phantom{-}z\bmatrix \lambda_1\\\lambda_2\\\lambda_3\endbmatrix-x\bmatrix 0\\0\\\lambda_1\endbmatrix\\-y\bmatrix \lambda_1\\\lambda_2\\\lambda_3\endbmatrix +x\bmatrix 0\\\lambda_1\\0\endbmatrix \\
\hphantom{-y\bmatrix \lambda_1\\\lambda_2\\\lambda_3\endbmatrix +}x\bmatrix\lambda_2\\\lambda_3\endbmatrix
\endbmatrix =B'\bmatrix \pmb \delta \bmatrix  z^2 \\y^2\\ yz \endbmatrix-x\lambda_1 \bmatrix t_{z^2}\\t_{y^2}\\t_{yz}\endbmatrix\\x\lambda_2\\ x\lambda_3\endbmatrix,\tag\tnum{calc'}$$where $B'$ is the matrix of (\tref{blah'}) which expresses the basis $\gamma_1,\dots,\gamma_5$ of $(X_{1,2})_{\pmb \delta}$ in terms of the basis $\ell_{1,2;1},\dots,\ell_{2,3;3}$ of $L_{1,2}$ as given in (\tref{basis}). 

We verify (\tref{calc'}). The bottom two rows are obvious. The top three rows of the right side of (\tref{calc'}) is equal to 
$$Q\bmatrix x\lambda_2t_{yz}+x\lambda_3t_{z^2}\\\pmb \delta yz-x\lambda_1t_{yz}\\\pmb \delta z^2-x\lambda_1t_{z^2}\endbmatrix=S_1+S_2,\text{ with}$$
$$S_1= xQ\bmatrix 0&t_{yz}&t_{z^2}\\-t_{yz}&0&0\\-t_{z^2}&0&0\endbmatrix Q\bmatrix x\\y\\z\endbmatrix\quad\text{and}\quad S_2=\pmb \delta Q\bmatrix 0\\yz\\z^2\endbmatrix.$$ Observe that $S_1=S_1'+S_1''$ with
$$S_1'=xQ\bmatrix -t_{xz}&0&0\\ -t_{yz}&0&0\\-t_{z^2}&0&0\endbmatrix Q\bmatrix x\\y\\z\endbmatrix=x\bmatrix 0&0&0\\0&0&0\\-\pmb\delta&0&0\endbmatrix Q\bmatrix x\\y\\z\endbmatrix=\bmatrix 0\\0\\-x\pmb \delta \lambda_1\endbmatrix,
\ \text{and}$$
$$S_1''=xQ\bmatrix t_{xz}&t_{yz}&t_{z^2}\\0&0&0\\0&0&0\endbmatrix Q\bmatrix x\\y\\z\endbmatrix=xQ\bmatrix \pmb\delta z\\0\\0\endbmatrix.$$
Thus, $$S_1''+S_2=z\pmb\delta \bmatrix \lambda_1\\\lambda_2\\\lambda_3\endbmatrix$$
and the top three rows of the right side of (\tref{calc'}) is
$$S_1'+(S_1''+S_2)= \bmatrix 0\\0\\-x\pmb \delta \lambda_1\endbmatrix+z\pmb\delta \bmatrix \lambda_1\\\lambda_2\\\lambda_3\endbmatrix,$$ and this is equal to the top three rows of the left side of (\tref{calc'}). Rows 4, 5, and 6 of (\tref{calc'}) are treated in the same manner.  
\qed \enddemo
\endexample

\SectionNumber=\nonempty\tNumber=1
\heading Section \number\SectionNumber. \quad Non-empty disjoint sets of orbits.
\endheading

    We now turn to projects (\tref{EK2}) -- (\tref{EK5}) when $d=3$.  That is, we identify non-empty disjoint subsets of  $\Bbb I_n^{[3]}(\pmb k)$ which are closed under the action of $\operatorname{GL}_{2n+1}\pmb k\times \operatorname{GL}_{3}\pmb k$. 

  Let $\pmb k$ be a field and $P$ be the polynomial ring $P=\pmb k[x,y,z]$. We first record the action of the group $\operatorname{GL}_{2n+1}\pmb k\times \operatorname{GL}_3\pmb k$
on  the sets $\Bbb I_n^{[3]}(\pmb k)$ and $\Bbb X_n(\pmb k)$ from (\tref{p1}). If $A$ is an invertible $(2n+1)\times (2n+1)$ matrix with entries from $\pmb k$ and $\alpha$ is an automorphism of the three-dimensional vector space  $U=[P]_1$,  then $(A,\alpha)\in \operatorname{GL}_{2n+1}\pmb k\times \operatorname{GL}_{3}\pmb k$ carries the matrix $X=(x_{i,j})$ of $\Bbb X_n(\pmb k)$  to the matrix $A(\alpha(x_{i,j}))A^{-1}$ in $\Bbb X_n(\pmb k)$. If $I\in \Bbb I_n^{[3]}(\pmb k)$, then $I$ is generated by the maximal order Pfaffians of some $X$ in $\Bbb X_n(\pmb k)$ and $(A,\alpha)$ carries $I$ to the ideal generated by the maximal order Pfaffians of $(A,\alpha)$ of $X$. Notice that the subgroup  $\operatorname{GL}_{2n+1}\pmb k\times 1$ of  $\operatorname{GL}_{2n+1}\pmb k\times \operatorname{GL}_3\pmb k$ moves the elements of  $\Bbb X_n(\pmb k)$, but acts like the identity on $\Bbb I_n^{[3]}(\pmb k)$. 

Throughout the present section we are interested in $n\ge 3$. Indeed, if $n=2$, then the orbit structure of $\Bbb I_n^{[3]}(\pmb k)$ is not very interesting.
\proclaim{Observation  \tnum{n=2}}  If $\pmb k$ is a field of characteristic not equal to $2$ which is closed under the taking of square root, then $\Bbb I_2^{[3]}(\pmb k)$ consists of exactly one orbit under the action of $1\times \operatorname{GL}_{3}\pmb k$.\endproclaim

\demo{Proof} Let $(\phi)$ be the Macaulay Inverse system for some ideal $I$ in $\Bbb I_2^{[3]}(\pmb k)$. Observe that  the matrix
$$T_{\phi}=\bmatrix \phi(x^2)&\phi(xy)&\phi(xz)\\\phi(xy)&\phi(y^2)&\phi(yz)\\\phi(xz)&\phi(yz)&\phi(z^2)\endbmatrix$$ represents a non-degenerate symmetric bilinear form on the three vector space $U$, whose basis is $x,y,z$. It is well-known, see for example Thm.~XIV.3.1 on page 358 in \cite{\rref{lang}}, and easy to see, that one can choose a new basis for $U$ so that the   matrix for $T_{\phi}$, in the new basis, is the matrix of Example \tref{Ex1}. Thus, under a linear change of variables, $\operatorname{ann}(\phi)$ becomes equal to the ideal $\operatorname{BE}_2$. Actually, the standard theorem from linear algebra converts $T_\phi$ into a diagonal matrix. Our hypothesis about square roots converts  the non-degenerate diagonal matrix  into an identity matrix and $(x,y,z)\mapsto (\frac {x+\sqrt{-1}y}{\sqrt{2}}, \frac {x-\sqrt{-1}y}{\sqrt{2}},\sqrt{-1}z)$ converts the identity matrix to $$\bmatrix 0&1&0\\1&0&0\\0&0&-1\endbmatrix. \qed$$  \enddemo

\bigskip 
For each pair of integers $n$ and $\mu$, with $1\le n$ and $0\le \mu\le 3$, recall the set 
 $$\Bbb I_{n,\mu}^{[3]}(\pmb k)=\left\{I\in \Bbb I_n^{[3]}(\pmb k)\left\vert \matrix\format \l\\ \exists \text{ linearly independent linear forms $\ell_1,\dots,\ell_\mu$ in $P_1$}\\\text{with $\ell_1^n,\dots,\ell_\mu^n$ in $I$ and 
$\not\exists$ $\mu+1$ such forms}\endmatrix\right. \right\}$$of (\tref{p4}).
It is clear that $\Bbb I_n^{[3]}(\pmb k)$ is the disjoint union of $\bigcup\limits_{\mu=0}^3\Bbb I_{n,\mu}^{[3]}(\pmb k)$ and each $\Bbb I_{n,\mu}^{[3]}(\pmb k)$ is closed under the action of $\operatorname{GL}_{2n+1}\pmb k\times \operatorname{GL}_{3}\pmb k$. 

\proclaim {Theorem \tnum{T2}} If $n\ge 3$ and the characteristic of $\pmb k$ is zero, then $\Bbb I_{n,\mu}^{[3]}(\pmb k)$ is non-empty for $0\le \mu\le 3$.\endproclaim

  \remark{Remark} We do {\bf not} claim that every ideal in $\Bbb I_{n,\mu}^{[3]}(\pmb k)$ may be converted into any other ideal in $\Bbb I_{n,\mu}^{[3]}(\pmb k)$ by using $\operatorname{GL}_{2n+1}\pmb k\times \operatorname{GL}_{3}\pmb k$. \endremark

 \demo{Proof of Theorem {\rm\tref{T2}}} In Definition \tref{Inmu} we introduce ideals   $I_{n,\mu}(\pmb k)$ for $\mu$ equal to $0$, $1$, and $2$. The ideals $I_{n,\mu}(\pmb k)$ are shown to be in $\Bbb I_n^{[3]}(\pmb k)$ in Proposition \tref{SEK} and to be in $\Bbb I_{n,\mu}^{[3]}(\pmb k)$ in Proposition \tref{F-1-13}. In Proposition \tref{K} we exhibit an ideal $J_{n,n-1}$, which is in $\Bbb I^{[3]}_{n,3}(\pmb k)$. \qed\enddemo

The ideal $I_{n,2}$ will be defined to be equal to the Buchsbaum-Eisenbud ideal $\operatorname{BE}_n$ of Definition \tref{Hsubn}. To prove Theorem \tref{T2}, we first show that when the hypotheses of Theorem \tref{T2} are in effect, then   $\operatorname{BE}_n$ is in $\Bbb I_{n,2}^{[3]}(\pmb k)$. It is clear that $x^n$ and $y^n$ are in $\operatorname{BE}_n$. One need only show that if $\ell$ is a linear form, then  $\ell^n$ is in $\operatorname{BE}_n$ only if $\ell$ is a multiple of $x$ or $y$. We identify   homogeneous generators for the ideal   $\operatorname{BE}_n$ in Proposition \tref{0.4} and the Macaulay inverse system $(\phi_n)$ for $\operatorname{BE}_n$ in Proposition \tref{phi}.   These calculations work over any field. We complete the proof that $\operatorname{BE}_n$ is in $\Bbb I_{n,2}^{[3]}(\pmb k)$ in Proposition \tref{F-1-13} by solving the equation $\ell^n(\phi_n)=0$; this part of the proof is sensitive to the characteristic of $\pmb k$. In Definition \tref{Inmu} we modify $\phi_n$ twice producing $\phi_{n,\mu}$ for $\mu=0$ and  $\mu=1$; we set  $\phi_{n,2}=\phi_n$; and we define $I_{n,\mu}$ to be $\operatorname{ann}(\phi_{n,\mu})$.   In Proposition \tref{SEK}, we prove that the Gorenstein ideals $I_{n,0}$  and  $I_{n,1}$ are linearly presented by studying how the determinant of $T_{\phi_{n,\mu}}$ is related to $\det T_{\phi_n}$. In   Proposition \tref{F-1-13} we show that $I_{n,\mu}\in \Bbb I_{n,\mu}^{[3]}(\pmb k)$  by solving the equation $\ell^n(\phi_{n,\mu})=0$, as $\ell$ roams over the linear forms of $\pmb k[x,y,z]$.   We do not know how to modify $\phi_n$ to produce an element of $\Bbb I_{n,3}^{[3]}(\pmb k)$; however, we can use ideas from the study of the Weak Lefschetz property to see that the Gorenstein ideal $J_{n,n-1}=(x^n,y^n,z^n):(x+y+z)^{n-1}$ is linearly presented and is generated in degree $n$ (these calculations are very sensitive to the characteristic of $\pmb k$); hence is in $\Bbb I_n^{[3]}(\pmb k)$. It is clear that $x^n$, $y^n$, and $z^n$ all are in $J_{n,n-1}$; hence, $J_{n,n-1}$ is in $\Bbb I_{n,3}^{[3]}(\pmb k)$.

We first study the alternating matrices introduced by Buchsbaum and Eisenbud in Section 6 of \cite{\rref{BE}}. (Our indexing is slightly different than the indexing  of \cite{\rref{BE}}: the matrix that is called $H_{2n+1}$ in \cite{\rref{BE}} is   called $H_n(x,y,z)$ in the present paper.) 
\definition{Definition \tnum{Hsubn}}Let $x,y,z$ be elements of a ring $R$.
For each positive integer $n$, we define the $(2n+1)\times (2n+1)$ alternating matrix $H_{n}(x,y,z)$. The  non-zero entries of $H_{n}(x,y,z)$   above the main diagonal are
$$(H_n(x,y,z))_{i,j}=\cases x&\text{if $i$ is odd and $j=i+1$}\\y&\text{if $i$ is even and $j=i+1$}\\z&\text{if $j=2n+2-i$}.\endcases$$
When $P$ is the ring $\pmb k[x,y,z]$, for some field $\pmb k$, then we let $\operatorname{BE}_n$ denote the ideal of $P$ which is generated by the maximal order Pfaffians of $H_n(x,y,z)$. In other words, $\operatorname{BE}_n$ is generated by $\{B_i\mid 1\le i\le 2n+1\}$ where $B_i$ is the Pfaffian of the $2n\times 2n$ submatrix of $B=H_n(x,y,z)$ which is obtained by deleting row and column $i$. We call the matrix $H_n(x,y,z)$ a Buchsbaum-Eisenbud matrix and the ideal $\operatorname{BE}_n$ a Buchsbaum-Eisenbud ideal. \enddefinition 

\example{Example \tnum{E1.1}}The first few Buchsbaum-Eisenbud  matrices $H_1=H_1(x,y,z)$ and $H_2=H_2(x,y,z)$ are given in (\tref{Hn}).
When $n=1$, then $B_1=y$, $B_2=z$, and $B_3=x$. When $n=2$, then $B_1=y^2$, $B_2=xz$,  $B_3=xy+z^2$, $B_4=yz$, and $B_5=x^2$. 
\endexample

In Proposition \tref{0.4} we explicitly identify the maximal order Pfaffians of the Buchsbaum Eisenbud matrices. Our first step is to express the Pfaffians of one of these matrices  in terms of the Pfaffians of smaller  matrices of the same form. Recall our Pfaffian conventions from Subsection \number\prelim.\number\PC.
   
 \proclaim{Lemma \tnum{BE-gens}}  Let $x$, $y$, and $z$ be elements of a ring $P$ and $n\ge 3$ be an integer. If $B=H_n(x,y,z)$, $b=H_{n-1}(y,x,z)$, and $\tilde b=H_{n-2}(x,y,z)$, then 
$$
\matrix
\format \r &\ \c\ &\l&\c&\r &\ \c\ &\l\\ 
     B_1&=&yb_{2n-1},&\quad       &B_2     &=&zb_1, \quad B_i= xy\tilde b_{i-2}+z b_{i-1} \text{ for $ 3 \leq i \leq 2n-1$}, \\
  B_{2n}&=&zb_{2n-1},&\ \text{ and }\ &B_{2n+1}&=&xb_1. 
\endmatrix
$$
 \endproclaim
 
 \demo{Proof} Throughout this proof we use the fact that $b$ is the submatrix of $B$ obtained by deleting the first and last rows and columns. Expand the Pfaffians along the last column to obtain
   $$B_1 = yB_{1, 2n, 2n+1}=yb_{2n-1}\text{ and }B_{2n}=zB_{1, 2n, 2n+1}= z b_{2n-1}.$$
 Expand  the rest of the Pfaffians along the first row: 
  $B_2 = zB_{1, 2, 2n+1}=zb_1$,  $B_{2n+1}= xB_{1, 2,  2n+1}= xb_1$,  and for all $i$ with $3 \leq i \leq 2n-1$,
$$\split
B_i &{}=xB_{1, 2, i} + zB_{1, i, 2n+1} 
= xyB_{1, 2, i, 2n, 2n+1}+zB_{1, i, 2n+1}\\
&{}= xyb_{1, i-1, 2n-1}+zb_{i-1}= xy \tilde b _{i-2} +zb_{i-1}. \qed
\endsplit$$
\enddemo

\proclaim{Proposition \tnum{0.4}} If $x$, $y$, and $z$ are elements of a ring $P$ and $n$ is a positive integer, then
the ideal generated by the maximal order Pfaffians of the matrix $H_n(x,y,z)$ is generated by
$$\{ x^is_{n-i}\mid 1\le i\le n\}\cup \{s_n\}\cup\{ y^is_{n-i}\mid 1\le i\le n\},$$ where
 $s_i= \sum\limits_{j=0}^{\lfloor i/2 \rfloor} \binom{i- j}   jx^jy^jz^{i-2j}$.
\endproclaim

\example{Example \tnum{E1.2}} The first few   $s_i$ are
$s_0= 1$, $s_1=z$, $s_2=z^2+xy$, and $s_3=z^3+2xyz$.\endexample
 
\demo{Proof of Proposition {\rm\tref{0.4}}} Let $B=H_n(x,y,z)$. We prove the result by showing that    the maximal order Pfaffians of $B$ are given by 
$$\nopagebreak\eightpoint B_i= \cases  
x^{n+1-i}s_{i-1} &\text{if $i$ is even}\\
y^{n+1-i}s_{i-1} & \text{if $i$ is odd} 
\endcases   
\hskip0.2cm \text{and} \hskip0.2cm
B_{2n+2-i}= \cases
y^{n+1-i}s_{i-1} & \text{if $i$ is  even}\\
x^{n+1-i}s_{i-1} & \text{if $i$ is odd}, 
\endcases\tag\tnum{Cl2}$$for $ 1\leq i \leq n+1$. We establish (\tref{Cl2}) by induction on $n$. If $n$ is $1$ or $2$, then Examples \tref{E1.1} and \tref{E1.2} show that (\tref{Cl2}) holds. Henceforth, we assume $3\le n$ and we apply Lemma  \tref{BE-gens} with $b=H_{n-1}(y,x,z)$ and $\tilde b=H_{n-2}(x,y,z)$. Induction gives
 $$\split B_1&{}=yb_{2n-1}= y\cdot y^{n-1},\ \ B_2 = zb_1= z \cdot x^{n-1},\ \  B_{2n}= z b_{2n-1} =  z\cdot y^{n-1},\  \text{and}\\  B_{2n+1}&{}= xb_1=x\cdot x^{n-1}.\endsplit $$   
Now suppose $ 3 \leq i \leq n+1$. If $i$ is odd, then 
$$\split
B_i&= xy \tilde b_{i-2}+z b_{i-1} \\
&= xy (y^{n+1-i}s_{i-3})+z(y^{n-i+1}s_{i-2}) \qquad \text{by induction} \\
&=   y^{n+1-i}(xys_{i-3}+zs_{i-2})
\endsplit$$ and  $ 
B_{2n+2-i} = x^{n+1-i}(xys_{i-3}+zs_{i-2}).$
If $i$ is even, then $$B_i=  x^{n+1-i}(xys_{i-3}+zs_{i-2})\text{ and } 
B_{2n+2-i}= y^{n+1-i}(xys_{i-3}+zs_{i-2}).$$ We complete the proof by showing that
$$xys_{\alpha-1}+zs_{\alpha} = s_{\alpha+1}\tag\tnum{sts}$$ for all $\alpha$ with $1 \leq \alpha \leq n-1$. Indeed, we see  that 
$xys_{\alpha -1}+ zs_{\alpha}$ is equal to $$\allowdisplaybreaks\align  
&\tsize {} \phantom{{}={}}xy  \sum\limits_{j=0}^{\lfloor\frac{\alpha-1}{2}\rfloor} \binom{ \alpha-1- j} j  x^jy^jz^{\alpha-1-2j} + z \sum\limits_{j=0}^{\lfloor\frac\alpha 2\rfloor} \binom{ \alpha- j} j  x^jy^jz^{\alpha-2j}\\
&\tsize {}=
  \sum\limits_{j=0}^{\lfloor\frac{\alpha-1}{2}\rfloor} \binom{ \alpha-1- j} jx^{j+1}y^{j+1}z^{\alpha-1-2j} +  \sum\limits_{j=0}^{\lfloor\frac\alpha 2\rfloor} \binom{ \alpha- j } jx^jy^jz^{\alpha+1-2j}\\
&\tsize {}=  \sum\limits_{j=1}^{\lfloor\frac{\alpha+1}{2}\rfloor} \binom{\alpha- j }{j-1}x^{j}y^{j}z^{\alpha+1-2j} +  \sum\limits_{j=1}^{\lfloor\frac\alpha 2\rfloor} \binom {\alpha- j }j x^jy^jz^{\alpha+1-2j}+z^{\alpha +1}.\endalign$$
Notice that $\lfloor\frac{\alpha+1}{2}\rfloor=\lfloor\frac{\alpha}{2}\rfloor+\chi$, where
$$\chi=\cases 0&\text{if $\alpha$ is even}\\1&\text{if $\alpha$ is odd;}\endcases$$and therefore, 
$xys_{\alpha -1}+ zs_{\alpha}$ is equal to $$\allowdisplaybreaks\align  
&\tsize {} \phantom{{}={}} \chi x^{\lfloor\frac{\alpha+1}{2}\rfloor}y^{\lfloor\frac{\alpha+1}{2}\rfloor} + \sum\limits_{j=1}^{\lfloor\frac{\alpha}{2}\rfloor} \binom{\alpha- j }{j-1}x^{j}y^{j}z^{\alpha+1-2j} +  \sum\limits_{j=1}^{\lfloor\frac\alpha 2\rfloor} \binom {\alpha- j }j x^jy^jz^{\alpha+1-2j}+z^{\alpha +1}\\
&\tsize {}=   \chi x^{\lfloor\frac{\alpha+1}{2}\rfloor}y^{\lfloor\frac{\alpha+1}{2}\rfloor} + \sum\limits_{j=1}^{\lfloor\frac{\alpha}{2}\rfloor} \left[\binom{\alpha- j }{j-1}  +   \binom {\alpha- j }j\right] x^jy^jz^{\alpha+1-2j}+z^{\alpha +1}
\\
&\tsize {}=  \chi x^{\lfloor\frac{\alpha+1}{2}\rfloor}y^{\lfloor\frac{\alpha+1}{2}\rfloor} + \sum\limits_{j=1}^{\lfloor\frac{\alpha}{2}\rfloor}   \binom{\alpha +1- j }{j}   x^jy^jz^{\alpha+1-2j}+z^{\alpha +1}\\
&\tsize {}=   \sum\limits_{j=0}^{\lfloor\frac{\alpha+1}{2}\rfloor}   \binom{\alpha +1- j }{j}   x^jy^jz^{\alpha+1-2j}=s_{\alpha+1}. 
\endalign$$We have established (\tref{sts}); therefore the proof is complete. \qed
\enddemo

\proclaim{Proposition \tnum{phi}} Let $n$ be a positive integer, $\pmb k$ be a field, $U$ be a vector space of dimension three over $\pmb k$ with basis $x,y,z$, $P$ be the polynomial ring $P=\operatorname {Sym}^{\pmb k}_\bullet(U)=\pmb k[x,y,z]$, and $\operatorname{BE}_n$ be the Buchsbaum-Eisenbud ideal. Then the Macaulay inverse  system for $\operatorname{BE}_n$ is the $P$-submodule of $D_{\bullet}^{\pmb k}(U^*)$ which is generated by  $$\phi_n=   \sum\limits_{i=0}^{n-1}(-1)^ic_i{x^*}^{(n-1-i)}{y^*}^{(n-1-i)}{z^*}^{(2i)}\in D_{2n-2}^{\pmb k}(U^*),$$ where $c_i$ is the $i^{th}$ Catalan number
$c_i = \frac{1}{i+1} {2i \choose i}$. 
\endproclaim  

\remark{Note} The first few $\phi$'s are $$\phi_1=1,\quad \phi_2=x^*y^*-z^{*(2)},\quad \text{and}\quad \phi_3=x^{*(2)}y^{*(2)}-x^{*}y^{*}z^{*(2)}+2z^{*(4)}.$$ \endremark

\demo{Proof} The ideal $\operatorname{BE}_n$ is presented by the matrix $H_{n}(x,y,z)$, which has homogeneous linear entries. It follows that the socle degree of $\operatorname{BE}_n$ is $2n-2$; and therefore Macaulay's Theorem guarantees that $\operatorname{ann} \operatorname{BE}_n$ is a homogeneous cyclic $P$-submodule of $D_{\bullet}(U^*)$ generated by an element of $D_{2n-2}(U^*)$. As a consequence, any non-zero element of degree $2n-2$ in $\operatorname{ann} \operatorname{BE}_n$ is a generator of $\operatorname{ann} \operatorname{BE}_n$. To prove the result, it suffices to show that $\phi_n\in \operatorname{ann} \operatorname{BE}_n$. In light of   Proposition \tref{0.4}, it suffices to show that $x^{n-i}s_i(\phi_n)$ and $y^{n-i}s_i(\phi_n)$ are zero for $0\le i\le n$. The expressions $s_i$ and $\phi_n$ are both symmetric in $x$ and $y$; consequently, it suffices to show that $x^{n-i}s_i(\phi_n)=0$ for $0\le i\le n$.
We   compute  
$$\split x^{n-i}(\phi_n)&=x^{n-i} \left(  \sum_{k=0}^{n-1}(-1)^kc_k{x^*}^{(n-1-k)}{y^*}^{(n-1-k)}{z^*}^{(2k)} \right)\\
&\\
&=   \sum_{k=0}^{i-1}(-1)^kc_k{x^*}^{(i-1-k)}{y^*}^{(n-1-k)}{z^*}^{(2k)}.
\endsplit$$
It follows that 
$$  s_i(x^{n-i}\phi_n)= \sum{i-j \choose j} c_{k}(-1)^{k}{x^*}^{(i-1-k-j)}{y^*}^{(n-1-k-j)}{z^*}^{(2k-i+2j)},$$where the sum 
is taken over all pairs $(k,j)$ which satisfy:
$$\tsize 0\le j\le \lfloor \frac i2\rfloor,\quad 0\le k\le i-1,\quad 0\le i-1-k-j,\quad\text{and}\quad 0\le 2k-i+2j.$$ Replace $k$ with $i-1-\ell-j$ to obtain $$ \tsize  s_i(x^{n-i}\phi_n)= \sum\limits_{\ell=0}^{\lfloor \frac{i-2}{2} \rfloor}(-1)^{i-1-\ell} \left[ \sum\limits_{j=0}^{\lfloor \frac{i}{2} \rfloor}{i-j \choose j} (-1)^jc_{i-1-\ell-j}\right]{x^*}^{(\ell)}{y^*}^{(\ell+n-i)}{z^*}^{(i-2-2\ell)}. $$ The sum inside the brackets is zero due to Bennett's identity \cite{\rref{BG}}:
$$\sum_{j=0}^{\lfloor\frac m2\rfloor}(-1)^jc_{m-n-j}\binom {m-j}j =0\tag\tnum{BI}$$ when $m$ and $n$ are positive integers with $2n\le m$. The proof of (\tref{BI}) that is given in \cite{\rref{BG}} is based on generating functions.
The power series expansion for 
$F(z) = (1 -\sqrt{1 - 4z})/(2z)$ is $\sum_{i=0}^\infty c_i z^i$ and the coefficient $x^mz^{m-n}$ in the power series expansion of $\frac{F(z)}{1-(x-zx^2)}$ is equal to the left side of (\tref{BI}) .
\qed
\enddemo

\definition{Definition \tnum{Inmu}} Let $n$ be a positive integer, $\pmb k$ be a field, $U$ be a vector space of dimension three over $\pmb k$ with basis $x,y,z$, $P$ be the polynomial ring $P=\operatorname {Sym}^{\pmb k}_\bullet(U)=\pmb k[x,y,z]$, and $\mu$ be one of the integers $0$, $1$, or $2$. Define $\phi_{n,u}$ is the element $$\phi_{n,\mu}= \displaystyle \sum_{i=0}^{n-1}(-1)^ic_i{x^*}^{(n-1-i)}{y^*}^{(n-1-i)}{z^*}^{(2i)}+\chi(\mu\le 1){x^*}^{(2n-2)}+2\chi(\mu=0){y^*}^{(2n-2)} $$ of $D_{2n-2}^{\pmb k}U^*$  and define $I_{n,\mu}$ to be the ideal $I_{n,\mu}=\operatorname{ann}(\phi_{n,\mu})$ of $P$. \enddefinition

\remark{Remark} The symbol $\chi$ is defined  in (\tref{chi}). In particular,
$$\phi_{n,2}=\phi_n,\qquad \phi_{n,1}=\phi_n+{x^*}^{(2n-2)},\qquad\text{and}\qquad \phi_{n,0}=\phi_n+{x^*}^{(2n-2)}+2{y^*}^{(2n-2)},$$ for $\phi_n$ as defined in Proposition \tref{phi}. It is a consequence of Proposition \tref{phi}, that $\operatorname{BE}_n=I_{n,2}$.  \endremark

\proclaim{Proposition \tnum{SEK}} If $\pmb k$ is a field, $n$ is a positive integer, and $\mu$ is equal to $0$, $1$, or $2$, then the ideal $I_{n,\mu}$ of Definition {\rm\tref{Inmu}} is in the set $\Bbb I_{n}^{[3]}(\pmb k)$ of {\rm(\tref{p1})}.
\endproclaim

\demo{Proof}   Recall, from Definition \tref{T}, that  
$T_{\phi_{n,\mu}}$ is the $ N \times N$ matrix $(\phi_{n,\mu}( m_im_j))$, where $N = {n+1 \choose 2}$ and $\left\{m_i \right\}$ is a basis for $\operatorname {Sym}^{\pmb k}_{n-1}U$. Let $m_1=x^{n-1}$ and $m_2=y^{n-1}$.   Observe that the matrix $T_{\phi_{n,\mu}}$      has the form  $$T_{\phi_{n,\mu}}=\left( \matrix M_\mu& 0\\ 0&M' \endmatrix \right),
$$ where   $M_\mu$  is the  $2\times 2$ matrix with entries $(\phi_n( m_im_j))$ and $(\phi_{n,\mu}( m_im_j))$, respectively,   with $1\le i, j\le 2$. 
A quick calculation yields that 
$$
M_2=\left( \smallmatrix 0&1\\ 1&0 \endsmallmatrix \right),\   
M_1=\left( \smallmatrix 1&1\\ 1&0 \endsmallmatrix \right),\ \text{and}
\ M_0=\left( \smallmatrix 1&1\\ 1&2\endsmallmatrix \right).$$
The construction of the ideal $\operatorname{BE}_n=I_{n,2}$ puts this ideal in $\Bbb I_n^{[3]}(\pmb k)$; and therefore, according to Proposition \tref{J18}, $\det T_{\phi_{n,2}}\neq 0$; hence, $\det M'\neq 0$. It is clear that both matrices $M_{0}$ and $M_1$ have non-zero determinant and therefore $\det T_{\phi_{n,\mu}}\neq 0$, for $\mu$ equal to $0$ and $1$. It follows from Proposition \tref{J18} that $I_{n,0}$ and $I_{n,1}$ are both in $\Bbb I_n^{[3]}(\pmb k)$.
\qed
\enddemo

\proclaim{Proposition \tnum{F-1-13}} Let $\pmb k$ be a field of characteristic zero,  $n\ge 3$ be a positive integer,
  $\mu$ equal   $0$, $1$, or $2$, and $I_{n,\mu}\subseteq P$ be the ideal of Definition {\rm\tref{Inmu}}.
If  ${\ell=\alpha x+\beta y+ \gamma z}$ is an arbitrary linear form in $P$, with $\alpha$, $\beta$, and $\gamma$ in $\pmb k$, then 
$$\ell^n\in I_{n,\mu} \iff \cases \alpha=\gamma=0\quad \text{or}\quad \beta=\gamma=0&\text{when $\mu=2$}\\
\alpha=\gamma=0 &\text{when $\mu=1$}\\
\alpha=\beta=\gamma=0&\text{when $\mu=0$}.\endcases\tag\tnum{iff}$$In particular, the ideal $I_{n,\mu}$ is in the set $\Bbb I_{n,\mu}^{[3]}(\pmb k)$ of {\rm(\tref{p4})}.\endproclaim

\example{Example} The hypothesis $n\ge 3$ in Proposition \tref{F-1-13} is necessary. Indeed, if $\pmb k$ is the field $\Bbb Q[\sqrt{2}]$, then $x^2$, $y^2$, and $(x+y+\sqrt{2} z)^2$ all are in the ideal $\operatorname{BE}_2$. \endexample

\demo{Proof of Proposition {\rm\tref{F-1-13}}}We proved in Proposition \tref{SEK} that $I_{n,\mu}$ is in $\Bbb I_{n}^{[3]}(\pmb k)$, which is the disjoint union   $\bigcup\limits_{i=0}^3\Bbb I_{n,i}^{[3]}(\pmb k)$. To prove that $I_{n,\mu}$ is in $\Bbb I_{n,\mu}^{[3]}(\pmb k)$, it suffices to establish (\tref{iff}).  The direction ($\Leftarrow$) of (\tref{iff}) is obvious. We prove ($\Rightarrow$). Fix $\ell$ with $\ell^n\in I_{n,\mu}$. 

We first assume that $4\le n$. One calculates 
$$\ell^n=\displaystyle \sum_{a+b+c =n} {n \choose a,b,c} (\alpha x)^a(\beta y)^b (\gamma z)^c\quad\text{and}$$
$$\ell^n(\phi_{n,\mu})= \cases\displaystyle \sum_{a+b+c =n}\displaystyle \sum_{i=0}^{n-1}(-1)^ic_i {n \choose a,b,c} \alpha^a\beta^b\gamma^c{x^*}^{(n-1-i-a)} {y^*}^{(n-1-i-b)} {z^*}^{(2i-c)}\\+\chi(\mu\le 1)\alpha^n{x^*}^{(n-2)}+2\chi(\mu=0)\beta^n{y^*}^{(n-2)}.\endcases$$
Recall that the binomial coefficient $\binom n{a,b,c}$ is zero if any of the parameters $a$, $b$, or $c$ is negative. 
Each coefficient of $\ell^n(\phi_{n,\mu})$ is zero; so, in particular: \vphantom{\tnum{2.}\tnum{5.}\tnum{1.}
\tnum{3.}\tnum{7.}}
$$\eightpoint \matrix\format\l&\ \c&\quad\l\\ &\text{the coefficient of}&\text{in $\ell^n(\phi_{n,\mu})$ is}\\
(\tref{2.})&{x^*}^{(n-2)}&0=\beta^{n-2}\left(n \alpha\beta- {n \choose 2}  \gamma^2\right)+\chi(\mu\le 1)\alpha^n\\\vspace{3pt}
(\tref{5.})&{y^*}^{(n-2)}&0=\alpha^{n-2}\left( n \alpha \beta -  {n \choose 2}  \gamma^2 \right) +2\chi(\mu=0)\beta^n\\
(\tref{1.})& {z^*}^{(n-2)}&0=\sum\limits_{i=0}^{n-1}(-1)^ic_i {n \choose n-1-i,n-1-i,2i+2-n} \alpha^{n-1-i}\beta^{n-1-i}\gamma^{2i+2-n}\\&&\phantom{0}{}=(-1)^{n-1}c_{n-1}\gamma^n+\alpha\beta \kappa, \text{for some integer $\kappa$ in $\pmb k$} \\\vspace{3pt}
(\tref{3.})& {x^*}^{(n-3)}z^*&0=\beta^{n-3}\gamma
\left( 2{n \choose 3}   \gamma^2
 -   {n \choose 1,n-2,1} \alpha \beta\right)\\\vspace{3pt}
(\tref{7.})&x^*{y^*}^{(n-3)}&0=\alpha^{n-4}\left(2 {n \choose 4} \gamma^4 -  {n \choose n-3,1,2} \alpha\beta\gamma^2+  {n \choose 2} \alpha^{2}\beta^2\right)\\\vspace{3pt}
\endmatrix$$The hypothesis $4\le n$ ensures that the $5$ listed elements of $D_{n-2}^{\pmb k}U^*$ are distinct. 

We first show that $$\alpha\beta\gamma=0.\tag\tnum{abc}$$
 Indeed, if 
all three constants are non-zero, then 
we may combine (\tref{3.}) and (\tref{7.}) to see that
$(\alpha\beta,\gamma^2)$ is a point in the intersection 
$$\tsize 0=  2{n \choose 3}   Y
 -   {n \choose 1,n-2,1} X \ \text{and}\ 0=2 {n \choose 4} Y^2 -  {n \choose n-3,1,2} XY+  {n \choose 2} X^2.$$However the only intersection point is $(0,0)$. Thus, in every case, at least one of the constants $\alpha$, $\beta$, or $\gamma$ must be zero and (\tref{abc}) is established.

Next we show that $$\gamma=0.\tag\tnum{gma}$$ Indeed, (\tref{abc}) ensures that at least one of the constants are zero. Furthermore, we may 
apply (\tref{1.}) to see that if $\alpha=0$ or $\beta=0$, then $\gamma$ is also zero. We conclude that (\tref{gma}) holds.

Now that (\tref{gma}) holds, we apply (\tref{7.}) again to see that  
$$\alpha\beta=0.\tag\tnum{ab}$$ The proof is complete if $\mu=2$.

We now focus on $\mu=1$. We have shown that $0=\gamma=\alpha\beta$. If $\beta=0$, then apply (\tref{2.}) to conclude $\alpha=0$. Thus $\alpha$ must be zero if $\mu=1$ and the proof is complete in this case. 

Finally, we assume that $\mu=0$. We have shown that (\tref{gma}) and (\tref{ab}) hold. If $\alpha=0$, then (\tref{5.}) yields that $\beta$ is also zero. If $\beta=0$, then (\tref{2.}) yields that $\alpha$ is also zero. Thus all three constants are zero and the proof is complete in this case.

Now we treat the case $n=3$. The argument is similar. 
Each coefficient of $\ell^3(\phi_{3,\mu})$ is zero; so, in particular:  
$$\matrix\format\l&\ \c&\quad\l\\ &\text{the coefficient of}&\text{in $\ell^3(\phi_{3,\mu})$ is}\\
(\tref{2.}')&{x^*}^{(n-2)}&0= 3 \alpha\beta^2 -  3 \beta\gamma^2+\chi(\mu\le 1)\alpha^3\\\vspace{3pt}
(\tref{5.}')&{y^*}^{(n-2)}&0=3 \alpha^2\beta     
-  3\alpha\gamma^2 +2\chi(\mu=0)\beta^3\\
(\tref{1.}')& {z^*}^{(n-2)}&0=-  6 \alpha\beta\gamma +
   2   \gamma^3. \endmatrix$$
As before, we first establish (\tref{abc}). If $\mu$ is $1$ or $2$, then an easy argument yields that every simultaneous solution of (\tref{5.}$'$) and (\tref{1.}$'$) is also a solution of (\tref{abc}). If $\mu=0$, then one can show that 
every simultaneous solution of (\tref{2.}$'$), (\tref{5.}$'$) and (\tref{1.}$'$) is also a solution of (\tref{abc}). One now uses (\tref{1.}$'$) to show that (\tref{gma}) holds. 

If $1\le \mu$, then (\tref{ab}) follows from (\tref{5.}$'$), which  now is $0=3\alpha^2\beta$ since $\gamma=0$. If $\mu=0$, then one can use (\tref{2.}$'$) and (\tref{5.}$'$), which now are $0=3\alpha\beta^2+\alpha^3$ and $0=3 \alpha^2\beta+2\beta^3$ to conclude (\tref{ab}). Thus, (\tref{ab}) holds in all cases. The proof is complete when $\mu=2$.

To complete the proof when $\mu=1$, we use (\tref{2.}$'$), together with (\tref{gma}) and (\tref{ab}), to see that $\alpha=\gamma=0$. To complete the proof when $\mu=0$, we use (\tref{2.}$'$) and  (\tref{5.}$'$)  together with (\tref{gma}) and (\tref{ab}), to see that $\alpha=\beta=\gamma=0$.
\qed\enddemo

We could not modify the  generator of the Macaulay inverse system $\phi_n$ of the Buchsbaum Eisenbud ideal $\operatorname{BE}_n$ to produce an element of $\Bbb I_{n,3}^{[3]}(\pmb k)$; however the  ideal  $$J_{n,n-1}=(x^n,y^n,z^n):(x+y+z)^{n-1},$$ which  arises in the study of the Weak Lefschetz Property, (see, for example, Observation \tref{hype}) is  in  $\Bbb I_{n,3}^{[3]}(\pmb k)$ when $\pmb k$ has characteristic zero.

\proclaim{Proposition \tnum{K}} Let $n$ be a positive integer, $\pmb k$ be a field, $P$ be the polynomial ring $\pmb k[x,y,z]$, and $J_{n,n-1}$ be the  ideal $(x^n,y^n,z^n):(x+y+z)^{n-1}$ of $P$. If the characteristic of $\pmb k$ is  zero, then $J_{n,n-1}$ is in $\Bbb I_{n,3}^{[3]}(\pmb k)$. \endproclaim

\demo{Proof} It is clear that $J_{n,n-1}$ contains $x^n$, $y^n$, and $z^n$. We apply Proposition \tref{J18} to show that $J_{n,n-1}$ is in $\Bbb I_n^{[3]}(\pmb k)$. It suffices to show that $[P]_{2n-1}\subseteq J_{n,n-1}$ and $[J_{n,n-1}]_{n-1}$ is equal to $0$. It is clear that $[P]_{3n-2}\subseteq (x^n,y^n,z^n)$; and therefore $[P]_{2n-1}\subseteq J_{n,n-1}$. Furthermore,  Theorem 5 in \cite{\rref{RRR}} guarantees that the minimal generator degree of $\frac{J_{n,n-1}}{(x^n,y^n,z^n)}$ is at least $n$; and 
therefore, $[J_{n,n-1}]_{n-1}=0$. \qed \enddemo

\example{Example} A quick calculation shows that $$J_{1,0}=(x,y,z)\quad \text{and} \quad J_{2,1}=(x^2,y^2,z^2,z(x-y),y(x-z)).$$\endexample

\bigskip\noindent{\bf Acknowledgment.} The authors are grateful to L\'aszl\'o Sz\'ekely for his suggestions concerning the Catalan numbers and to Maria Evelina Rossi for conversations concerning divided powers.

\enddocument